\documentclass[12pt,twoside,english,a4paper]{article}
\pdfoutput=1
\usepackage{etex}
\usepackage{caption}
\usepackage[latin1]{inputenc}
\usepackage[intlimits] {amsmath}   
\usepackage{graphicx}
\usepackage{psfrag}
\usepackage{ifthen}
\usepackage{fancyhdr}
\usepackage{rotating}
\usepackage{multirow}
\usepackage{booktabs}              
\usepackage{amssymb}
\usepackage{amsbsy}
\usepackage{bm}                    
\usepackage{bbm}
\usepackage{babel}
\usepackage{theorem}
\usepackage{upgreek}  
\usepackage{pifont}   
\usepackage[active]{srcltx}   
\usepackage{subcaption}
\usepackage{suetterl}  
\newfont{\suetdbl}{suet14 scaled 2000}  
\usepackage{comment}
\usepackage{yfonts}  
\newfont{\gothdbl}{ygoth scaled 2000} 
\newfont{\frakdbl}{yfrak scaled 2000}
\newfont{\swabdbl}{yswab scaled 2000}
\usepackage{caption}
\usepackage{subcaption}
\usepackage{placeins} 
\usepackage[normalem]{ulem}  
\usepackage[textsize=footnotesize,textwidth=1.7cm]{todonotes}
\usepackage[]{defjs3}      
\usepackage{stmaryrd}
\usepackage{subcaption}
\usepackage[]{defjs3}  
\usepackage{amsfonts}
\usepackage{graphicx}
\usepackage{cite}

\usepackage[authoryear]{natbib}
\bibpunct{[}{]}{,}{n}{}{;}
\fancypagestyle{plain}{%
  \fancyhead{}%
  \fancyfoot[c]{\sffamily\thepage}%
}
\makeatletter                           
\def\cleardoublepage{\clearpage\if@twoside \ifodd\c@page\else
  \hbox{}
  \vspace*{\fill}
  \thispagestyle{empty}
  \newpage
  \if@twocolumn\hbox{}\newpage\fi\fi\fi}
\makeatother
\def\Curl{\operatorname{Curl}}      
\def\curl{\operatorname{curl}}       
\newcommand{\Ce}{\mathbb{C}_{\mathrm{e}}}
\newcommand{\Cc}{\mathbb{C}_{\mathrm{c}}}
\newcommand{\Cmicro}{\mathbb{C}_{\mathrm{micro}}}
\newcommand{\Cmacro}{\mathbb{C}_{\mathrm{macro}}}
\newcommand{\Lc}{L_{\mathrm{c}}}
\newcommand{\IL}{\mathbb{L}}
\newcommand{\R}{\mathbb{R}}
\newcommand{\Bdis}{\bP}
\newcommand{\dis}{P}
\begin{document}
\unitlength1.0cm
\frenchspacing


\thispagestyle{empty}
\ce{\bf \large
Lagrange and $H(\curl,\B)$ based Finite Element formulations for} 
\ce{\bf \large the relaxed micromorphic model}

\vspace{4mm}
\ce{ J\"org Schr\"oder$^{1,*}$, Mohammad Sarhil$^{1}$, Lisa Scheunemann$^2$ and  Patrizio Neff$^3$}

\vspace{4mm}
\ce{$^1$Institute of Mechanics, University of Duisburg-Essen}
\ce{Universit\"atsstr. 15, 45141 Essen, Germany}
\ce{\small e-mail: 
j.schroeder@uni-due.de
}    

\vspace{4mm}
\ce{$^2$Chair of Applied Mechanics, TU Kaiserslautern,}
\ce{67653 Kaiserslautern, Germany}

\vspace{4mm}
\ce{$^3$Faculty of Mathematics, University of Duisburg-Essen,}
\ce{ 45141 Essen, Germany}

\vspace{4mm}
\begin{center}
{\bf \large Abstract}
\bigskip

{\footnotesize
\begin{minipage}{14.5cm}
\noindent
Modeling the unusual mechanical properties of  metamaterials is a challenging topic for the mechanics community and enriched continuum theories are promising computational tools for such materials. The so-called relaxed micromorphic model has shown many advantages in this field. In this contribution, we present the significant aspects related to the relaxed micromorphic model realization with the finite element method. The  variational problem is derived and different FEM-formulations for the two-dimensional case are presented. These are a nodal standard formulation $H^1(\B) \times H^1(\B)$ and a nodal-edge formulation $H^1(\B) \times H(\curl, \B)$, where the latter employs the N\'ed\'elec space. However, the implementation of higher-order N\'ed\'elec elements is not trivial and requires some technicalities which are demonstrated. We discuss the convergence behavior of Lagrange-type and tangential-conforming finite element discretizations. Moreover, we analyze the characteristic length effect on the different components of the model and reveal how the size-effect property is captured via this characteristic length. 
\end{minipage}
}
\end{center}


\sect{\hspace{-5mm}. Introduction}
\label{sec:intor}
 Metamaterials are receiving tremendous attention in both academia and industry due to their unconventional mechanical properties. These are not solely governed by the bulk mechanical properties but also by the geometry of the unit cells which can be designed to attain the desired functionality, see \cite{FisHilEbe:2020:mmo,LeeSinTho:2012:mnm,SurGaoDuLiXioFanLu:2019:AEM,YaZhoLiaJiaWu:2018:mma,Zad:2016:mmm}. Moreover, the recent advances of  the additive manufacturing (AM, or 3D printing) techniques are empowering the fabrication process of three-dimensional architected  metamaterials, c.f. \cite{JiaLi:2018:3pa,MonKuaArnQi:2020:rai,PloPan:2019:rod,LeiHonZhaHamCheLuQi:2019:3po}. To simplify the design process of novel metamaterials, suitable computational tools are needed to capture their unprecedented effective mechanical properties. The classical Cauchy-Boltzmann theory and the first-order homogenization procedures often fail to describe the mechanical macroscopic behavior of mechanical metamaterials since they exhibit the size-effect phenomenon, i.e. small specimens are stiffer than big specimens, and therefore other generalized theories are needed such as the classical Mindlin-Eringen micromorphic theory \cite{Min:1964:msi,SuhEri:1964:nto,EriSub:1964:nto,Eri:1968:mom,LeiMah:2015:coh,LeiMah:2015:tfh}, the Cosserat theory \cite{CosCos:1909:tof,Nef:2006:tcc,NefJeoMueRam:2010:lce}, gradient elasticity  \cite{MinEsh:1968:ofsg,AltAif:1997:osa,FisKlaMerSteMue:2011:iao} or others. 

The relaxed micromorphic model, which we adopt in this work, has been introduced recently in \cite{NefGhiMadPlaRos:2014:aup,GhiNefMadPlaRos:2015:trl,NefGhiLazMad:2015:trl}. It keeps the full kinematics of the  micromorphic theory but employs the matrix Curl operator of a non-symmetric second-order micro-distortion field for the curvature measurement.  The relaxed micromorphic model reduces the complexity of the classical micromorphic theory by decreasing the number of material parameters and has shown many advantages such as the separation of the material parameters into scale-dependent and scale-independent ones, see for example \cite{DagBarGhiEidNefMad:2020:edo}.  Furthermore, it has already been used to obtain the main mechanical characteristics (stiffness, anisotropy, dispersion) of the targeted metamaterials for many well-posed dynamical and statical problems, e.g. \cite{MadNefGhiPlaRos:2015:bgi,MadNefGhiPlaRos:2015:wpi,MadNefGhiRos:2016:rat,MadNefBar:2016:cbg,NefMadBarDagAbrGhi:2017:rwp,BarMadDagAbrGhiNef:2017:taf,MadNefBarGhi:2017:aro,MadColMinBilOuiNef:2018:mpc,BarTalDagAivNefMad:2019:rmm,AivTalAgoDaoNefMad:2020:fan}. Recently, the scale-independent short range elastic parameters in the relaxed micromorphic model were determined for artificial periodic microstructures in \cite{NefEidMad:2019:ios} which are used to capture the band-gaps as  a dynamical property of  mechanical metamaterial in \cite{DagBarGhiEidNefMad:2020:edo}.  Analytical solutions of the relaxed micromorphic compared to the solutions of other generalized continua for some essential boundary value problems, i.e. pure shear, bending, torsion and uniaxial tension, are discussed in \cite{RizHueKhaGhiMadNef:2021:aso1,RizKhaGhiMadNef:2021:aso2,RizHueMadNef:2021:aso3,RizHueMadNef:2021:aso4}, emphasizing the validity of the relaxed micromorphic model for small sizes  (bounded stiffness), where most of the other generalized continua exhibit unphysical stiffness properties. 

As a result of employing the matrix Curl operator of the micro-distortion field for the curvature measurement, the relaxed micromorphic model seeks the solution of the micro-distortion in $H(\curl,\B)$, while the displacement solution is still in $H^1(\B)$.  The appropriate finite elements of such case must be conforming in $H(\curl,\B)$  (tangentially conforming). The first formulation of edge elements was presented in \cite{RavTho:1977:amf}. In fact, the name ``edge'' elements was used because the degree of freedoms (dofs) are associated only with edges for the first-order approximation. $H(\curl,\B)$-conforming finite elements of  first kind were introduced in \cite{Ned:1980:mfe} and second kind in \cite{Ned:1986:anf}, which are comparable with $H(\div,\B)$-conforming  elements of first kind in \cite{RavTho:1977:amf} and second kind in \cite{BreDouMar:1984:tfo}. An extension to elements with curved edges, based on covariant projections, was developed by \cite{CroSilHur:1988:cpe}. A general implementation of  N\'ed\'elec elements of  first kind is presented in \cite{OlmBadMAr:2019:oag} and a detailed review about $H(\div,\B)$- and $H(\curl,\B)$-conforming finite elements is available in \cite{KirLogRogTer:2012:cau} and \cite{RogKirAnd:2009:eao}.  Furthermore, hierarchical $H(\curl)$-conforming finite elements are used to solve Maxwell boundary and eigenvalue problems in \cite{schzag:2005:hon}. A  $H^1(\B) \times H(\curl,\B)$ finite element formulation for a simplified anti-plane shear case of the relaxed micromorphic model utilizing a scalar displacement field and a vectorial micro-distortion field is available in \cite{SkyNeuMueSchNef:2021:CM}.

 In our work, we demonstrate the main technologies related to the finite element realization of the theoretically-sound relaxed micromorphic model. The proper finite element approximation of the micro-distortion field is the N\'ed\'elec space which utilizes tangential-conforming vectorial shape functions. We provide a comprehensive description of the construction of  $H^1(\B) \times H(\curl, \B)$ elements with N\'ed\'elec  formulation of  first kind on triangular and quadrilateral meshes.  Six finite elements are built, which are different in the approximation space of the micro-distortion: two triangular elements with first- and second-order N\'ed\'elec formulation, two quadrilateral elements with first- and second-order N\'ed\'elec formulation, and two nodal triangular elements with standard  first- and second-order Lagrangian formulation. The paper is organized as follows. In Section \ref{sec:model}, we introduce the relaxed micromorphic model and derive the variational problem with the resulting strong forms and  the associated boundary conditions which are modulated in a physical point of view by the so-called consistent coupling condition. We cover in Section \ref{sec:fem} the main components of the implementation of standard nodal and nodal-edge elements. Two numerical example are introduced in Section \ref{sec:fem}. The first numerical example is designed to check the convergence behavior of the different finite elements when the solution is discontinuous in the micro-distortion field. We investigate the influence of the characteristic length in a second example which covers the size-effect property. We conclude the paper in Section \ref{sec:con}. 

\sect{\hspace{-5mm}. The relaxed micromorphic model}
\label{sec:model} 

The relaxed micromorphic model is a continuum model which describes the kinematics of a material point using a displacement vector $\bu\colon\B\subseteq\R^3\to\R^3$ and a non-symmetric micro-distortion field  $\Bdis\colon\B\subseteq\R^3\to\R^{3\times3}$. Both are defined  for the static case by the minimization of potential
\begin{equation}
\label{eq:pot}
\Pi (\bu,\Bdis) = \int_\B W\left(\nabla \bu,\Bdis,\Curl \Bdis \right) \, - \overline{\bbf}\cdot{\bu} \, - \overline{\bM} : \Bdis \,\textrm{d}V\ - \int_{\partial \B_t} \overline{\bt} \cdot \bu  \,\textrm{d}A  \longrightarrow\ \min\,,
\end{equation}
with $(\bu,\Bdis)\in H^1(\B)\times H(\curl,\B)$. The vector $\overline\bbf$ and the tensor $\overline{\bM}$ describe, respectively, the given body force and the body moment, while $\overline\bt$ is the traction vector acting on the boundary $\partial \B_t \subset \partial \B$. The elastic energy density $W$ reads  
\begin{equation}
\begin{aligned}
W\left(\nabla \bu,\Bdis,\Curl \Bdis \right) =& \, \dfrac{1}{2} ( \symb{ \nabla \bu - \Bdis} : \Ce : \symb{ \nabla \bu - \Bdis} +   \sym \Bdis : \Cmicro: \sym \Bdis  \,, \\
&+ \skewb{ \nabla \bu - \Bdis} : \Cc : \skewb{ \nabla \bu - \Bdis} +   \mu \, \Lc^2 \, \textrm{Curl} \Bdis : \IL :\textrm{Curl} \Bdis ) \,. \\
\end{aligned}
\end{equation}
Here, $\Cmicro,\Ce$ are fourth-order positive definite standard elasticity tensors, $\Cc$ is a fourth-order positive semi-definite rotational coupling tensor, $\IL$ is a fourth-order tensor, $\Lc$ is a non-negative parameter describing the characteristic length scale and $\mu$ is a typical shear modulus.  The characteristic length parameter plays a significant role in the relaxed micromorphic model. This parameter is related to the size of the microstructure  and determines its influence on the macroscopic mechanical behavior. A relation of the relaxed micromorphic model to the classical Cauchy model was shown in \cite{NefEidMad:2019:ios} for limiting values of $\Lc$, which we can also observe in our numerical examples.

The variation of the potential with respect to the displacement field, i.e. $\delta_{\bu} \Pi = 0$, with 
\begin{equation}
\delta_{\bu} \Pi =  \int_{\B} \{ \underbrace{\Ce : \symb{ \nabla \bu - \Bdis}+  \Cc : \skewb{ \nabla \bu - \Bdis}}_{ \textstyle =:\Bsigma} \} :  {\nabla \delta \bu} - \overline\bbf \cdot \delta \bu  \, \textrm{d}V - \int_{\partial \B_t} \overline\bt \cdot \delta \bu  \, \textrm{d}A\,,
\end{equation} 
leads after integration by parts and applying the divergence theorem to the weak form 
 \begin{equation}
\delta_{\bu} \Pi =  \int_{\B}   \{ \div \Bsigma + \overline\bbf \} \cdot \delta \bu \, \textrm{d}V  = 0 \, ,
\end{equation}  
where  $\Bsigma$ denotes the non-symmetric force stress tensor.
The associated strong form with the related boundary conditions reads 
 \begin{equation}
   \label{eq:sf1}
 \div\Bsigma+\overline\bbf = \bzero \,  \quad \textrm{with} \quad \bu = \overline{\bu} \quad \textrm{on} \quad \partial \B_u\, \quad  \textrm{and} \quad \overline\bt = \Bsigma \cdot \bn \quad \textrm{on} \quad  \partial \B_t \,, 
\end{equation}  
satisfying  $\partial \B_u \cap \partial \B_t = \emptyset $ and  $\partial \B_u \cup \partial \B_t = \partial \B $ and $\bn$ is the outward normal on the boundary. In a similar way, the variation of the potential with respect to the micro-distortion field, i.e. $\delta_{\Bdis} \Pi = 0$, with 
\begin{equation}
\label{eq:de1}
 \delta_{\Bdis} \Pi =   \int_{\B} \{ \Bsigma  - \underbrace{\Cmicro : \sym \Bdis}_{\textstyle =:\Bsigma_\textrm{micro} } + \overline\bM\} :  \delta {\Bdis}  -  \underbrace{\mu \, \Lc^2  (\IL : \Curl \Bdis)}_{\textstyle =:\bbm} : \Curl \delta {\Bdis} \, \textrm{d} V\,, 
\end{equation}
yields after integration by parts and applying Stokes' theorem  
\begin{equation}
\begin{aligned}
 \delta_{\Bdis} \Pi = &  \int_{\B}  \{ \Bsigma - \Bsigma_\textrm{micro}  - \Curl \bbm + \overline\bM\} :  \delta {\Bdis} \, \textrm{d} V +  \int_{ \partial \B}  \{ \sum_{i=1}^3  \left( \bbm^i  \times \delta {\Bdis}^i  \right) \cdot \bn  \}  \;\; \textrm{d}A  
 = 0\,, 
\end{aligned} 
\end{equation}
where $\Bsigma_\textrm{micro}  $ and $\bbm$ are  the micro- and moment stresses, respectively, and $\bbm^i$ and $\delta \Bdis^i$ denote the row vectors of the associated second-order tensors. Using the identity of the scalar triple vector product
\begin{equation}
(\ba \times \bb) \cdot \bc = (\bc \times \ba) \cdot \bb  = (\bb \times \bc) \cdot \ba \,,
\end{equation}
allows for the reformulation
\begin{equation}
 \int_{ \partial \B}  \{ \sum_{i=1}^3  \left( \bbm^i  \times \delta {\Bdis}^i  \right) \cdot \bn  \}  \;\; \textrm{d}A =
 \int_{ \partial \B_P}  \{ \sum_{i=1}^3  \left(   \delta {\Bdis}^i  \times \bn  \right) \cdot \bbm^i  \}  \;\; \textrm{d}A   -  \int_{ \partial \B_m}  \{ \sum_{i=1}^3  \left(  \bbm^i  \times \bn  \right) \cdot \delta {\Bdis}^i   \}  \;\; \textrm{d}A  \,.
\end{equation}
The associated strong form reads  
 \begin{equation}
   \label{eq:sf2}
 \Curl \bbm = \Bsigma - \Bsigma_\textrm{micro} + \overline\bM \,,
 \end{equation}
 with the related boundary conditions  
 \begin{equation}
\sum_{i=1}^3  {\Bdis}^i  \times \bn  = \overline{\bt}_p \quad \textrm{on} \quad \partial \B_\dis \quad \textrm{and by definition} \quad \sum_{i=1}^3    \bbm^i  \times \bn = \bzero \quad \textrm{on} \quad  \partial \B_m\,, 
 \end{equation}
 where $\partial \B_\dis \cap \partial \B_m = \emptyset $ and $\partial  \B_\dis \cup \partial \B_m= \partial \B $.  
A dependency between the displacement field and the micro-distortion field on the boundary was proposed by \cite{NefEidMad:2019:ios} and subsequently considered in \cite{SkyNeuMueSchNef:2021:CM,RizHueMadNef:2021:aso3,RizHueMadNef:2021:aso4,DagRizKhaLewMadNef:2021:tcc}. This so-called consistent coupling condition is defined by
 \begin{equation}
  \Bdis \cdot \Btau = \nabla \bu \cdot \Btau \, \Leftrightarrow \,  \Bdis \times \bn = \nabla \bu \times \bn  \quad \textrm{on}\quad  \partial \B_\dis = \partial \B_u  \,,
 \end{equation}
where $\Btau$ is a tangential vector on the Dirichlet boundary.  This condition relates the projection of the displacement gradient on the tangential plane of the boundary to the respective parts of the micro-distortion.

The first strong form in Equation (\ref{eq:sf1}) represents a generalized balance of linear momentum (force balance) while the second strong form in Equation (\ref{eq:sf2}) outlines a generalized balance of angular momentum (moment balance). The generalized moment balance invokes the Cosserat theory with the $\Curl \Curl$ operator rising from the matrix $\Curl$ operator of the second-order moment stress $\bbm$.  In comparison to the classical micromorphic model, see \cite{Eri:1968:mom,NefGhiMadPlaRos:2014:aup}, the relaxed micromorphic model uses the same kinematical measures but employs a curvature measure form the Cosserat theory, see \cite{NefJeoMueRam:2010:lce}. The micro-distortion field has the following general form for the three-dimensional case
 \begin{equation}
\Bdis = \left[ \begin{array}{c}
(\Bdis^{1})^T \\
(\Bdis^{2})^T \\
(\Bdis^{3})^T \\
\end{array}\right] = \left[\begin{array}{c c c}
 \dis_{11}  &  \dis_{12}  & \dis_{13}  \\   
 \dis_{21}  &  \dis_{22}  & \dis_{23}  \\ 
 \dis_{31}  &  \dis_{32}  & \dis_{33}  \\ 
\end{array}\right] \quad \textrm{with} \quad \Bdis^{i} = \left[ \begin{array}{c}
\dis_{i1} \\
\dis_{i2} \\ 
\dis_{i3}
\end{array} \right],
 \end{equation}
where $\Bdis^{i}$  are the row vectors of $\Bdis$. 
We let the Curl operator act on the row vectors of the micro-distortion field $\Bdis$, i.e.,
 \begin{equation}
\Curl \Bdis = \left[ \begin{array}{c}
(\curl \Bdis^{1})^T \\
(\curl \Bdis^{2})^T \\
(\curl \Bdis^{3})^T \\
\end{array}\right] = \left[\begin{array}{c|c|c}
 \dis_{13,2} - \dis_{12,3}	& \dis_{11,3} - \dis_{13,1}  &  \dis_{12,1} - \dis_{11,2} \\
 \dis_{23,2} - \dis_{22,3}	& \dis_{21,3} - \dis_{23,1}  &  \dis_{22,1} - \dis_{21,2} \\
 \dis_{33,2} - \dis_{32,3}	& \dis_{31,3} - \dis_{33,1}  &  \dis_{32,1} - \dis_{31,2} 
\end{array}\right].
 \end{equation}
For the two-dimensional case, the micro-distortion field and its Curl operator are reduced to 
 \begin{equation}
\label{eq:t:2dCurlBdis}
\Bdis = \left[ \begin{array}{c}
(\Bdis^{1})^T \\
(\Bdis^{2})^T \\
\bzero^T 
\end{array}\right] = \left[\begin{array}{c c c }
 \dis_{11}  &  \dis_{12}   & 0\\   
 \dis_{21}  &  \dis_{22}   & 0\\
 0 & 0 & 0 
\end{array}\right] \quad \textrm{and} \quad 
\Curl \Bdis = 
\left[\begin{array}{c|c|c}
0	&   0  &  \dis_{12,1} - \dis_{11,2} \\
0   &   0  &  \dis_{22,1} - \dis_{21,2} \\
0 & 0 & 0
\end{array}\right].
 \end{equation}
\sect{\hspace{-5mm}. Approximation spaces}
\label{sec:fem} 
\ssect{\hspace{-5mm}. Nodal elements $(\bu,\Bdis) \in H^1(\B) \times H^1(\B)$}
\label{sec:standardelements}
We introduce the formulation of a  standard nodal element utilizing Lagrange-type shape functions for both displacement and micro-distortion field, see for example \cite{Wri:2008:nfe}. Let us assume that there are $n^u$ nodes in each element for the discretization of the displacement field $\bu$ and $n^\dis$ nodes for micro-distortion field $\Bdis$ in two dimensions. Geometry and displacement field are approximated employing the related Lagrangian shape functions $N^u_I$ defined in the parameter space with the natural coordinates $\Bxi = \{ \xi, \eta\}$ by
\begin{equation}
\label{eq:t:approx_chi_u}
 \bX_h = \sum_{I=1}^{n^u} N^u_I(\Bxi)  \bX_I\,, \qquad \bu_h = \sum_{I=1}^{n^u} N^u_I(\Bxi)  \bd^u_I\,,
\end{equation}
where $\bX_I$  are the coordinates of the displacement node  $I$ and  $\bd^u_I$ are its displacement degrees of freedom. The deformation gradient is obtained in the physical space by 
\begin{equation}
\nabla \bu_h = \sum_{I=1}^{n^u} \bd^u_I \otimes \nabla N^u_I(\Bxi)\, \quad \textrm{with} \quad 
\nabla {N^u_I(\Bxi)} = \bJ^{-T} \cdot \nabla_\Bxi N^u_{I}\,,
\end{equation}
where $\bJ = \frac{\partial \bX}{\partial \Bxi  }$ is the Jacobian, $\nabla$  and $\nabla_\Bxi $ denote the gradient operators with respect to $\bX$ and $\Bxi$, respectively.  
The micro-distortion field $\Bdis$ for the 2D case is approximated using the relevant scalar shape functions $N^\dis_I$  
\begin{equation}
\Bdis_h^1 =  \left[ \begin{array}{c}
\dis_{11} \\
\dis_{12} \\
\end{array} \right] =   \sum_{I=1}^{n^\dis}   N^\dis_I(\Bxi) \bd^{\dis^1}_I\,, \quad
\Bdis_h^2 =  \left[ \begin{array}{c}
\dis_{21} \\
\dis_{22} \\
\end{array} \right] =   \sum_{I=1}^{n^\dis}   N^\dis_I(\Bxi) \bd^{\dis^2}_I\,,
\end{equation}
where  $\bd^{\dis^1}_I$ and $\bd^{\dis^2}_I$ are the micro-distortion row vectors degrees of freedom of node $I$. In order to calculate the Curl of  $\Bdis$, the gradient of the row vectors in the physical space can be calculated by
\begin{equation}
\nabla \Bdis_h^{i} = \bJ^{-T} \cdot \nabla_\Bxi \Bdis_h^{i} \quad \textrm{for} \quad i=1,2 
\end{equation}
and the rotation of the vector $\Bdis_h^{i}$ reads 
\begin{equation}
\curl^{2D}{\Bdis_h^{i}} = (\dis_h)_{i2,1} - (\dis_h)_{i1,2} \quad \textrm{for} \quad i=1,2 \,.
\end{equation}

\begin{figure}[ht]
\center
	\unitlength=1mm
	\begin{picture}(120,50)
	\put(0,5){\def\svgwidth{11cm}{\small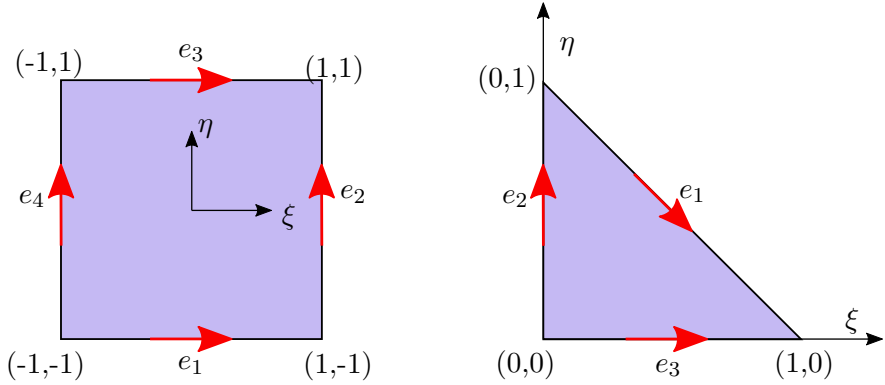}}
	\end{picture}
	\caption{Lowest-order ($k=1$) N\'ed\'elec elements: triangle $\left[ \mathcal{ND}^\triangle \right]^{2}_1$ (right) and quadrilateral $\left[ \mathcal{ND}^\square \right]^{2}_1$ (left). Definition of the individual edges $e_i$. The red arrows indicate the orientation of the tangential flux. }
	\label{Figure:nedelec_elements}
\end{figure}

\ssect{\hspace{-5mm}. Nodal-edge elements $(\bu,\Bdis) \in H^1(\B) \times H(\curl, \B)$ \\} 
\label{sec:mixelements}
The here presented formulation uses different spaces to describe the micro-distortion field. The geometry and the displacement field are approximated in the standard Lagrange space as in Equation (\ref{eq:t:approx_chi_u}).  For the micro-distortion field, its solution is in $ H(\curl,\B)$ and the suitable finite element space is known as N\'ed\'elec space, see \cite{Ned:1980:mfe,Ned:1986:anf}. In this work, we choose the N\'ed\'elec space of  first kind. For more details the reader is referred to  \cite{KirLogRogTer:2012:cau,RogKirAnd:2009:eao,BofBreFor:2014:mfe,Mon:1993:ano}. N\'ed\'elec formulations use vectorial shape functions which satisfy the tangential continuity at element interfaces. The lowest-order two-dimensional N\'ed\'elec elements are depicted in Figure  \ref{Figure:nedelec_elements}. 

\begin{figure}[ht]
	\unitlength=1mm
	\center
	\unitlength=1mm
	\center		  
		  	  \begin{subfigure}[b]{0.31\textwidth}
	\begin{picture}(50,45)
	\put(0,0){\def\svgwidth{5cm}{\small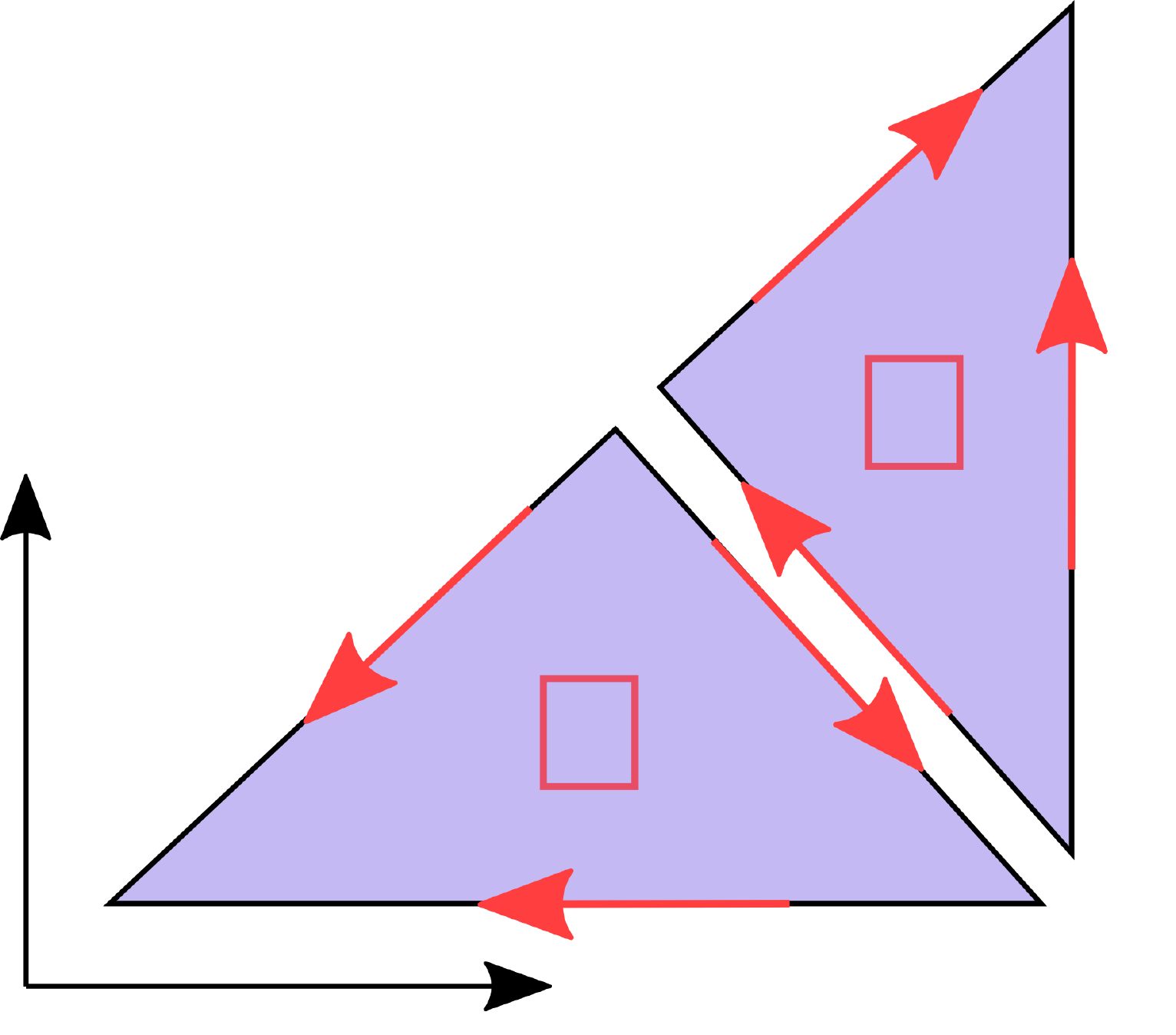}}
	\end{picture}
	\caption{local orientations of dofs}
		  \end{subfigure}
		  	  \begin{subfigure}[b]{0.35\textwidth}
\begin{tabular}{ c|c c c}
element & $\alpha_1$ & $\alpha_2$ & $\alpha_3$ \\ \hline
1 & -1 & +1 & -1 \\ 
2 & +1 & -1 & +1 
\end{tabular}
	\caption{orientation parameter values}
		  \end{subfigure}		  
 \begin{subfigure}[b]{0.31\textwidth}
	\begin{picture}(50,50)
	\put(0,0){\def\svgwidth{5cm}{\small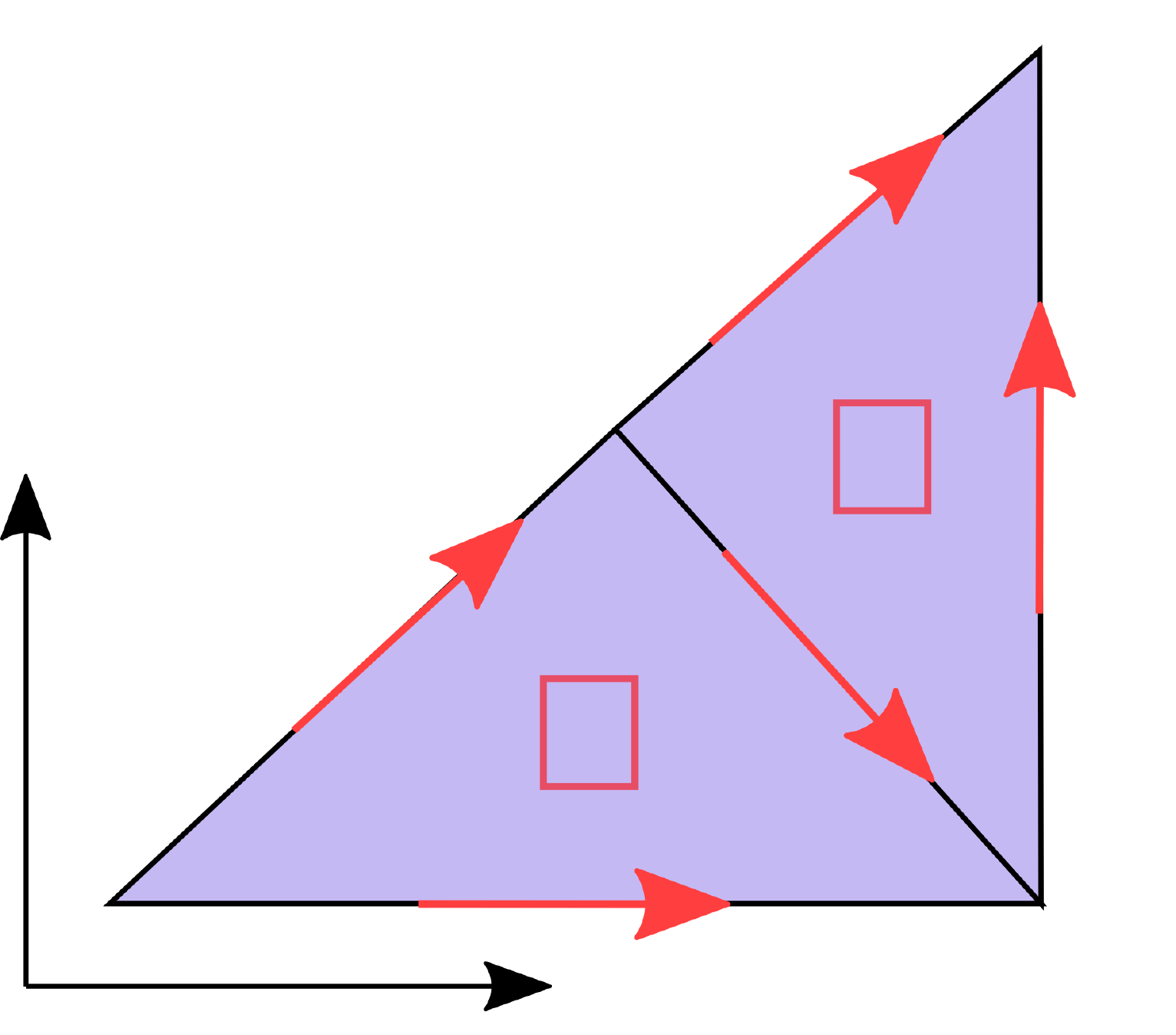}}
	\end{picture}
	\caption{global orientations of dofs}
		  \end{subfigure} 
		  	\caption{Example of assembling of two neighboring elements with satisfying the orientation consistency via the orientation parameter $\alpha_I$.  }
\label{fig:orientation_parameter}
\end{figure}

 Triangular  N\'ed\'elec elements  of order $k$ are based on the space 
\begin{equation}
\label{eq:t:shape_functions_form_t}
\left[ \mathcal{ND}^\triangle \right]^{2}_k = (\IP_{k-1})^2 \oplus \bS_k \quad \textrm{with} 
 \quad \bS_k = \{ \bp   \in (\tilde{\IP}_k)^2 \,|\, \bp \cdot \Bxi = 0 \}\,,
\end{equation}
 where  $\IP_{k-1}$ is the linear space of polynomials of degree $k-1$ or less and $\tilde{\IP}_{k}$ is the linear space of homogeneous polynomials of degree $k$. Equivalently, the space can be characterized by 
 \begin{equation}
\left[ \mathcal{ND}^\triangle \right]^{2}_k = (\IP_{k-1})^2 \oplus \tilde{\IP}_{k-1} \left[ \begin{array}{c}
- \eta  \\
\xi
\end{array} \right]. 
\end{equation}
  The dimension of this linear space is $
k(k+2) $.  Quadrilateral N\'ed\'elec elements of order $k$ are based on the linear space 
 \begin{equation}
 \label{eq:t:shape_functions_form_q}
\left[ \mathcal{ND}^\square \right]^{2}_k = \left[ \begin{array}{c}
 Q_{k-1,k}  \\
 Q_{k,k-1}
\end{array} \right] \quad \textrm{where} \quad Q_{m,n} = \textrm{span}\{ \xi^i  \eta^j  \,|\, i \leq m,  j \leq n \}, 
 \end{equation}   
with $
  \dim \left( \left[ \mathcal{ND}^\square \right]^{2}_k \right) = 2k(k+1) 
$.
The vectorial shape functions $\bv^k$ in the parametric space are obtained by constructing a linear system of equations based on a set of inner and outer dofs. For the 2D case, the outer dofs of an edge $e_i$ are defined by the integral 
   \begin{equation}
   \label{eq:t:edge_dof}
 m^{e_i}_{j} (\bv^k ) = \int_{e_i} (\bv^k \cdot \bt_i)  \, r_{j} \, \textrm{d} s \,, \quad \forall \, r_j \in \IP_{k-1}(e_i)\,
   \end{equation}  
   where $r_{j}$ is a polynomial $\IP_{k-1}$  along the edge $e_i$ and $\bt_i$ is the normalized tangential vector of the edge $e_i$. The inner dofs are introduced for triangular elements by 
\begin{equation}   
   \label{eq:t:inner_dof_t}
 m^\textrm{inner}_{i}  (\bv^k )   =  \int_{\B_e} \bv^k \cdot \bq_{i} \, \textrm{d} a \,, \quad \forall \, \bq_i \in (\IP_{k-2} (\B_e))^2\, ,
   \end{equation}
   while they are given for quadrilateral elements by
   \begin{equation}   
     \label{eq:t:inner_dof_q}
 m^\textrm{inner}_{i}  (\bv^k )    =  \int_{\B_e} \bv^k   \cdot \bq_{i} \, \textrm{d} a \,, \quad \forall \, \bq_i  \, 
   \in
    \left[ \begin{array}{c}
 Q_{k-1,k-2}  (\B_e) \\
 Q_{k-2,k-1}  (\B_e)
\end{array} \right].
   \end{equation}
 The scalar-valued and vectorial functions $r_j$ and $\bq_i$ are linearly independent polynomials which are chosen as Lagrange polynomials in our work. For lowest-order element ($k=1$), only outer dofs occur. For higher-order elements ($k \ge 2$), the number of outer dofs are increased and additional inner dofs are introduced. E.g. for the $\left[ \mathcal{ND}^\triangle \right]^{2}_2$ with a dimension 8, we have 6 outer dofs and 2 inner ones. The derivations of the $H(\curl,\B)$-conforming vectorial shape functions is shown in Appendix \ref{app:Nedelec_shape_functions}.  Mapping the vectorial shape functions $ \bv^k_I$ from the parametric space to $\hat\Bpsi^k_I $ in the physical space must conserve the tangential continuity property. This is guaranteed by using the covariant Piola transformation, see for example \cite{RogKirAnd:2009:eao}, which reads
\begin{equation}
\hat{\Bpsi}^k_I =  \, \bJ^{-T} \cdot \bv^k_I\, \quad \textrm{and} \quad \curl{\hat{\Bpsi}^k_I} = \frac{1 }{\det{\bJ}} \bJ \cdot \curl_{\Bxi}\bv^k_I\,. 
\end{equation}
For our implementation of the $H(\curl,\B)$-conforming elements, we modify the mapping to enforce the required orientation of the degrees of freedom at the inter-element boundaries and to attach a direct physical interpretation to the Neumann-type boundary conditions. Hence, two additional parameters, $\alpha$ and $\beta$, appear for the vectorial shape functions associated with edge dofs 
\begin{equation}
\label{eq:t:mpiola}
{{\Bpsi}}^k_I  = \alpha_I \beta_I {\hat\Bpsi}^k_I   \quad \textrm{and} \quad \curl{{\Bpsi}^k_I} = \alpha_I \beta_I  \curl{{\hat\Bpsi}^k_I} \,, 
\end{equation}

where $ \alpha_I = \pm 1$ is the orientation consistency function  which ensures that on an edge, belonging to two neighboring finite elements, a positive tangential flux direction is defined. Therefore, a positive tangential direction is defined based on a positive $x$-coordinate. A tangential component pointing in negative $x$-direction is multiplied by a value $\alpha_I = -1$ to obtain the overall positive tangential flux on each edge. If the tangential direction has no projection on $x$-axis, then the same procedure is employed on $y$-direction. Figure \ref{fig:orientation_parameter} illustrates an example of calculating the orientation parameter values of two neighboring elements. 

 The normalization parameter $\beta_I$  enforces that the sum of the vectorial shape functions ${\Bpsi}^k_I $ at a common edge  scalar multiplied with the associated tangential vector has to be equal one in the physical space. Furthermore,  the sum of the shape functions belonging to one edge scalar multiplied with the tangential vector of the other edges must vanish. These conditions are reflected by 
\begin{equation}
\label{eq:t:cond}
\Btau_I \cdot \sum_J {{\Bpsi}}_J^k  \bigg\rvert_{E_I}  \equiv 1  \quad \, \textrm{if}  \quad  I=J  \quad \textrm{and} \, \quad
\Btau_I \cdot \sum_J {{\Bpsi}}_J^k  \bigg\rvert_{E_I}  \equiv 0 \quad \textrm{if} \quad  I \ne J  \,.
\end{equation}
Here, $\sum_J {{\Bpsi}}_J^k  \bigg\rvert_{E_I} $ is the sum of shape vectors related to outer dofs of an edge $E_J$ evaluated on the edge $E_I$ and $\Btau_I$ is the normalized tangential vector of an edge $E_I$ where $E$ denotes the edges in the physical space.  
Based on Equations (\ref{eq:t:mpiola})$_1$ and (\ref{eq:t:cond})$_1$, we compute straightforward the parameters $\beta_I$. In detail we get for the first- and second-order elements 
\begin{equation}
\beta_I = L_I   \quad \textrm{and} \quad   \beta_I = \frac{L_I}{2}\,,
\end{equation}
respectively, where $L_I$ denotes the length of the edge $E_I$ in the physical space. For the 2D case, the rotation of the vectorial shape functions only has one active component out of the plane which reads 
\begin{equation}
 \curl^{2D}{\Bpsi^k_I} = \frac{ \alpha_I  \beta_I }{\det{\bJ}} \curl^{2D}_{\Bxi}\bv^k_I\,.
\end{equation}
The micro-distortion field $\Bdis$ is  approximated by the vectorial dofs $\bd^\dis_I$ representing its tangential components at the location $I=1,...,n^\dis$. The micro-distortion field and its Curl are interpolated as 
\begin{equation}
\label{eq:t:dis_restoring}
\Bdis_h =   \sum_{I=1}^{n^\dis}   \bd^\dis_I \otimes \Bpsi^k_I \, , \quad 
\Curl \Bdis_h  =  \sum_{I=1}^{n^\dis}   \bd^\dis_I  \otimes \curl{{\Bpsi}^k_I} \,.
\end{equation}
The non-vanishing components of the Curl operator of the micro-distortion field for the 2D case are obtained by 
\begin{equation}
\left[\begin{array}{c}
\curl^{2D}{\Bdis_h^{1}} \\
  \curl^{2D}{\Bdis_h^{2}} \\
\end{array}\right]\,
 =  \sum_{I=1}^{n^\dis}   \bd^\dis_I  \curl^{2D}{{\Bpsi}^k_I} =  \left[\begin{array}{c}
\sum_{I=1}^{n^\dis}   (d^\dis_I)_1  \curl^{2D}{{\Bpsi}^k_I} \\
 \sum_{I=1}^{n^\dis}   (d^\dis_I)_2  \curl^{2D}{{\Bpsi}^k_I}  \\
\end{array}\right]. 
\end{equation}

\begin{figure}[ht]
	\unitlength=1mm
	\center
		  	  \begin{subfigure}[b]{0.32\textwidth}
	\begin{picture}(40,40)
	\put(0,0){\def\svgwidth{5cm}{\small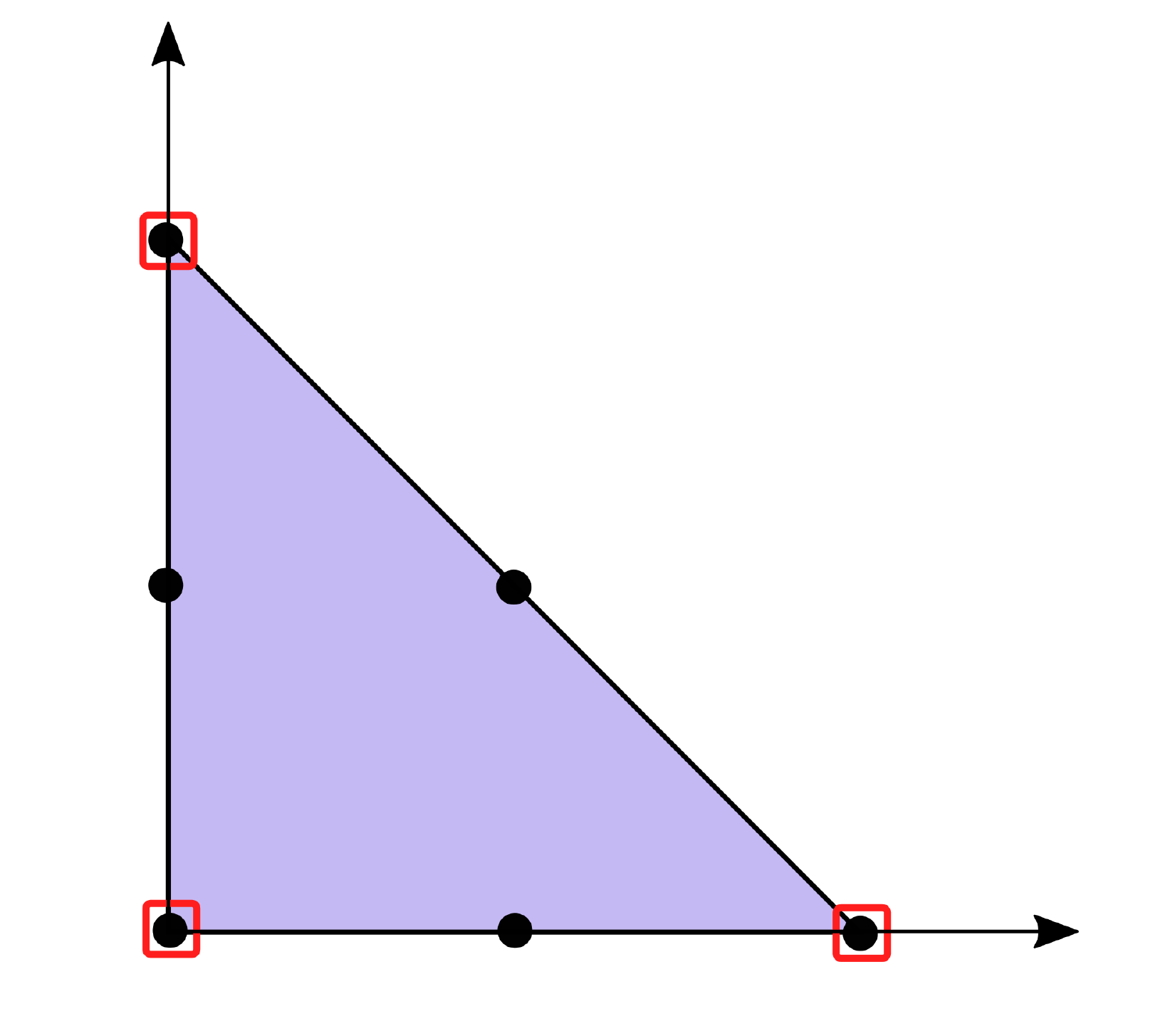}}
	\end{picture}
	\caption{T2T1}
		  \end{subfigure}
\begin{subfigure}[b]{0.32\textwidth}
	\begin{picture}(50,45)
	\put(0,0){\def\svgwidth{5cm}{\small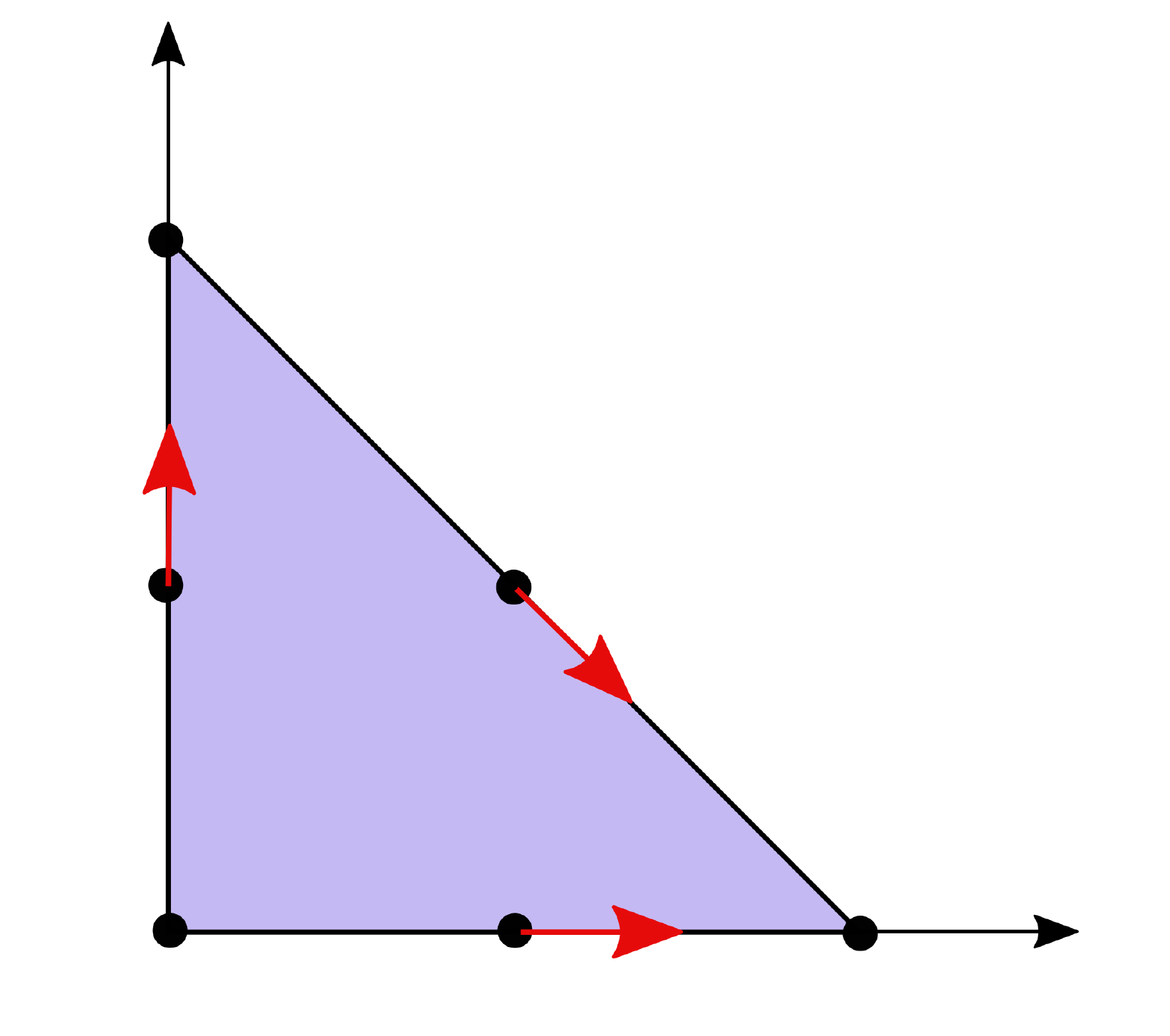}}
	\end{picture}
	\caption{T2NT1}
		  \end{subfigure}
  \begin{subfigure}[b]{0.32\textwidth}
	\begin{picture}(50,50)
	\put(0,0){\def\svgwidth{5cm}{\small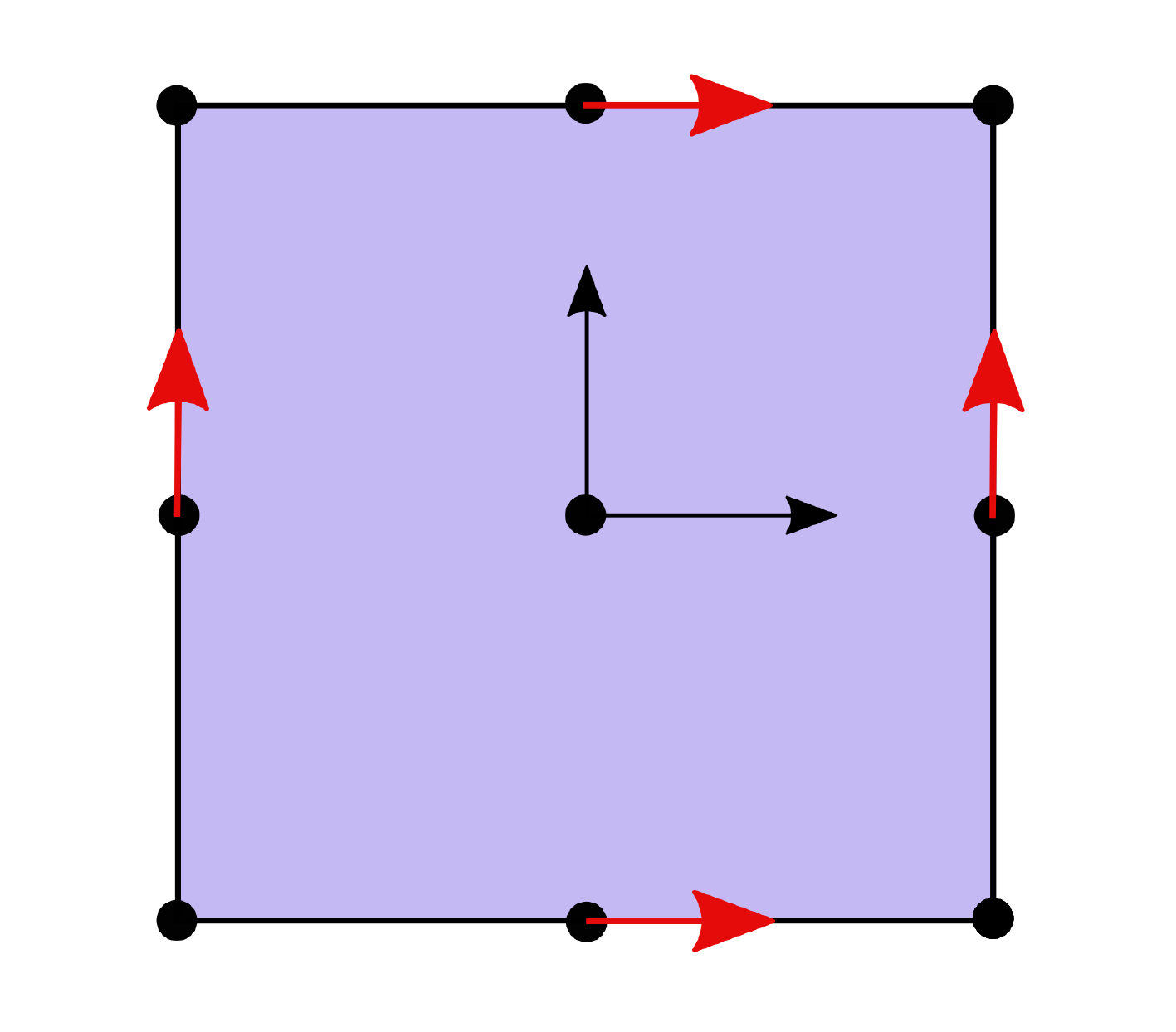}}
	\end{picture}
	\caption{Q2NQ1}
		  \end{subfigure}
  \begin{subfigure}[b]{0.32\textwidth}
	\begin{picture}(40,40)
	\put(0,0){\def\svgwidth{5cm}{\small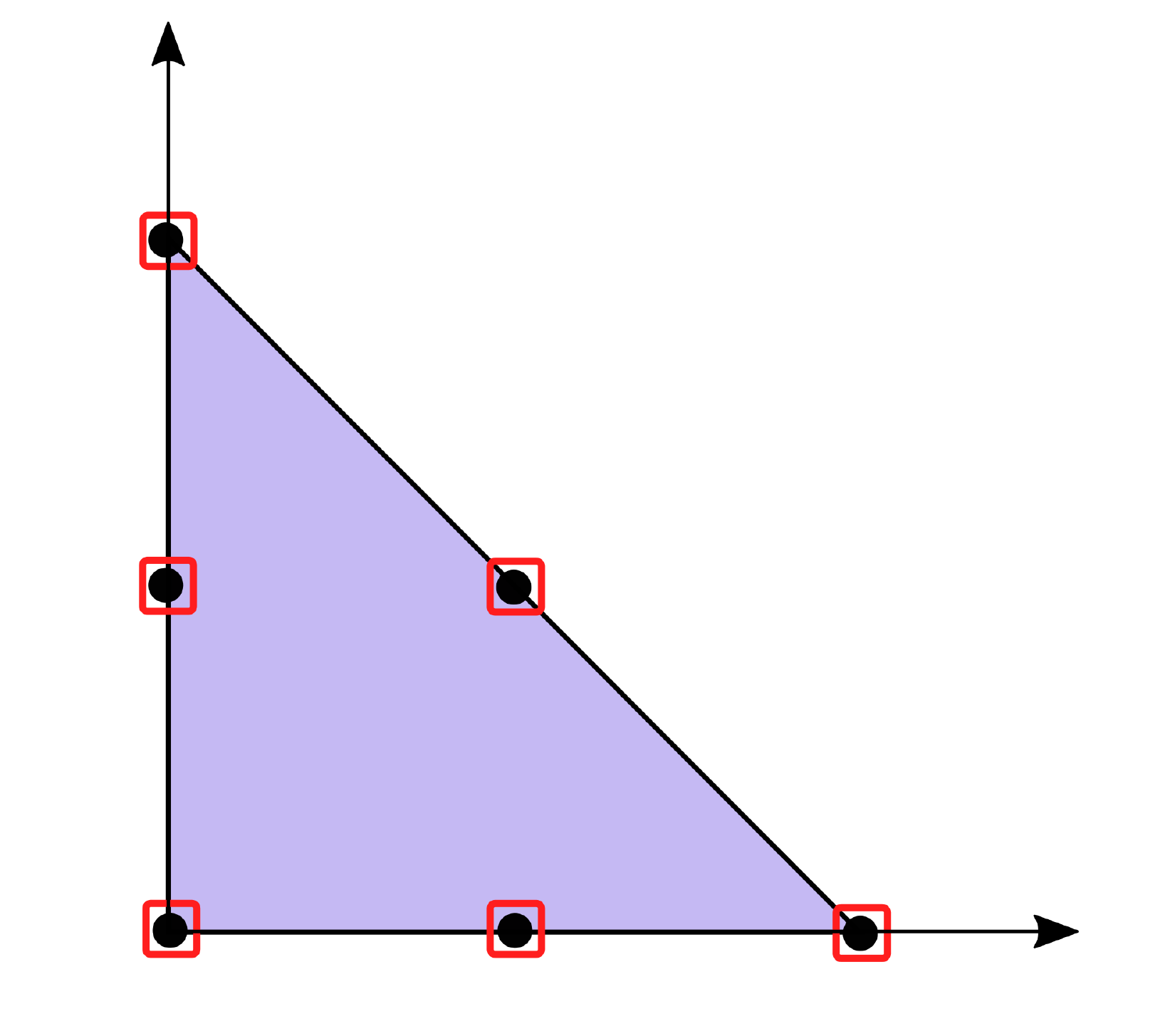}}
	\end{picture}
	\caption{T2T2}
		  \end{subfigure}	
		  	  \begin{subfigure}[b]{0.32\textwidth}
	\begin{picture}(50,45)
	\put(0,0){\def\svgwidth{5cm}{\small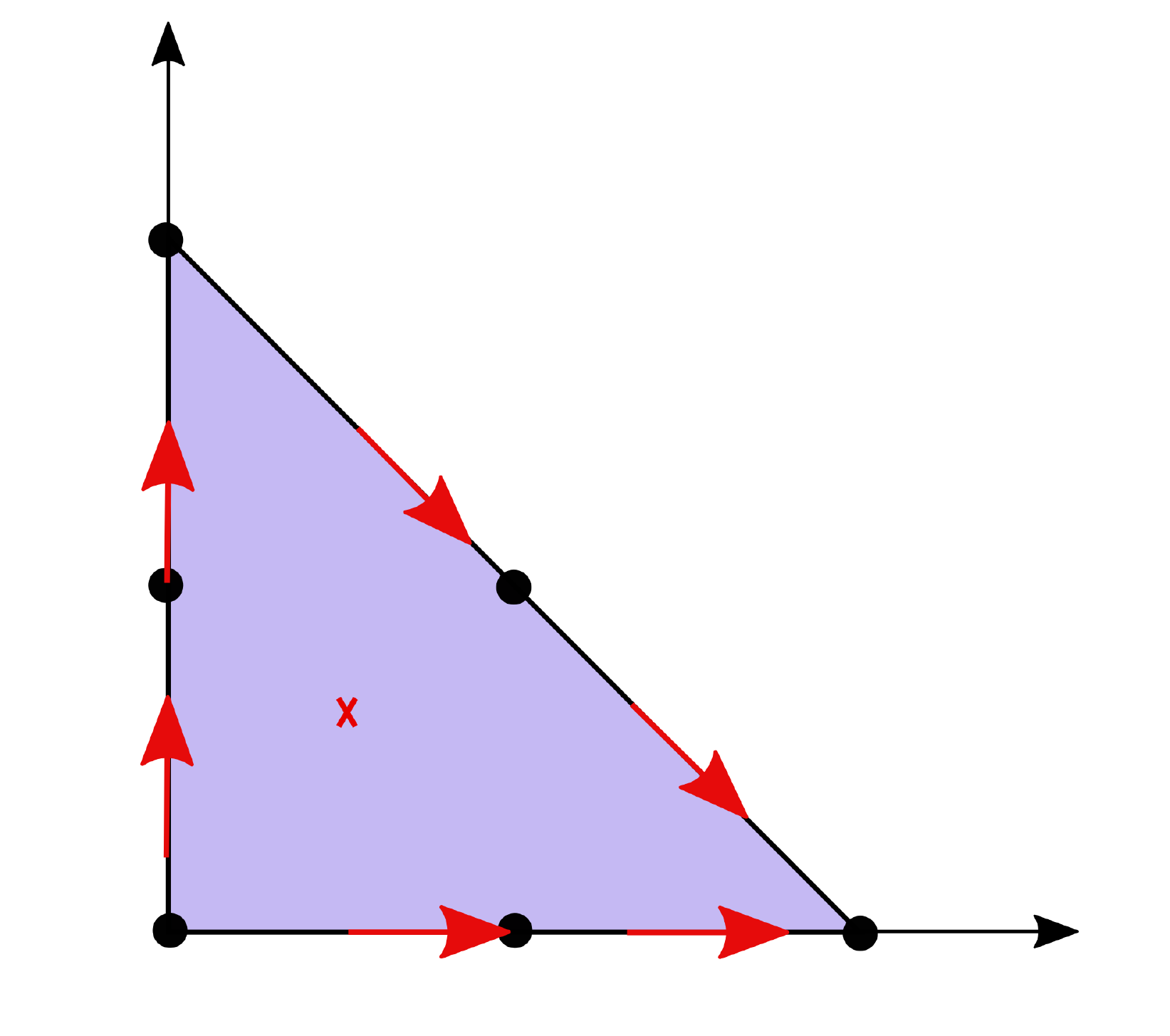}}
	\end{picture}
	\caption{T2NT2}
		  \end{subfigure}
 \begin{subfigure}[b]{0.32\textwidth}
	\begin{picture}(50,50)
	\put(0,0){\def\svgwidth{5cm}{\small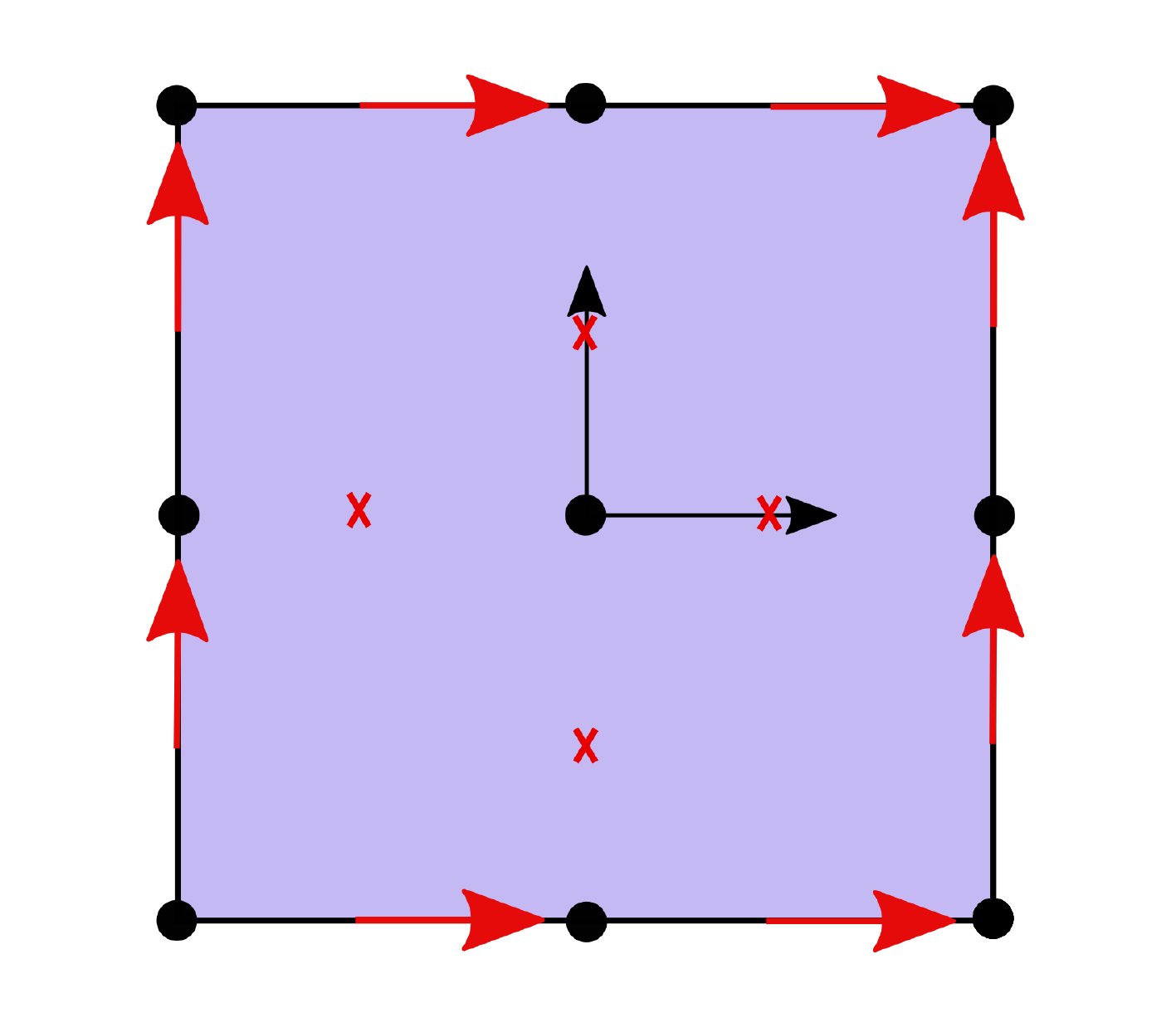}}
	\end{picture}
	\caption{Q2NQ2}
		  \end{subfigure} 
		  	\caption{The implemented finite elements in the parameter space. Black dots represent the displacement nodes while red squares stand for micro-distortion field nodes associated with tensorial dofs. Red arrows and crosses indicate the edge and inner vectorial dofs, respectively, of the micro-distortion field used in N\'ed\'elec formulation.}
\label{Fig:Finite_elements}
\end{figure}

\ssect{\hspace{-5mm}. Implemented finite elements \\} 
In this work, we present four nodal-edge elements based on the formulation in Section \ref{sec:mixelements} and two standard nodal elements based on Section \ref{sec:standardelements}.  All  implemented finite elements employ scalar quadratic shape functions of Lagrange-type for the displacement field approximation with the notation  T2 for triangles and Q2 for quadrilaterals.  The micro-distortion field is approximated using different formulations introduced in Sections \ref{sec:standardelements} and \ref{sec:mixelements}. For the standard nodal elements, Lagrange-type ansatz functions are used resulting in the element types T2T1 (linear ansatz for $\Bdis$) and  T2T2 (quadratic ansatz for $\Bdis$). Different nodal-edge elements  are built utilizing first- and second-order N\'ed\'elec formulations with tangential-conforming shape functions denoted as NT1 and NT2 for triangular elements and QT1 and QT2 for quadrilateral elements. The micro-distortion dofs in the standard nodal elements T2T1 and T2T2 are tensorial with $2 \times 2$ entries while the nodal-edge elements use vectorial dofs for the micro-distortion field which represent the tangential components. The full micro-distortion tensor is restored based on Equation (\ref{eq:t:dis_restoring}).   The used finite elements are depicted in the parameter space in Figure \ref{Fig:Finite_elements}. 

\section{Numerical examples} 
\label{sec:num}
\begin{figure}[ht]
\center
	\unitlength=1mm
	\begin{picture}(100,60)
	\put(0,0){\def\svgwidth{9.6cm}{\small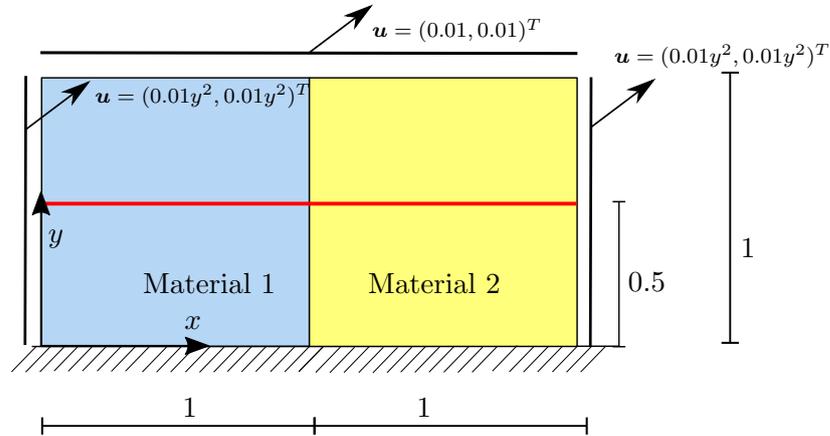}}
	\end{picture}
	\caption{2D rectangular bvp with two different materials. Inspection line, see Figures \ref{fig:example1:convergence1} and \ref{fig:example1:convergence2}, can be seen in red color. }
	\label{Figure:Geo1}
\end{figure}

\begin{figure}[ht]
     \centering
     \begin{subfigure}[b]{0.45\textwidth}
         \centering
         \includegraphics[width=\textwidth]{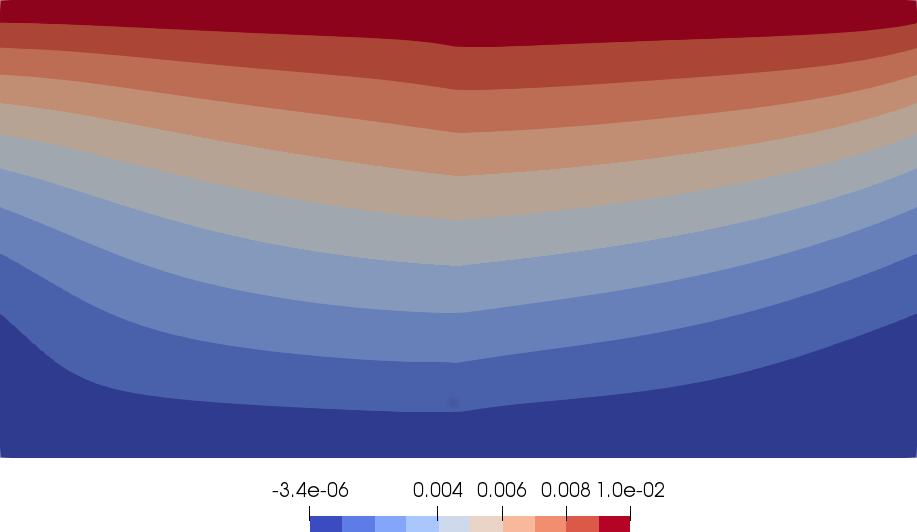}
			\caption{$u_1$}
     \end{subfigure}
     \centering
     \begin{subfigure}[b]{0.45\textwidth}
         \centering
         \includegraphics[width=\textwidth]{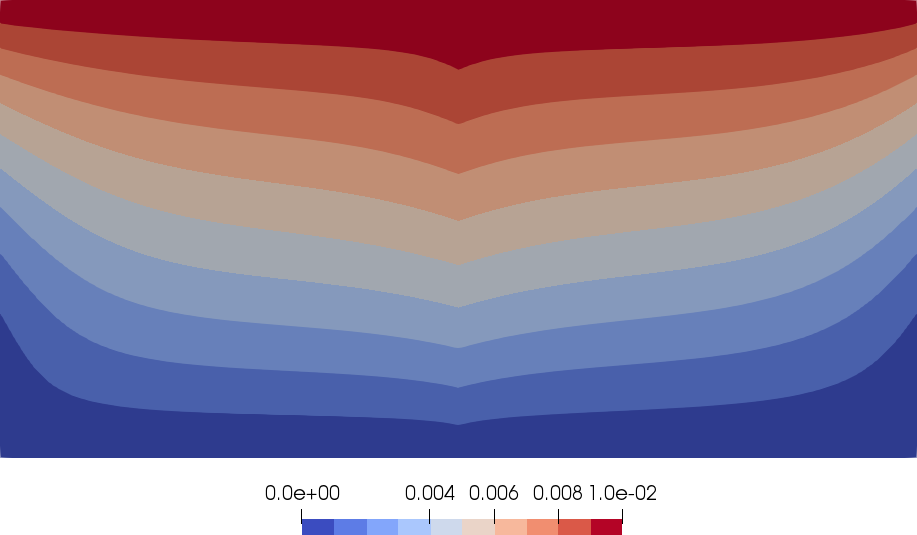}
			\caption{$u_2$}
     \end{subfigure}
     \begin{subfigure}[b]{0.45\textwidth}
         \centering
         \includegraphics[width=\textwidth]{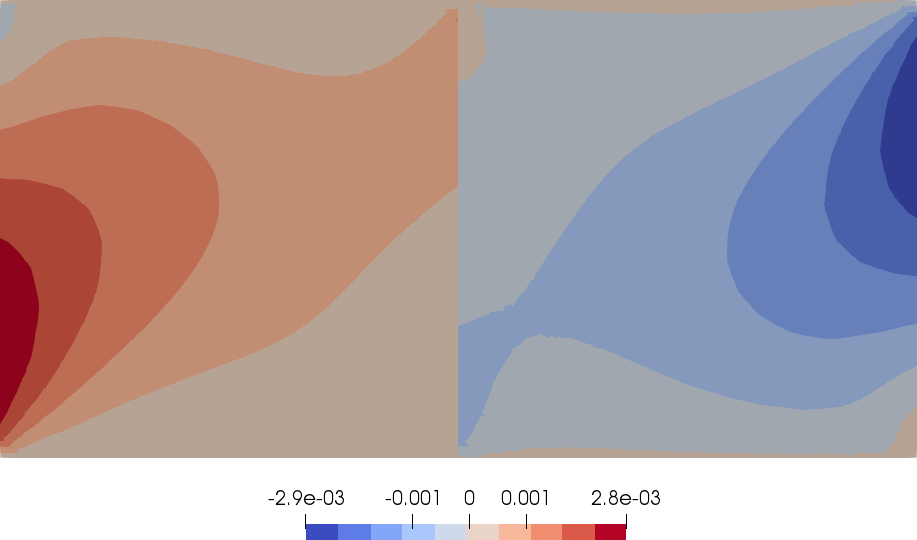}
			\caption{$\dis_{11}$}
     \end{subfigure}
     \begin{subfigure}[b]{0.45\textwidth}
         \centering
         \includegraphics[width=\textwidth]{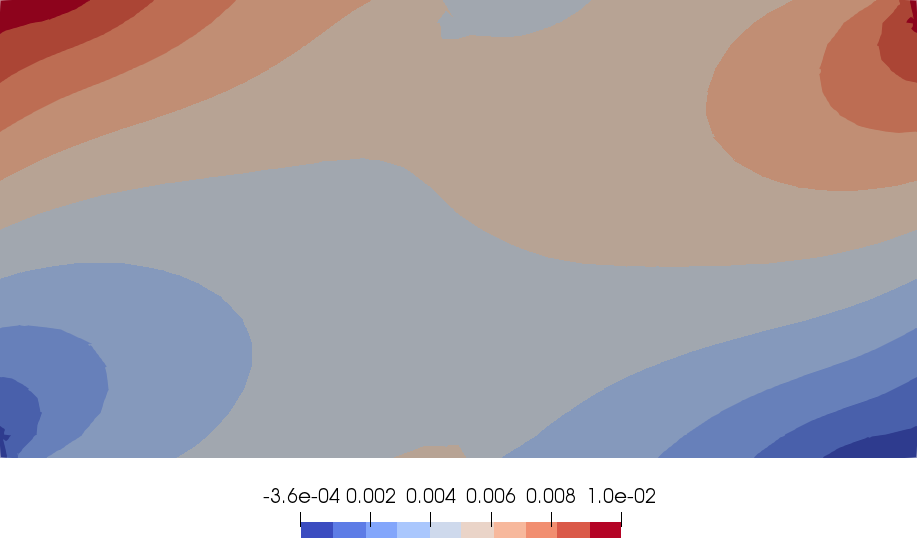}
			\caption{$\dis_{12}$}
     \end{subfigure}
     \begin{subfigure}[b]{0.45\textwidth}
         \centering
         \includegraphics[width=\textwidth]{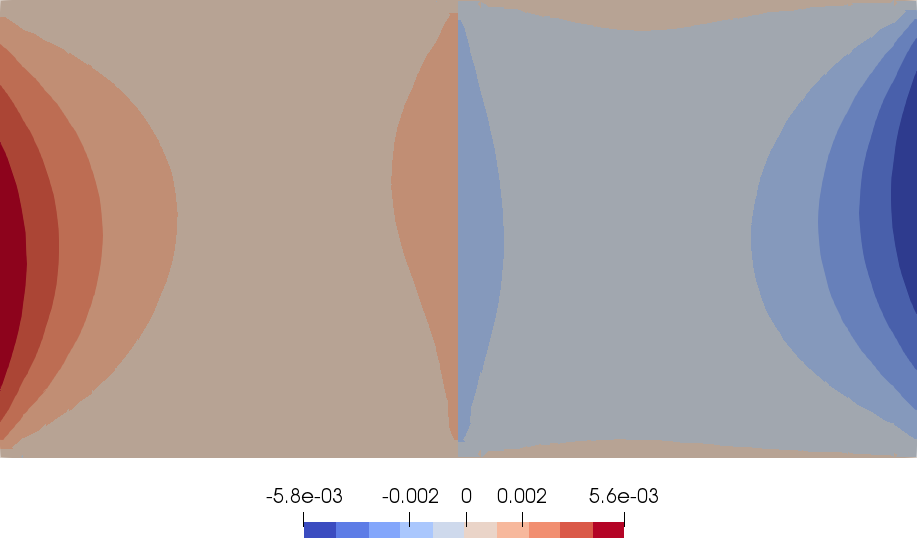}
			\caption{$\dis_{21}$}
     \end{subfigure}
     \begin{subfigure}[b]{0.45\textwidth}
         \centering
         \includegraphics[width=\textwidth]{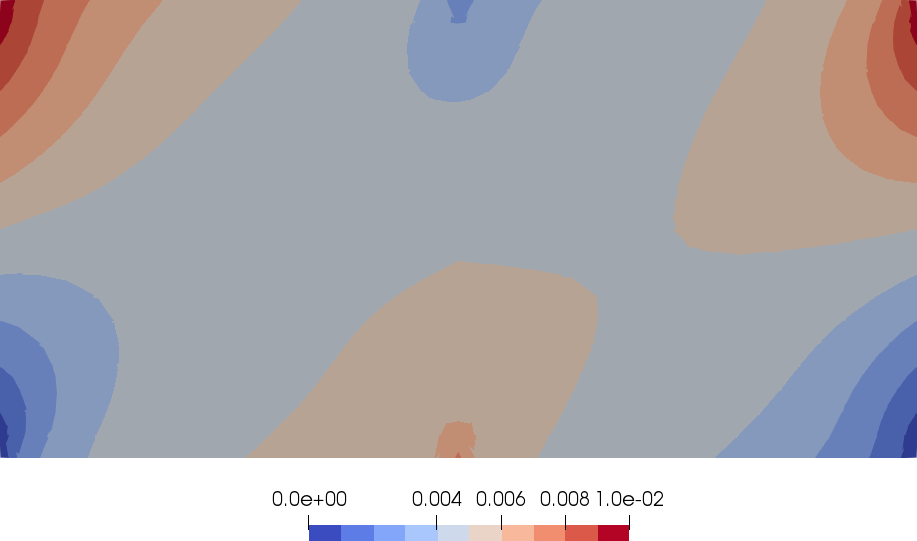}
			\caption{$\dis_{22}$}
     \end{subfigure}
        \caption{The displacement and the micro-distortion fields of the first boundary problem.}
\label{fig:example1:solution}
\end{figure}

For our numerical examples, we neglect the body forces and moments, i.e. $\overline\bbf =\bzero$ and $\overline\bM = \bzero$, and assume isotropic material behavior which can be described by the set of material parameters  $\lambda_\textrm{micro},\mu_\textrm{micro},\lambda_e,\mu_e,\mu_c ,\mu$ and $\Lc$, where $\lambda_*$ and $\mu_*$ denote the Lam\'e coefficients.  Furthermore, we assume $\IL$ as a fourth-order identity tensor $\II$. Throughout all examples, we consider the Cosserat modulus $\mu_c = 0$, cf.  \cite{NefGhiLazMad:2015:trl,Nef:2006:tcc}, leading to the symmetry of the force stress tensor. The simulations presented in this paper are performed within AceGen and AceFEM programs, which are developed and maintained by  Jo\v{z}e Korelc (University of Ljubljana). The interested reader is referred to \cite{KorWri:2016:aofem}.
\ssect{\hspace{-5mm}. Discontinuous solution: an interface between two different materials} 
\label{bvp:1}
In the first boundary value problem (bvp), we consider a rectangular domain $\B$ with length $l=2$ and height $h=1$ which consists of a side-by-side arrangement of two different materials, see Figure \ref{Figure:Geo1}. The bottom edge is fixed in both directions, $\overline\bu=(0,0)^T$, and we apply a displacement to the upper edge, $\overline\bu=(0.01,0.01)^T$, while the left and right edges are subjected to displacement given by $\overline\bu=(0.01y^2,0.01y^2)^T$. For the micro-distortion field, the consistent coupling  boundary condition $\Bdis \cdot \Btau = \nabla \overline\bu \cdot \Btau$ is enforced on the entire boundary $\partial \B$. The material parameters are given in Table \ref{tab:example1:MP}. Due to the different material parameters of the domains, a discontinuous solution in the micro-distortion is expected which allows us to compare and evaluate the behavior of the implemented finite elements.  \\
\begin{table}[ht]
\center
\begin{tabular}{|c|c|}
\hline
Material 1 & Material 2 \\
\hline
\begin{minipage}[t]{7cm}
	$ \lambda_\textrm{micro} = 555.55 , \quad \mu_\textrm{micro} = 833.33$ \\ $\lambda_e = 486.11, \quad \mu_e = 729.17$\\
	$\mu_c = 0, \quad \mu=833.33$ \\ 
	$\IL = \II, \quad \Lc=1$ 
\end{minipage} & 
\begin{minipage}[t]{7cm}
	$ \lambda_\textrm{micro} = 1111.11 , \quad \mu_\textrm{micro} = 1667.67$ \\ $\lambda_e = 972.22, \quad \mu_e = 1458.33$\\
	$\mu_c = 0, \quad \mu=1666.67$ \\
	 $\IL = \II, \quad \Lc=1$ 
\end{minipage} \\
\hline
\end{tabular}
\caption{Material parameters of the first boundary problem, see Figure \ref{Figure:Geo1}.}
\label{tab:example1:MP}
\end{table}
In Figure \ref{fig:example1:solution}, the displacement and micro-distortion fields obtained for a discretization with 1600 N2NT2 finite elements are shown.  Note that the elements solution is plotted in this work without the usual averaging or smoothing on the element edges in order to investigate the possible discontinuities. The solution of the displacement field is continuous through the interface  between the two domains. The tangential components of the micro-distortion $\Bdis \cdot \be_2 = (\dis_{12},\dis_{22})^T$ are continuous on the interface while the normal components $\Bdis \cdot \be_1 = (\dis_{11},\dis_{21})^T$ exhibit discontinuities. 

\begin{figure}[ht]
	\unitlength=1mm
	\center
	\unitlength=1mm
	\center		  
	  	  \begin{subfigure}[b]{0.45\textwidth}
         \includegraphics[width=\textwidth]{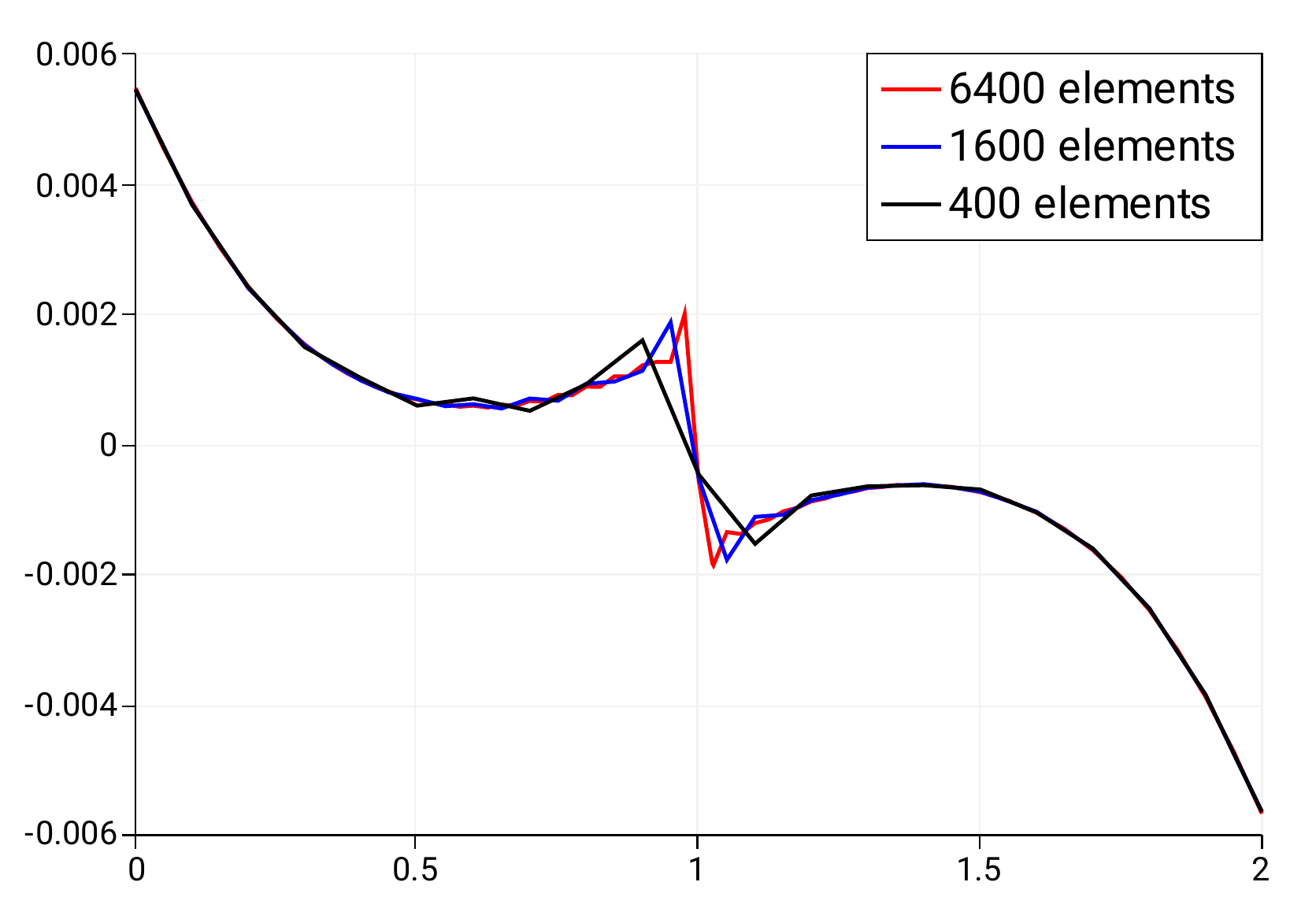}
         \put(-38,-2){$x$}
            \put(-75,25){$\dis_{21}$}
	\caption{T2T1}
		  \end{subfigure}
	  	  \begin{subfigure}[b]{0.45\textwidth}
         \includegraphics[width=\textwidth]{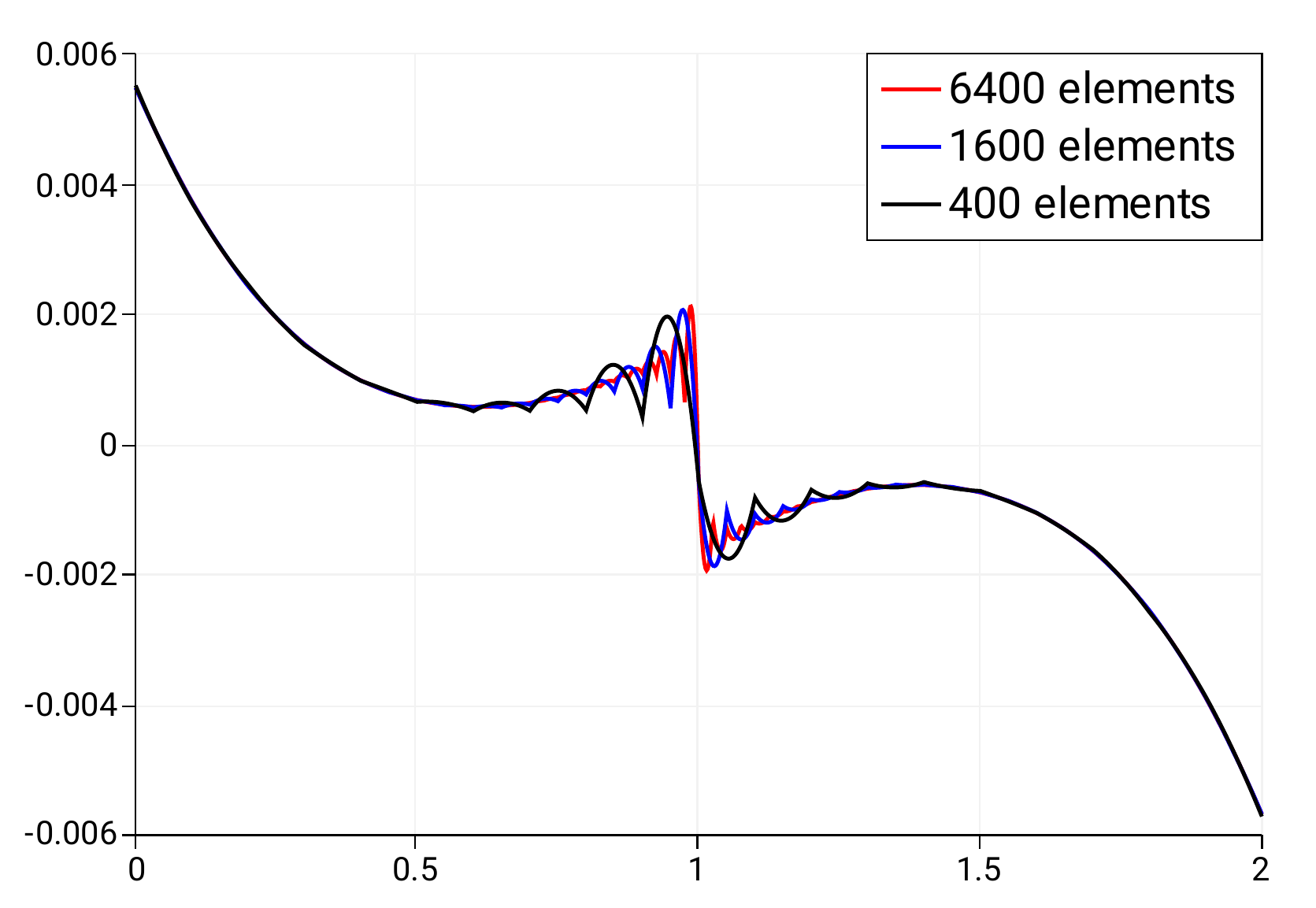}
         \put(-38,-2){$x$}
            \put(-75,25){$\dis_{21}$}
	\caption{T2T2}
		  \end{subfigure}
		  	\caption{Illustration of $\dis_{21}$ along the inspection line $y=0.5$ using the nodal elements with different mesh densities.}
\label{fig:example1:convergence1}
\end{figure}

\begin{figure}[ht]
	\unitlength=1mm
	\center		  
	  	  \begin{subfigure}[b]{0.45\textwidth}
         \includegraphics[width=\textwidth]{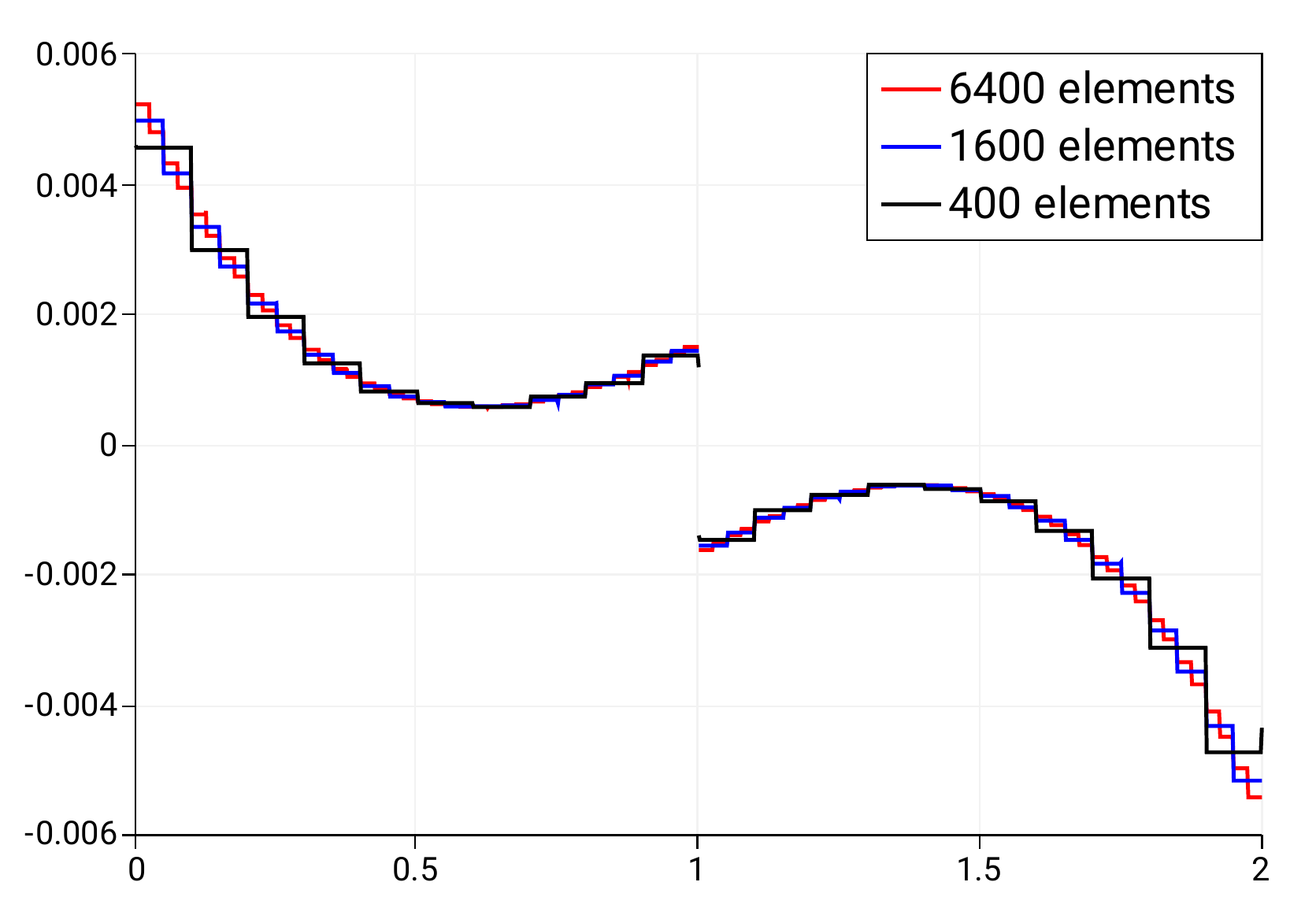}
         \put(-38,-2){$x$}
            \put(-75,25){$\dis_{21}$}
	\caption{T2NT1}
		  \end{subfigure}
	  	  \begin{subfigure}[b]{0.45\textwidth}
         \includegraphics[width=\textwidth]{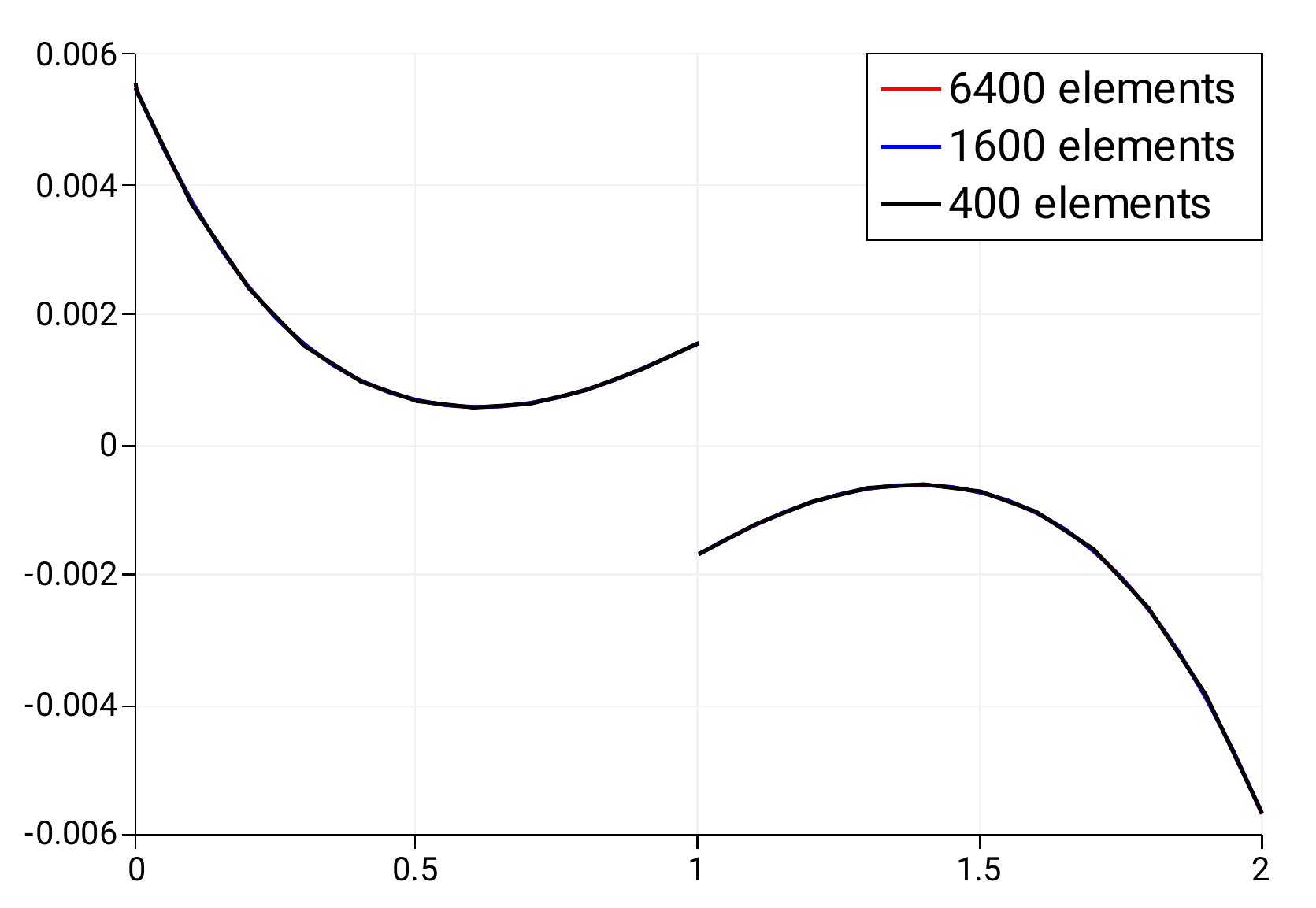}
         \put(-38,-2){$x$}
            \put(-75,25){$\dis_{21}$}
	\caption{T2NT2}
		  \end{subfigure}
		  	  	  \begin{subfigure}[b]{0.45\textwidth}
         \includegraphics[width=\textwidth]{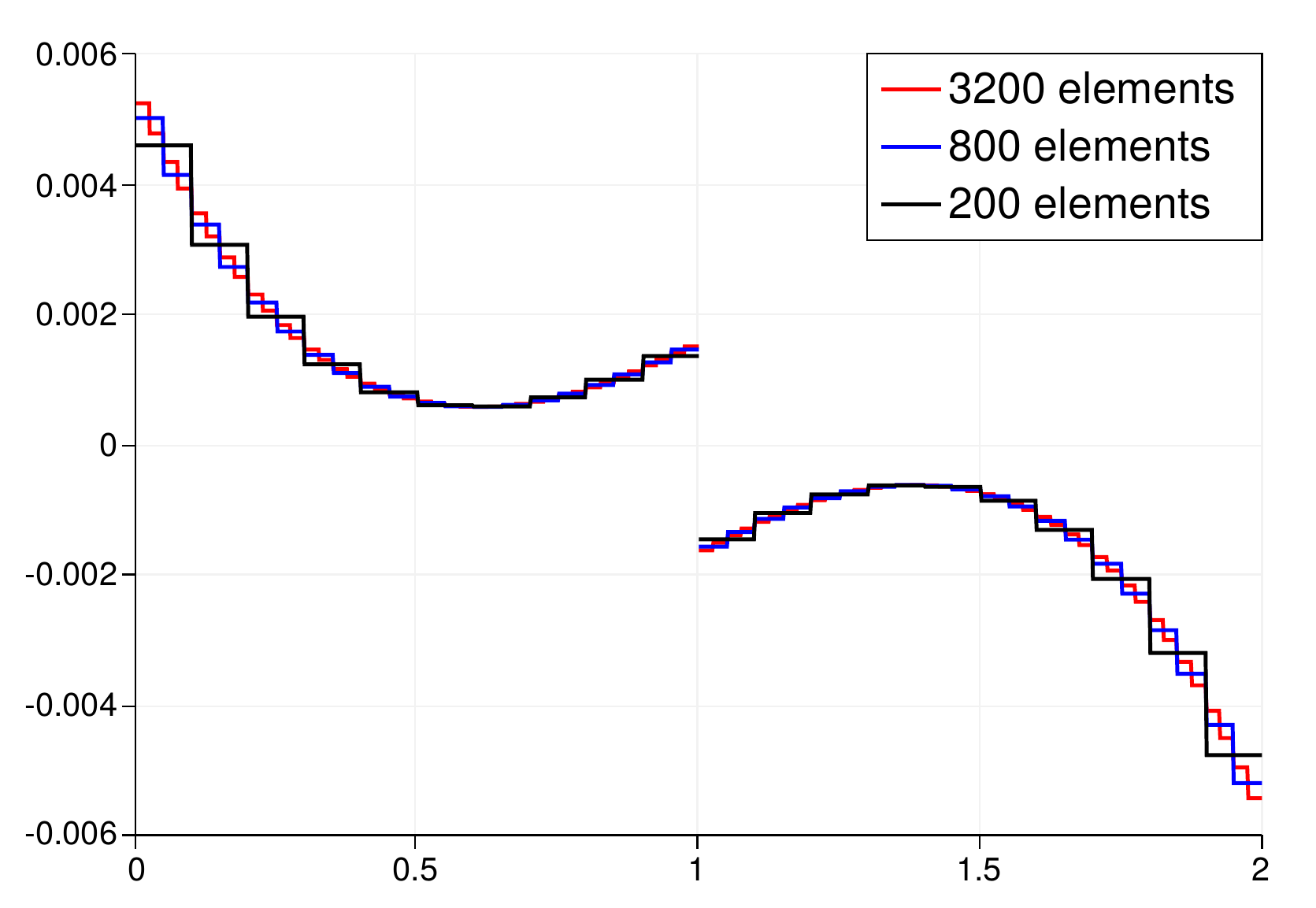}
         \put(-38,-2){$x$}
            \put(-75,25){$\dis_{21}$}
	\caption{Q2NQ1}
		  \end{subfigure}
		  	  	  \begin{subfigure}[b]{0.45\textwidth}
         \includegraphics[width=\textwidth]{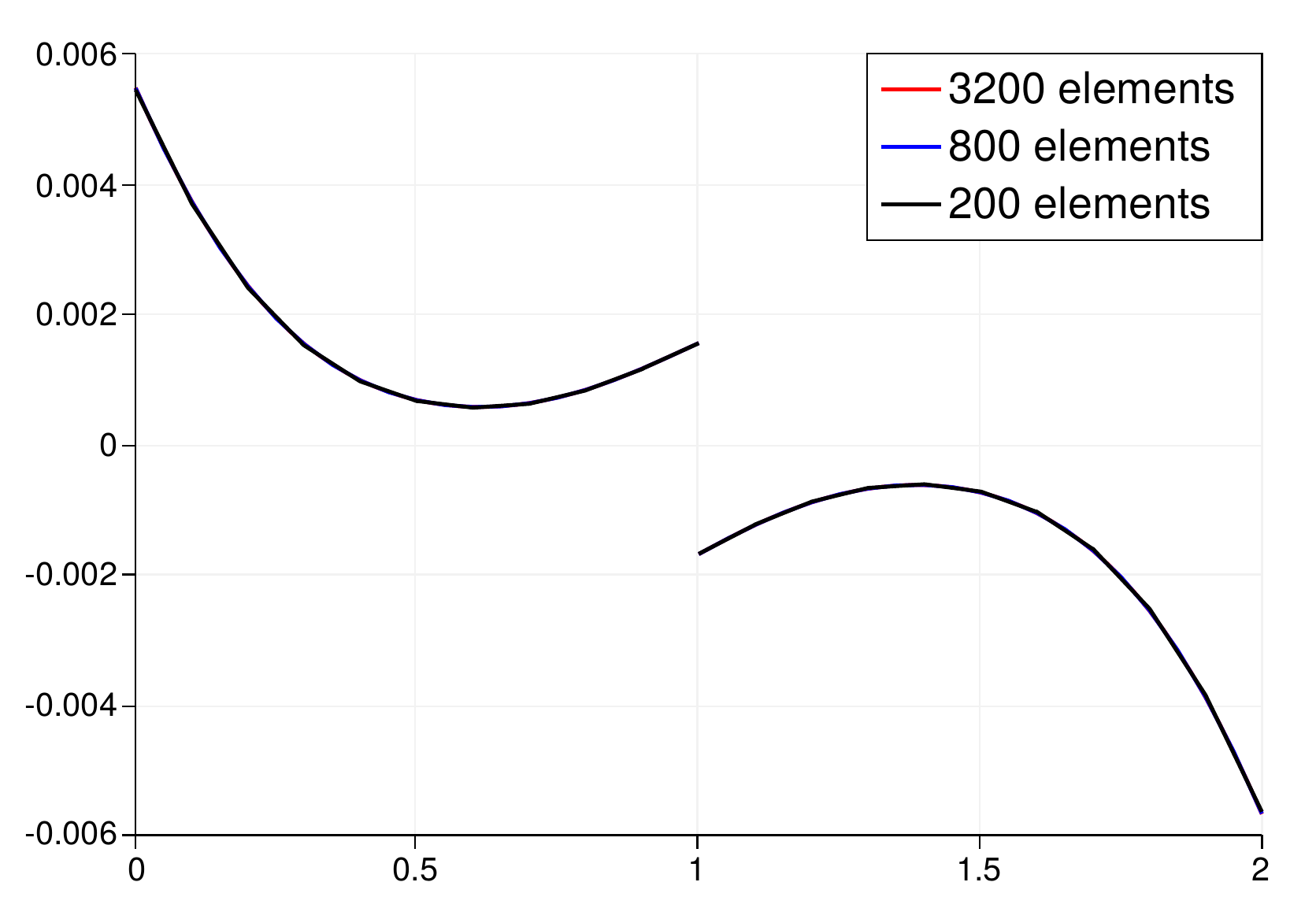}
         \put(-38,-2){$x$}
            \put(-75,25){$\dis_{21}$}
	\caption{Q2NQ2}
		  \end{subfigure}		  
		  	\caption{Illustration of $\dis_{21}$ along the inspection line $y=0.5$ using $H^1(\B) \times H(\curl, \B)$ elements with different mesh densities.}
\label{fig:example1:convergence2}
\end{figure}
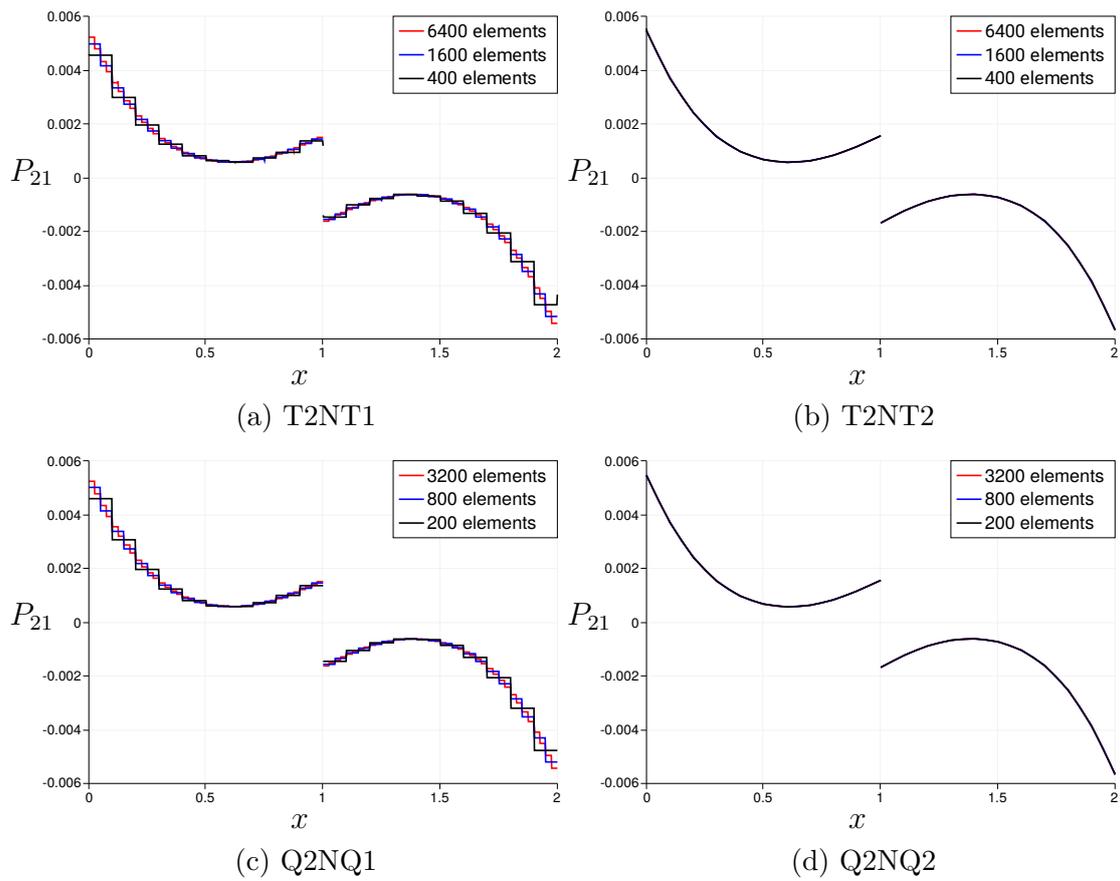

The component  $\dis_{21}$ along the inspection line $y=0.5$ is plotted in Figure \ref{fig:example1:convergence1}   using $H^1(\B) \times H^1(\B)$ elements  resulting in a continuous solution. This causes a transition zone at the interface which needs to be resolved by increasing the mesh density tremendously in order to approximate the discontinuous solution at the interface. Enhancing the approximation space of the micro-distortion field to the second-order Lagrange space does not improve the convergence behavior to the discontinuous  solution at the interface. However, the discontinuous solution of $\dis_{21}$ can be captured by the $H^1(\B) \times H(\curl, \B)$ elements, see Figure \ref{fig:example1:convergence2}. The second-order N\'ed\'elec formulations T2NT2 and Q2NQ2 shows an instant convergence with a coarse mesh while first-order N\'ed\'elec formulations T2NT1 and Q2NQ1 reach a good accuracy with intermediate/coarse mesh showing the expected linear behavior for the component within one element. 

\FloatBarrier

\ssect{\hspace{-5mm}. Characteristic length analysis: pure shear problem} 

\begin{figure}[ht]
\center
\begin{subfigure}[b]{0.49\textwidth}
	\unitlength=1mm
	\begin{picture}(100,70)
	\put(10,0){\def\svgwidth{6cm}{\small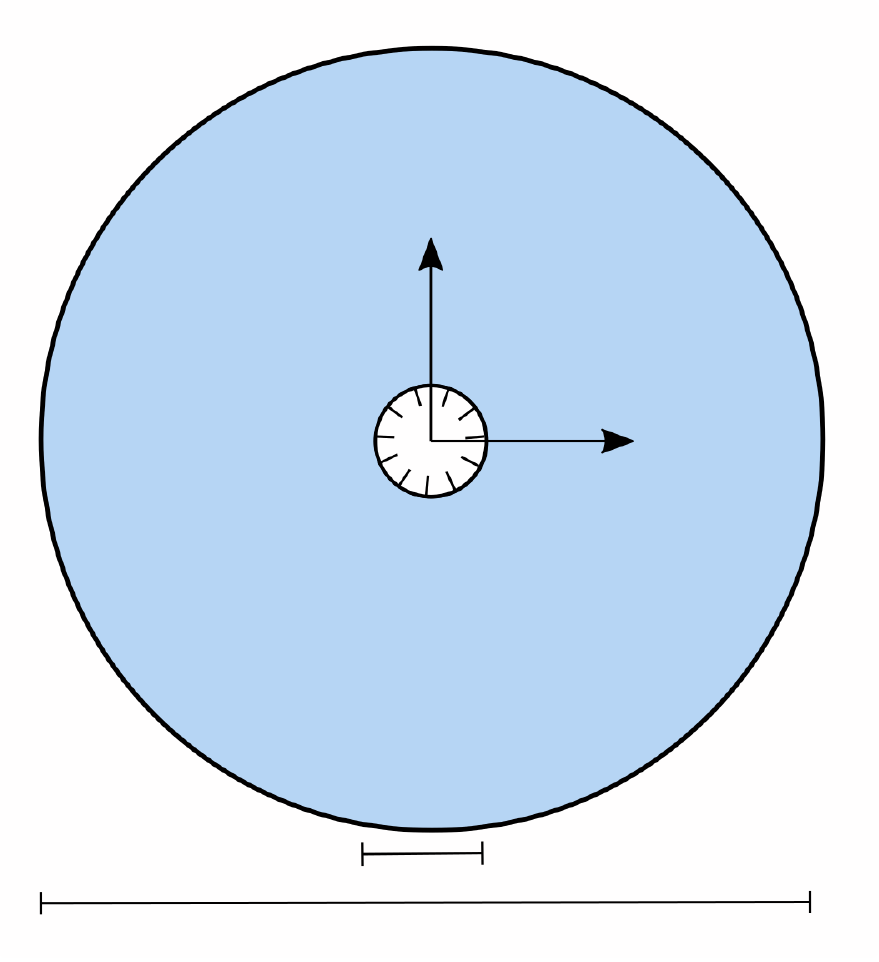}}
	\end{picture}
\label{fig:example2:BVPA}
\caption{case A}
\end{subfigure}
\begin{subfigure}[b]{0.49\textwidth}
	\unitlength=1mm
	\begin{picture}(100,70)
	\put(10,0){\def\svgwidth{6cm}{\small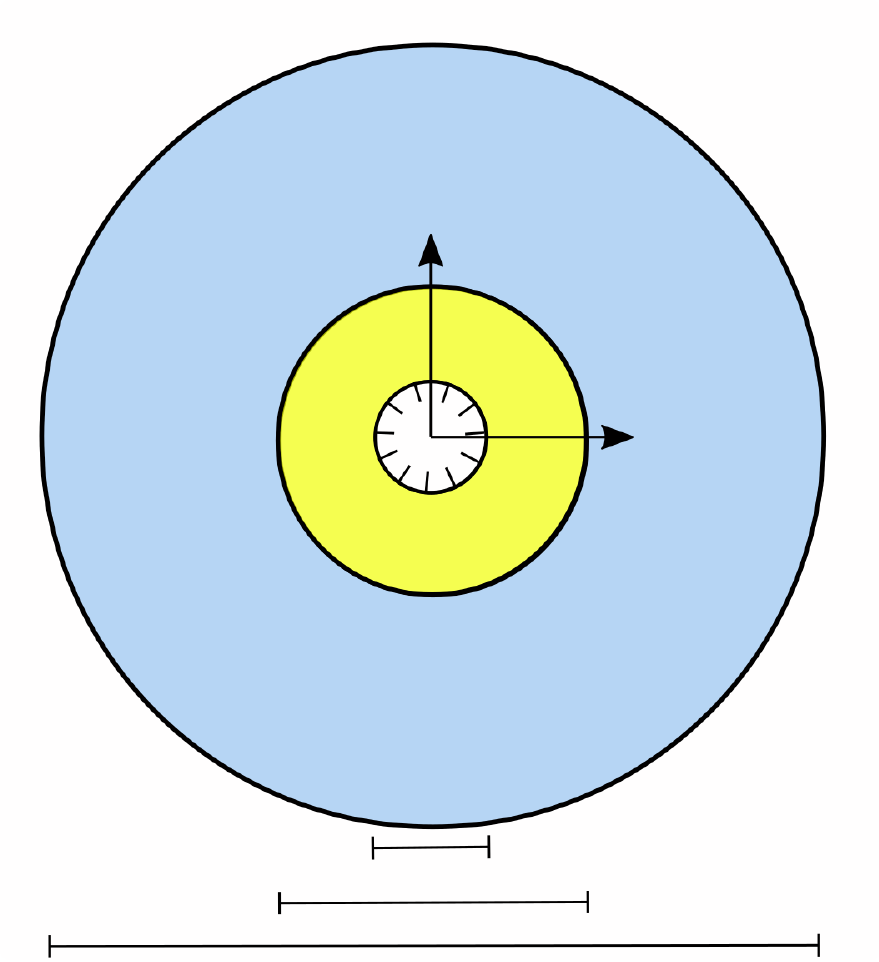}}
	\end{picture}
\label{fig:example2:BVPB}
\caption{case B}
\end{subfigure}
	\caption{ The geometry of the second boundary value problem. }
	\label{Figure:Geo2}
\end{figure}

We introduce a second boundary value problem, see Figure  \ref{Figure:Geo2}, consisting of a circular domain $\B$ with a radius $r_\textrm{o}= 25$ and a circular hole at its center with a radius $r_\textrm{i} = 2$. We fix the displacement field $ \overline{\bu} = \bzero$ on the inner boundary $\partial \B_\textrm{i}$ and we rotate the outer boundary $\partial \B_\textrm{o}$ counter clockwise with $ \bar{\bu} = (-\frac{\Delta}{r_\textrm{o}} y,  \frac{\Delta}{r_\textrm{o}} x)^T$ where $\Delta = 0.01$. For the micro-distortion field, we apply the consistent coupling  boundary condition  ($\Bdis \cdot \Btau = \nabla  \overline{\bu} \cdot \Btau$)  on all boundaries  $\partial \B = \partial \B_\textrm{i} \cup  \partial \B_\textrm{o} $.   Two different cases are discussed in the following. For case A, a single material is assumed whereas for case B two materials are considered. The second material is located as a ring with an outer radius $r_\textrm{m} = 10$ and an inner radius $r_\textrm{i} = 2$. The material parameters are shown in Table \ref{tab:example2:MP}. For the analysis of the influence of the characteristic length $\Lc$, the characteristic length will be varied.  

\begin{table}[ht]
\center
\begin{tabular}{|c|c|}
\hline
Material 1 & Material 2 \\
\hline
\begin{minipage}[t]{7cm}
	$ \lambda_\textrm{micro} = 555.55 , \quad \mu_\textrm{micro} = 833.33$ \\ $\lambda_e = 486.11, \quad \mu_e = 729.17$\\
	$\mu_c = 0, \quad \mu=833.33$ \\ 
	$\IL = \II, \quad \Lc \in \{0.001, 5, 1000\}$ 
\end{minipage} & 
\begin{minipage}[t]{7cm}
	$ \lambda_\textrm{micro} = 2777.78 , \quad \mu_\textrm{micro} = 4166.67$ \\ $\lambda_e = 2430.555, \quad \mu_e = 3645.85$\\
	$\mu_c = 0, \quad \mu=4166.67$ \\
	 $\IL = \II, \quad \Lc \in \{0.001, 5, 1000\}$ 
\end{minipage} \\
\hline
\end{tabular}
\caption{Material parameters of the second boundary problem, see Figure \ref{Figure:Geo2}.}
\label{tab:example2:MP}
\end{table}

\begin{figure}[ht]
     \centering
     \begin{subfigure}[b]{0.4\textwidth}
         \centering
         \includegraphics[width=\textwidth]{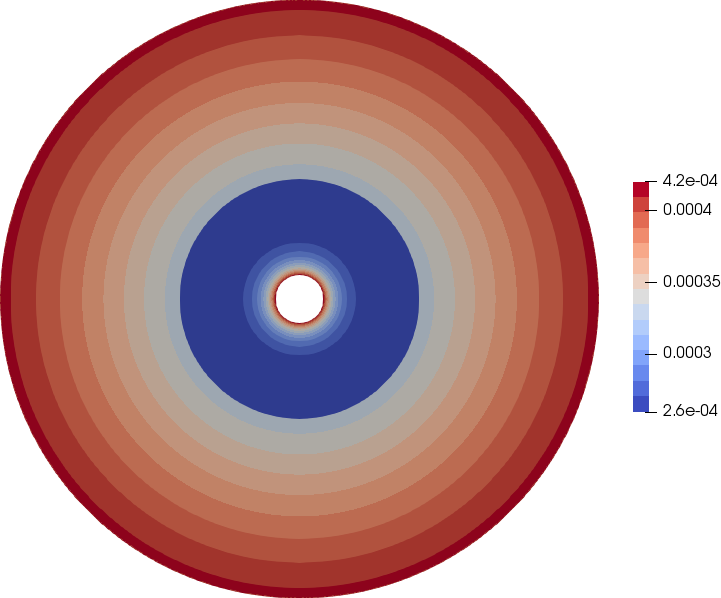}
			\caption{$\dis_{\theta r}$}
     \end{subfigure}
     \hspace{2 cm}
     \centering
     \begin{subfigure}[b]{0.4\textwidth}
         \centering
         \includegraphics[width=\textwidth]{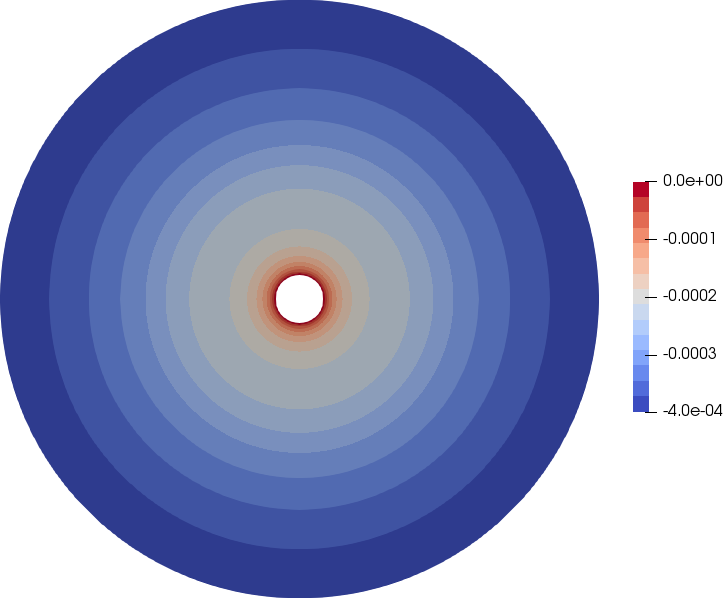}
			\caption{$\dis_{r \theta}$}
     \end{subfigure}
        \caption{The non-vanishing micro-distortion components of the second boundary problem using 30256 Q2NQ2 elements for $\Lc=5$.}
\label{fig:example2:solution}
\end{figure}

\begin{figure}[ht]
\unitlength=1mm
	\center
	  	  \begin{subfigure}[b]{0.45\textwidth}
      \includegraphics[width=\textwidth]{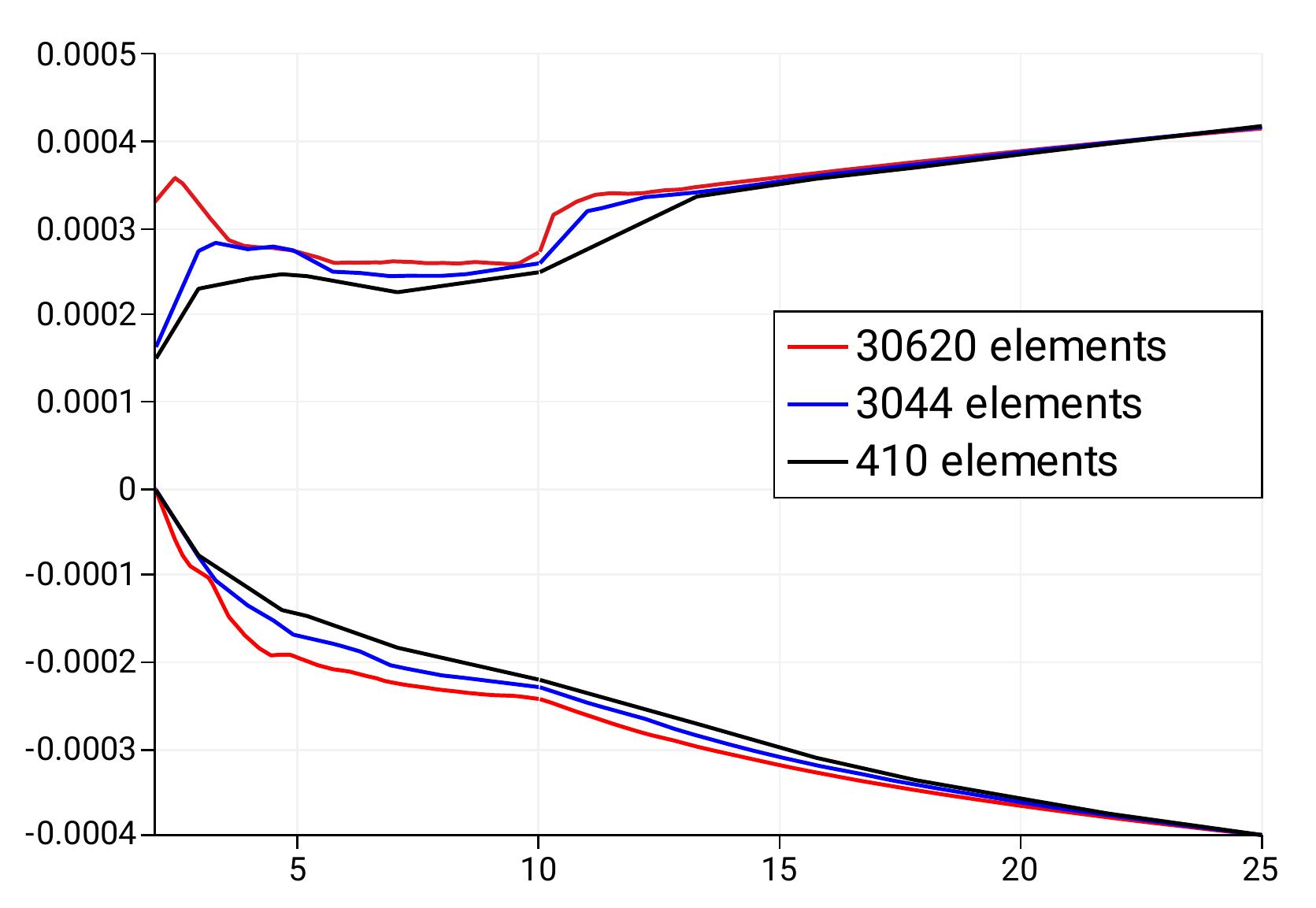}
             \put(-35.903717, -1.2734547){\color[rgb]{0,0,0}\makebox(0,0)[lb]{\smash{$r$}}}%
     \put(-37.903717,12.2734547){\color[rgb]{0,0,0}\makebox(0,0)[lb]{\smash{$\dis_{r \theta}$}}}%
 \put(-37.903717,42.2734547){\color[rgb]{0,0,0}\makebox(0,0)[lb]{\smash{$\dis_{\theta r}$}}}%
	\caption{T2T1}
		  \end{subfigure} 
	  \begin{subfigure}[b]{0.45\textwidth}
      \includegraphics[width=\textwidth]{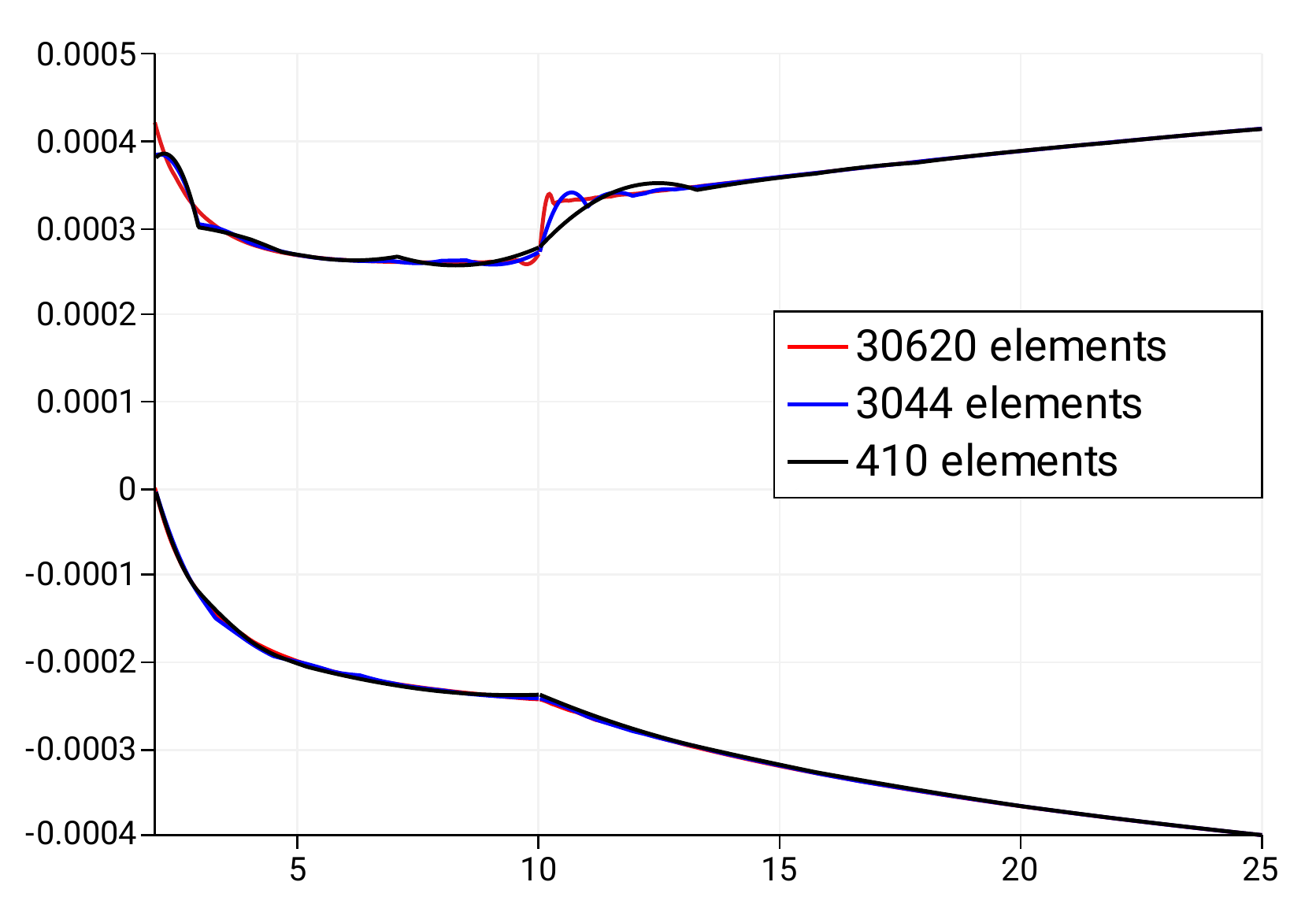}
             \put(-35.903717, -1.2734547){\color[rgb]{0,0,0}\makebox(0,0)[lb]{\smash{$r$}}}%
     \put(-37.903717,12.2734547){\color[rgb]{0,0,0}\makebox(0,0)[lb]{\smash{$\dis_{r \theta}$}}}%
 \put(-37.903717,42.2734547){\color[rgb]{0,0,0}\makebox(0,0)[lb]{\smash{$\dis_{\theta r}$}}}%
	\caption{T2T2}
		  \end{subfigure}
	\caption{The non-vanishing components of $\bP$ along the radius for the  $H^1(\B) \times H^1(\B)$ elements using three mesh densities and $\Lc=5$.} 
		  	 \label{fig:example2:convergence1}
\end{figure}

\begin{figure}[ht]
	\unitlength=1mm
	\center
	  	  \begin{subfigure}[b]{0.45\textwidth}
      \includegraphics[width=\textwidth]{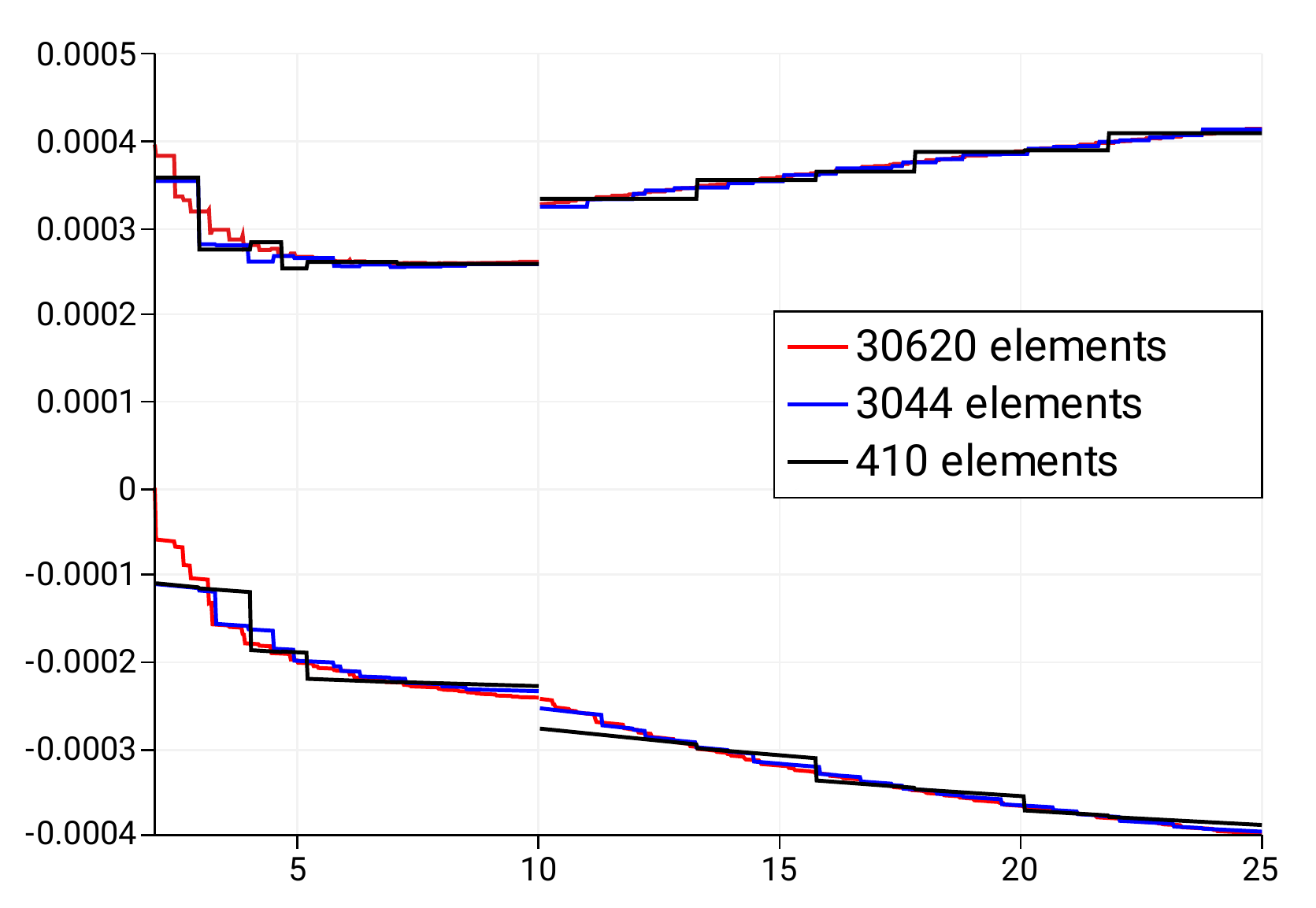}
             \put(-39.903717, -1.2734547){\color[rgb]{0,0,0}\makebox(0,0)[lb]{\smash{$r$}}}%
     \put(-37.903717,12.2734547){\color[rgb]{0,0,0}\makebox(0,0)[lb]{\smash{$\dis_{r \theta}$}}}%
 \put(-37.903717,42.2734547){\color[rgb]{0,0,0}\makebox(0,0)[lb]{\smash{$\dis_{\theta r}$}}}%
	\caption{T2NT1}
		  \end{subfigure} 
	  \begin{subfigure}[b]{0.45\textwidth}
      \includegraphics[width=\textwidth]{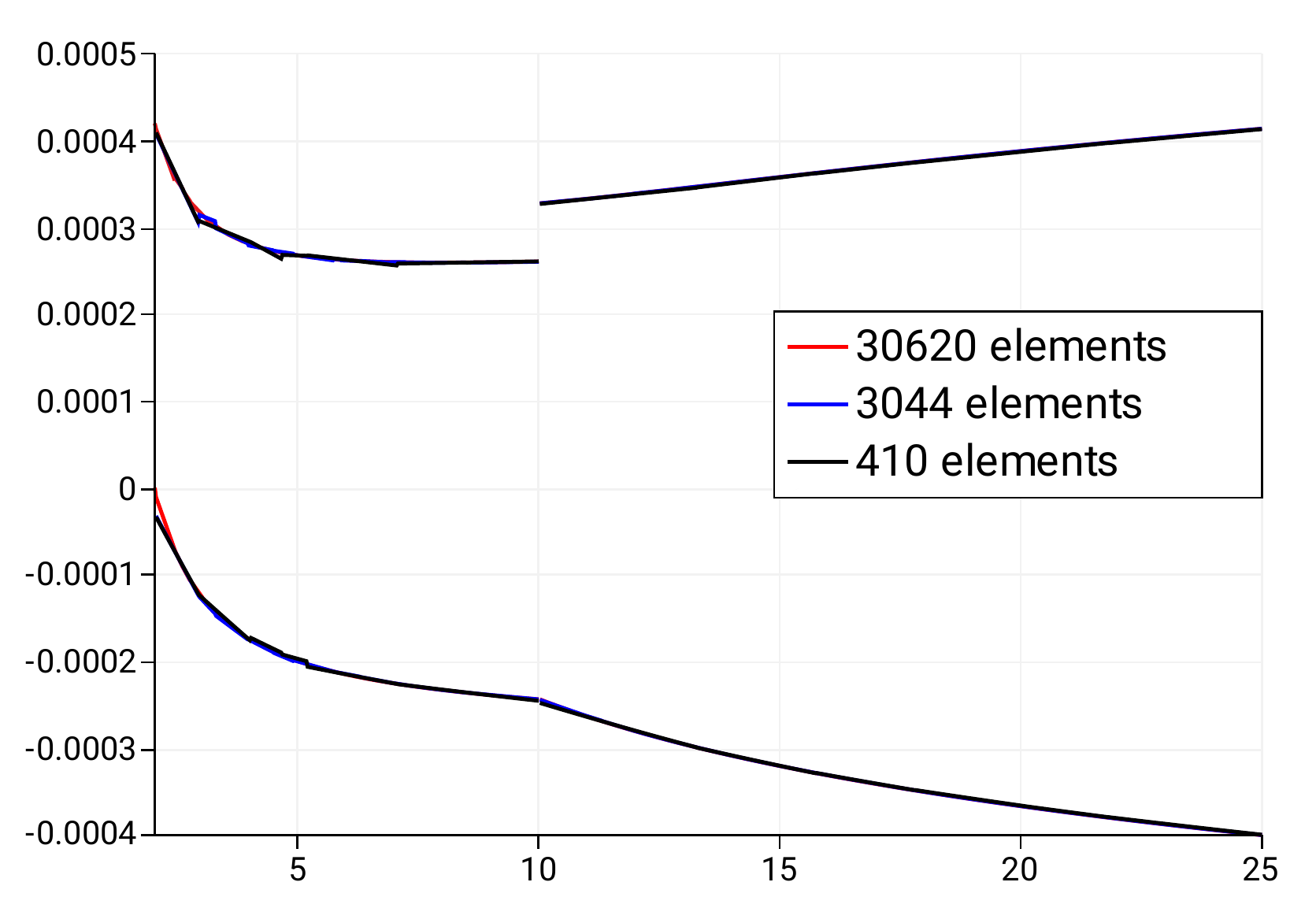}
             \put(-39.903717, -1.2734547){\color[rgb]{0,0,0}\makebox(0,0)[lb]{\smash{$r$}}}%
     \put(-37.903717,12.2734547){\color[rgb]{0,0,0}\makebox(0,0)[lb]{\smash{$\dis_{r \theta}$}}}%
 \put(-37.903717,42.2734547){\color[rgb]{0,0,0}\makebox(0,0)[lb]{\smash{$\dis_{\theta r}$}}}%
	\caption{T2NT2}
		  \end{subfigure}
	  	  \begin{subfigure}[b]{0.45\textwidth}
      \includegraphics[width=\textwidth]{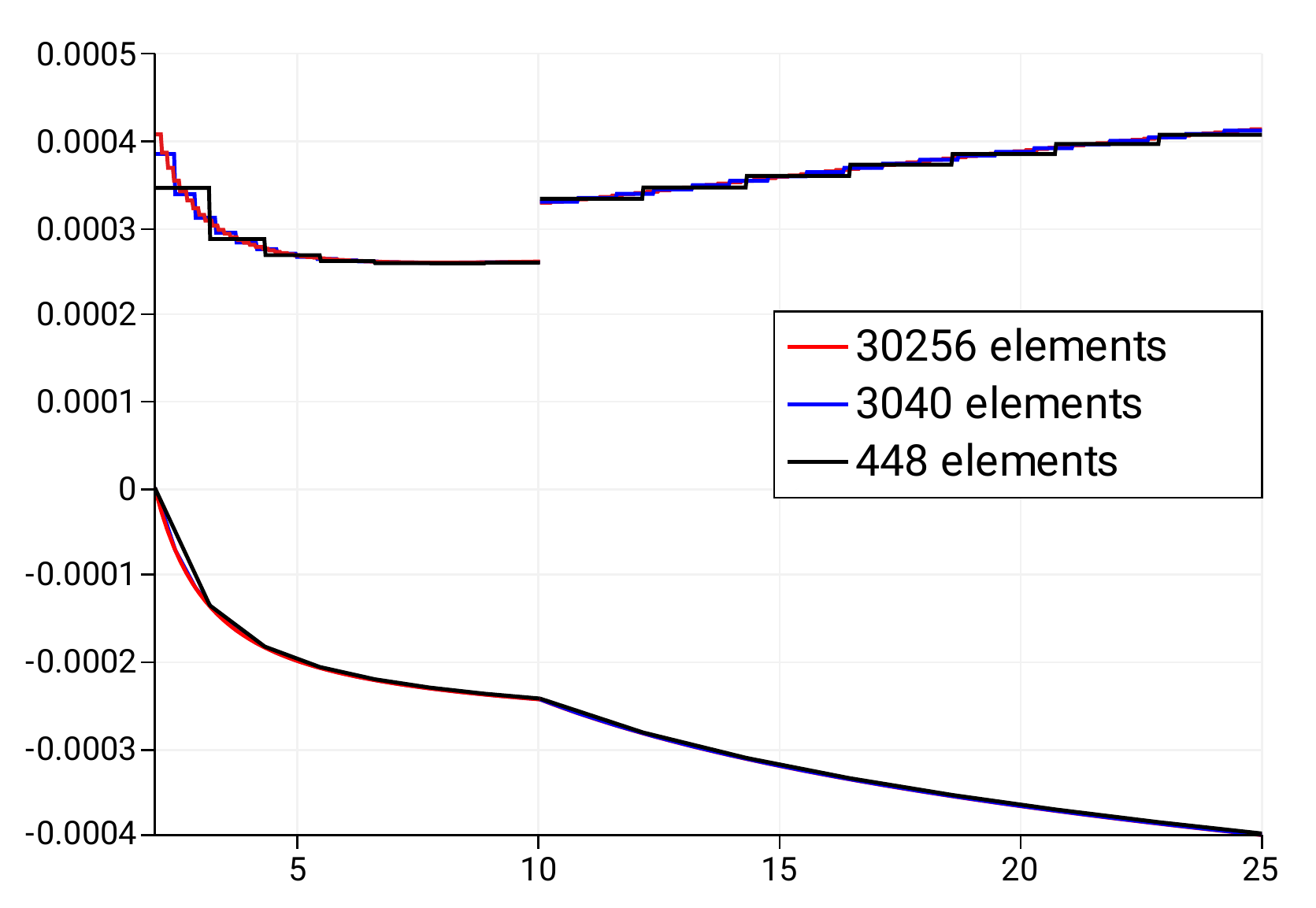}
             \put(-39.903717, -1.2734547){\color[rgb]{0,0,0}\makebox(0,0)[lb]{\smash{$r$}}}%
     \put(-37.903717,12.2734547){\color[rgb]{0,0,0}\makebox(0,0)[lb]{\smash{$\dis_{r \theta}$}}}%
 \put(-37.903717,42.2734547){\color[rgb]{0,0,0}\makebox(0,0)[lb]{\smash{$\dis_{\theta r}$}}}%
	\caption{Q2NQ1}
		  \end{subfigure} 
	  \begin{subfigure}[b]{0.45\textwidth}
      \includegraphics[width=\textwidth]{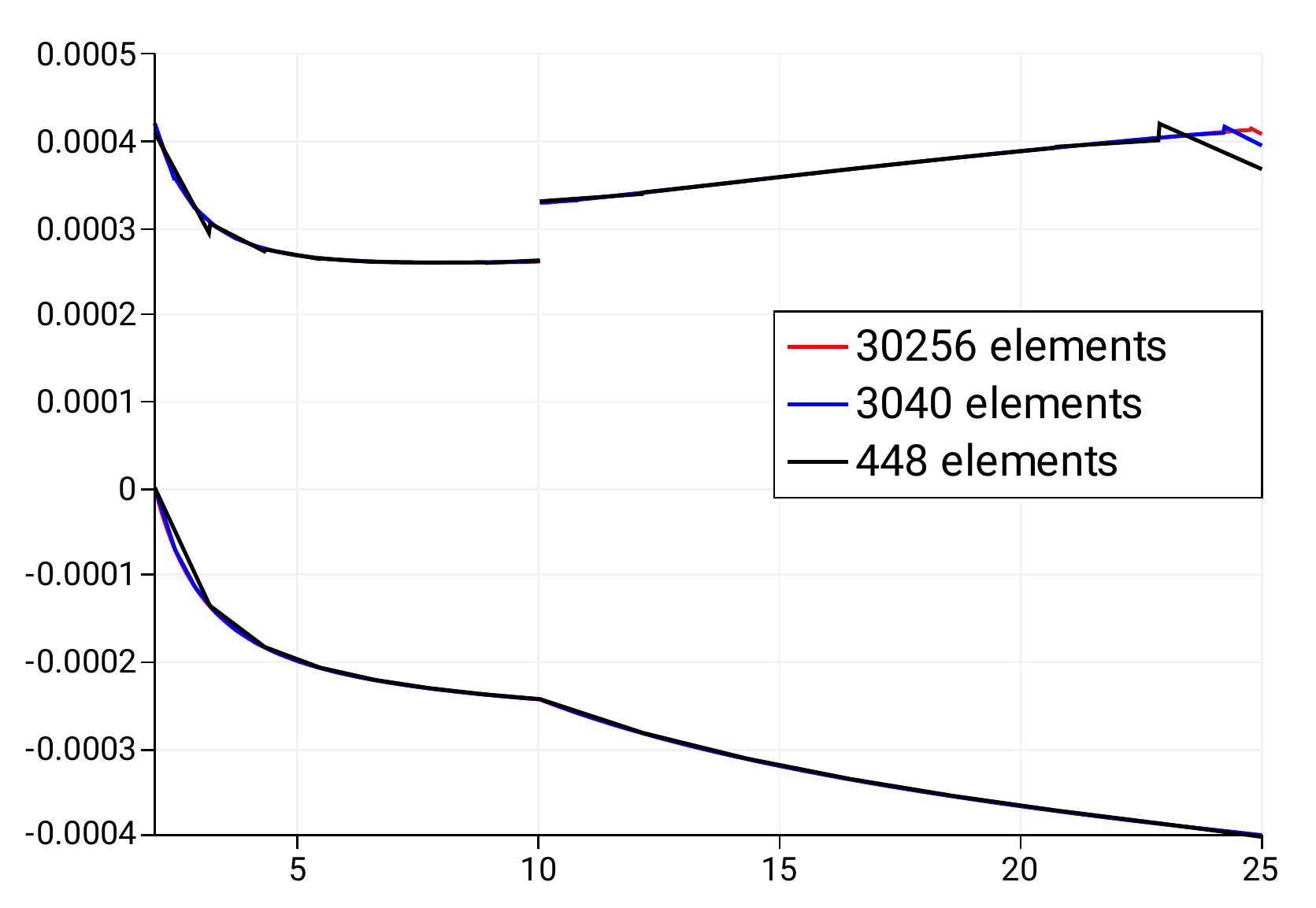}
             \put(-39.903717, -1.2734547){\color[rgb]{0,0,0}\makebox(0,0)[lb]{\smash{$r$}}}%
     \put(-37.903717,12.2734547){\color[rgb]{0,0,0}\makebox(0,0)[lb]{\smash{$\dis_{r \theta}$}}}%
 \put(-37.903717,42.2734547){\color[rgb]{0,0,0}\makebox(0,0)[lb]{\smash{$\dis_{\theta r}$}}}%
	\caption{Q2NQ2}
		  \end{subfigure}
	\caption{The non-vanishing components of $\bP$ along the radius for $H^1(\B) \times H(\curl,\B)$ elements using three mesh densities and $\Lc=5$.} 
	\label{fig:example2:convergence2}
\end{figure}
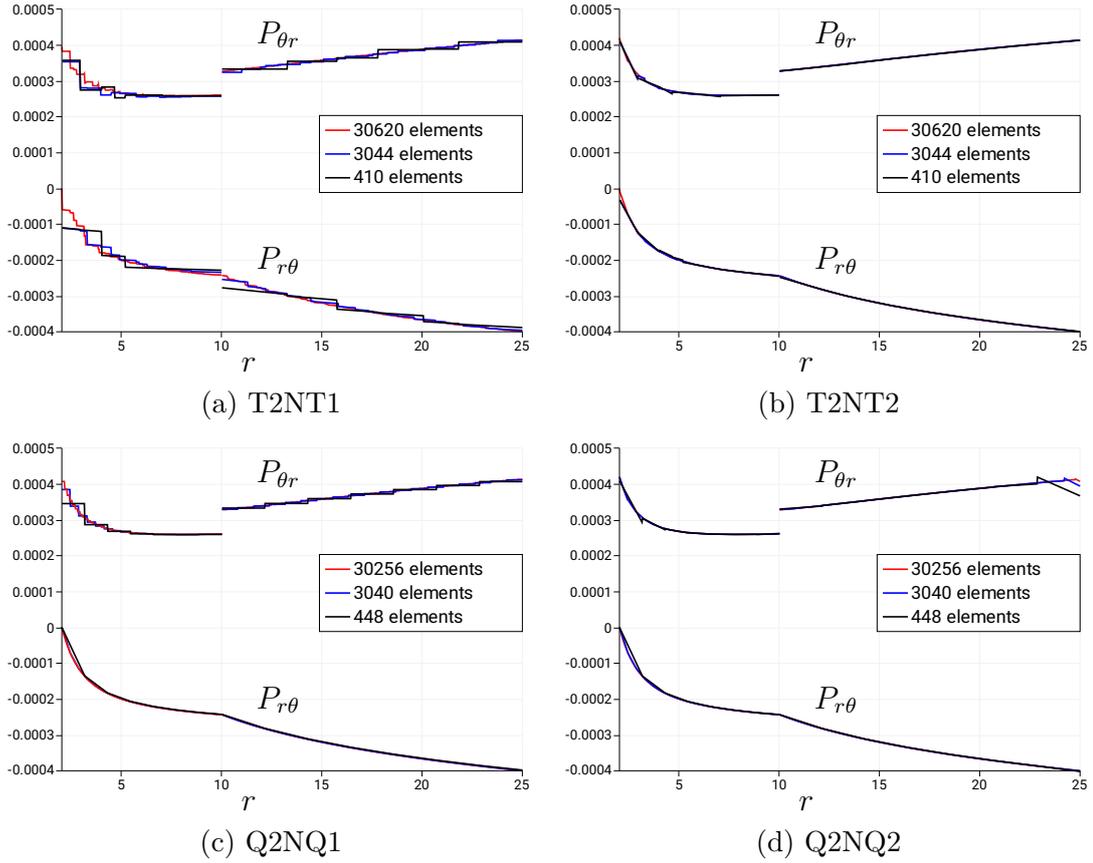

The problem results in a rotationally-symmetric solution where only the shear components ($u_{r,\theta}, u_{\theta,r}, \dis_{r \theta}, \dis_{\theta r} \neq 0$)  are non-vanishing. The convergence behavior of the different elements is investigated for case B and $\Lc = 5$ using three different mesh densities
(410, 3044 and 30620 triangular elements and 448, 3040 and 30256 quadrilateral elements).  Since the micro-distortion field is in $H(\curl,\B)$, the tangential shear component $\dis_{r \theta}$ has to be continuous while the radial shear component $\dis_{\theta r}$ exhibits a jump, see Figure \ref{fig:example2:solution}, where the Q2NQ1 element is used. Similar to Section \ref{bvp:1}, the $H^1(\B) \times H^1(\B)$ elements are unable to capture this discontinuity in $\dis_{\theta r}$, which is shown in Figure \ref{fig:example2:convergence1}. 
 Actually, $H^1(\B) \times H^1(\B)$ elements approximate the discontinuous  $H(\curl,\B)$ solution only when a very fine discretization is used. The discontinuous solution of the micro-distortion field is demonstrated in Figure  \ref{fig:example2:convergence2}  using the  $H^1(\B) \times H(\curl,\B)$ elements. The higher-order N\'ed\'elec formulation in T2NT2 and Q2NQ2 elements exhibit very satisfactory results already with the coarse mesh.

\begin{figure}[ht]
\center
	\unitlength=1mm
  	\begin{picture}(110,60)
	\put(0,0){\def\svgwidth{11 cm}{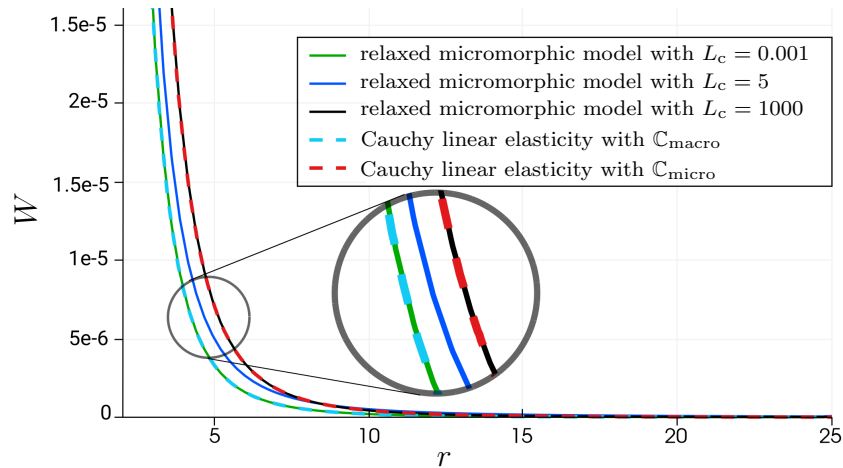}}
	      	   \put(52.903717, -0.2734547){\color[rgb]{0,0,0}\makebox(0,0)[lb]{\smash{$r$}}}%
   \put(0.203717, 30.2734547){\color[rgb]{0,0,0}\makebox(0,0)[lb]{\rotatebox[origin=c]{90}{ \smash{$W$}}}}%
	\end{picture}
\caption{Elastic energy $W$ along the radius.}
\label{fig:example2:energy}
	\end{figure}

\begin{figure}[ht]
\center
	\unitlength=1mm
  	\begin{picture}(110,70)
	\put(-2,0){\def\svgwidth{11.1 cm}{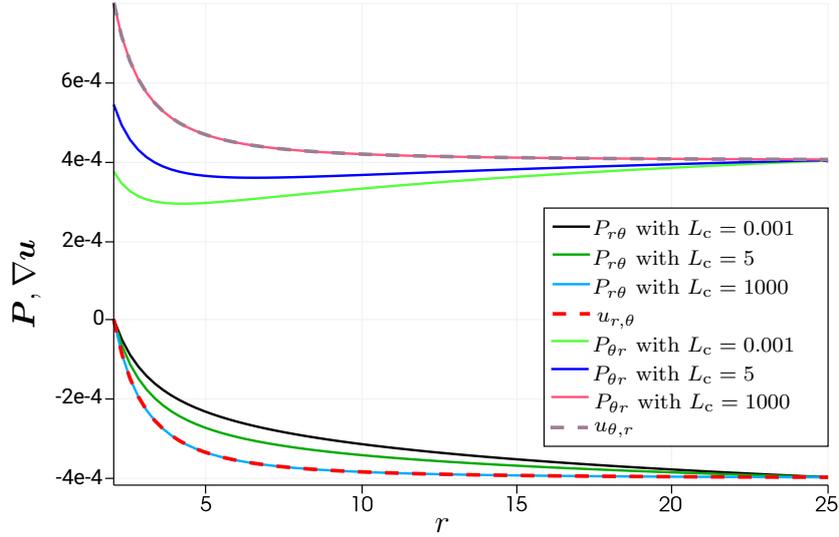}}
	      	   \put(52.903717, -0.2734547){\color[rgb]{0,0,0}\makebox(0,0)[lb]{\smash{$r$}}}%
   \put(0.203717, 25.2734547){\color[rgb]{0,0,0}\makebox(0,0)[lb]{\rotatebox[origin=c]{90}{ \smash{$\Bdis, \nabla \bu$}}}}%
	\end{picture}
	\caption{The non-vanishing components of $\Bdis$ and $\nabla \bu$ along the radius. $\nabla \bu$ is not influenced  by the value of the characteristic length $\Lc$.}
\label{fig:example2:chi}
	\end{figure}

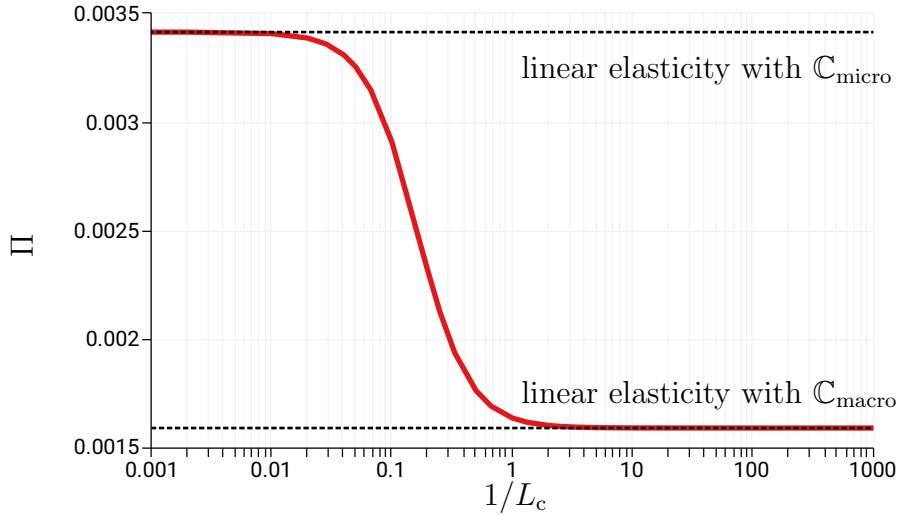
\begin{figure}[ht]
\center
	\unitlength=1mm
         \includegraphics[width=0.7 \textwidth]{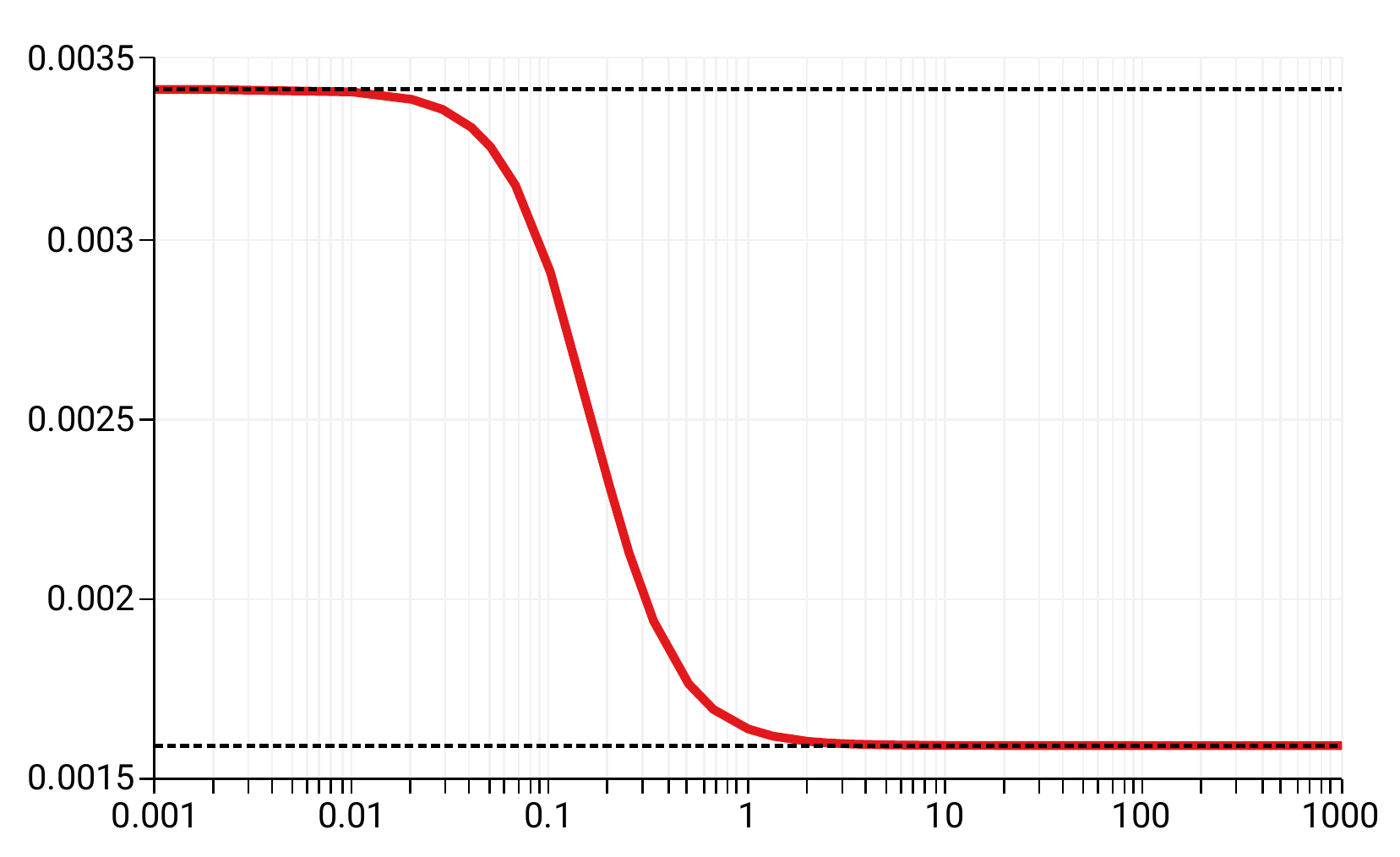} 
 \put(-55.903717, -1.2734547){\color[rgb]{0,0,0}\makebox(0,0)[lb]{\smash{$1/\Lc$}}}%
  \put(-50.903717, 12.2734547){\color[rgb]{0,0,0}\makebox(0,0)[lb]{\smash{linear elasticity with $\Cmacro$}}}%
   \put(-50.903717, 55.2734547){\color[rgb]{0,0,0}\makebox(0,0)[lb]{\smash{linear elasticity with $\Cmicro$}}}%
   \put(-115.203717, 30.2734547){\color[rgb]{0,0,0}\makebox(0,0)[lb]{\rotatebox[origin=c]{90}{ \smash{$\Pi$}}}}%
		  	\caption{Total potential varying the characteristic length.}
\label{fig:example2:potential}
\end{figure}

In the following, we will analyze the influence of a variation of the length-scale parameter $\Lc$ on the response of the relaxed micromorphic model for case A. The relation of the relaxed micromorphic model to the classical Cauchy theory has been discussed in detail in \cite{NefEidMad:2019:ios,BarMadDagAbrGhiNef:2017:taf} for the limiting case $\Lc \rightarrow 0$ and $\Lc \rightarrow \infty$. The case $\Lc \rightarrow 0$ relates to a macroscopic view on the material with microstructure,  with the relaxed micromorphic model being equivalent to a linear elasticity model with stiffness $\Cmacro$ defined as the Reuss lower-bound of $\Ce$ and $\Cmicro$, i.e. $\Cmacro := (\Ce^{-1} + \Cmicro^{-1})^{-1}$.  The case  $\Lc \rightarrow \infty$ resembles an infinite zoom into the material, where an equivalence to linear elasticity with $\Cmicro$ can be derived, cf. \cite{NefEidMad:2019:ios}.  In the latter case, it can be shown that it holds that $\Bdis = \nabla \bu$.  In our study, we approximate the limiting cases by $\Lc=10^{-3}$ and $\Lc=10^{3}$, respectively. Figure \ref{fig:example2:energy} shows the elastic energy $W$ along the radius and Figure \ref{fig:example2:chi} illustrates the non-vanishing components of $\Bdis$ together with the respective displacement gradient components using 30624 Q2NQ2 elements.  In Figure \ref{fig:example2:potential}, we plot the total potential of the relaxed micromorphic model varying  the characteristic length parameter $\Lc$. The figures clearly show the above described behavior. The bounded behavior of the relaxed micromorphic model for small sizes, i.e. $\Lc \rightarrow 0$, is an important advantage which most other generalized continua miss. For the linear elasticity model, a standard T2 nodal element is implemented.

\begin{figure}[ht]
\center
	\unitlength=1mm
  	\begin{picture}(120,70)
	\put(0,1){\def\svgwidth{11 cm}{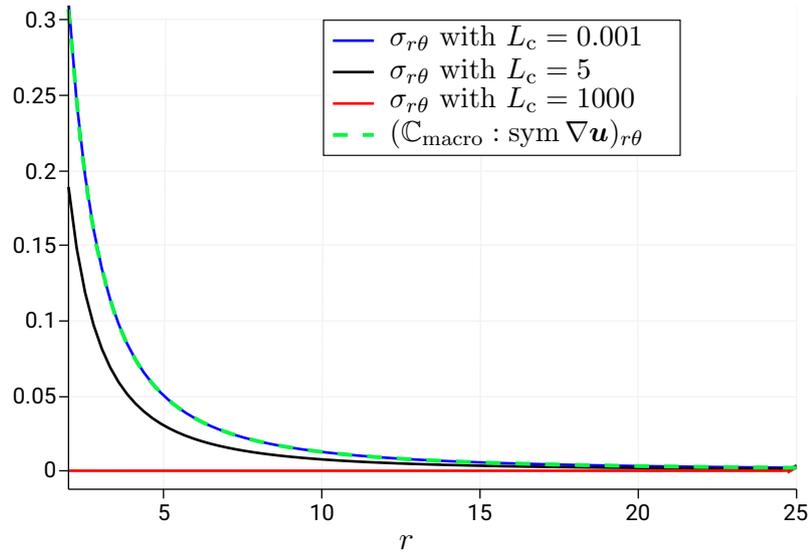}}
	      	   \put(52.903717, -0.2734547){\color[rgb]{0,0,0}\makebox(0,0)[lb]{\smash{$r$}}}%
	\end{picture}
\caption{Force stress shear component $\sigma_{r \theta} = \sigma_{\theta r}$ plotted along the radius.}
\label{fig:example2:cauchystress}
	\end{figure}

\begin{figure}[ht]
\center
	\unitlength=1mm
  	\begin{picture}(120,70)
	\put(0,1){\def\svgwidth{11 cm}{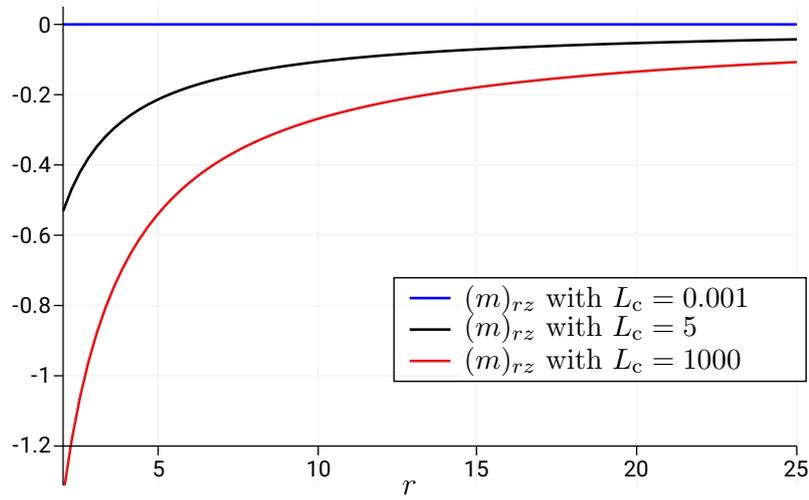}}
	      	   \put(52.903717, -0.2734547){\color[rgb]{0,0,0}\makebox(0,0)[lb]{\smash{$r$}}}%
	\end{picture}
\caption{Non-zero component of the moment stress $m_{r z}$ plotted along the radius.}
\label{fig:example2:momentstress}
	\end{figure}

\begin{figure}[ht]
\center
	\unitlength=1mm
  	\begin{picture}(120,70)
	\put(0,1){\def\svgwidth{11 cm}{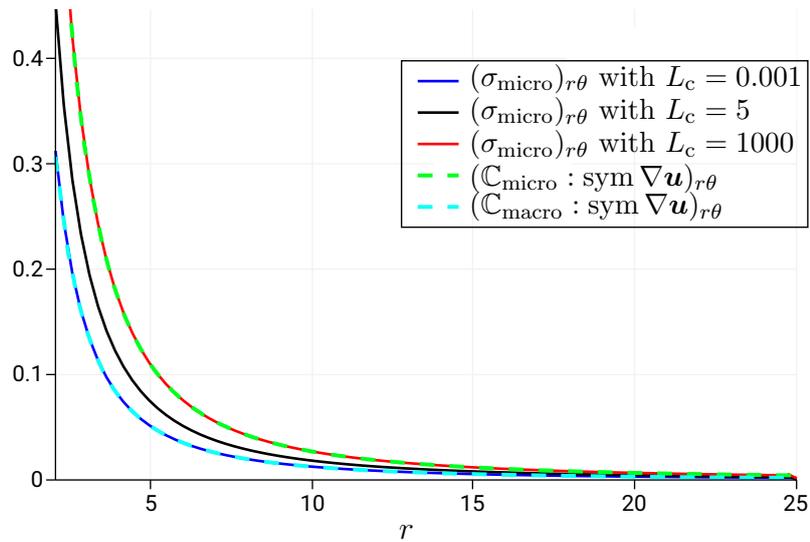}}
	      	   \put(52.903717, -0.2734547){\color[rgb]{0,0,0}\makebox(0,0)[lb]{\smash{$r$}}}%
	\end{picture}
\caption{Micro-stress shear component $(\sigma_\textrm{micro})_{r \theta} = (\sigma_\textrm{micro})_{\theta r}$ plotted along the radius.}
\label{fig:example2:microstress}
	\end{figure}

Next, we investigate the behavior of the different stresses $\Bsigma$, $\Bsigma_\textrm{micro}$ and $\bbm$ under a variation of $\Lc$. The force stress tensor $\Bsigma$ shown in Figure  \ref{fig:example2:cauchystress} vanishes for large value of the characteristic length, $\Lc=1000$, while it is bounded from above by the classical linear elasticity stress with elasticity tensor $\Cmacro$ for $\Lc=0.001$. The only non-vanishing component of the moment stress $m_{r z}$ is shown in Figure \ref{fig:example2:momentstress} ($m_{\theta z} = 0$), which behaves opposite to the force stress when varying $\Lc$. It is nearly zero for $\Lc=0.001$  and it rises for growing $\Lc$. The micro-stress shown in  Figure \ref{fig:example2:microstress} is confined between the linear elasticity stress with elasticity tensor $\Cmicro$ from above and the one with $\Cmacro$  from below for large and small values of the characteristic length, respectively. 
Summarizing the previous findings shortly, increasing the characteristic length diminishes force stress and raises the micro- and  moment stresses, while both force and moment stresses vanish for large and small values of the characteristic length, respectively, the micro-stress is always present.

\FloatBarrier

\sect{\hspace{-5mm}. Conclusions}
\label{sec:con}
The relaxed  micromorphic model is a generalized continuum model which can suitably reproduce the macroscopic effective  properties of mechanical metamaterials.  First,  we derived the variational problem with the relevant weak and strong forms and the associated boundary conditions.  We put together the main components of  standard nodal and nodal-edge finite element formulations of the relaxed micromorphic model.  The standard nodal elements $H^1(\B) \times H^1(\B)$ are incapable to achieve satisfactory results for a discontinuous solution unlike $H^1(\B) \times H(\curl, \B)$ elements which capture the jumps of the normal components of the micro-distortion field and therefore converges   efficiently. We reveal the role of the characteristic length which governs the scale-dependency property of the relaxed micromorphic model. For $\Lc  \rightarrow 0$, the model is equivalent to the standard Cauchy linear elasticity model with $\Cmacro$ defined as the Reuss lower-limit of elasticity tensors $\Ce$ and $\Cmicro$ while the model is corresponding to  Cauchy linear elasticity model with $\Cmicro$ with $\Bdis = \nabla \bu$ for  $\Lc  \rightarrow \infty$.  Furthermore, we have shown the dependency of the different stress measurements  on the characteristic length.  The force stress is at maximum for $\Lc  \rightarrow 0$  and it vanishes for  $\Lc  \rightarrow \infty$ but the moment stress behaves in the opposite way. The micro-stress varies between Cauchy linear elasticity stresses with  $\Cmicro$  and $\Cmacro$ for   $\Lc  \rightarrow \infty$ and  $\Lc  \rightarrow 0$, respectively. 
 
\vspace{2 cm}

{\bf Acknowledgment} \\
Funded by the Deutsche Forschungsgemeinschaft (DFG, German research Foundation) -  Project number 440935806 (SCHR 570/39-1, SCHE 2134/1-1, NE 902/10-1) within the DFG priority program 2256. The authors gratefully acknowledge Jo\v{z}e Korelc for the development and ongoing support when using AceGen and AceFEM.



\bibliographystyle{plainnat}
{\footnotesize
\bibliography{micmag_01}
}
\begin{appendix}
\sect{Appendix}
\label{app:Nedelec_shape_functions}
\numberwithin{equation}{section}
\setcounter{equation}{0}
\ssect{Construction of triangular N\'ed\'elec shape functions}
The parameter elements are defined in the natural coordinates $\Bxi = (\xi,\eta)^T$. Triangular elements are defined on the domain $\B_e^\triangle = \{0 \le \xi \le 1, 0 \le \eta \le 1 - \xi\}$.  The finite elements with edge numbering are shown in Figure  \ref{Figure:nedelec_elements}.

\begin{figure}[ht]
     \begin{subfigure}[b]{0.32\textwidth}
         \centering
         \includegraphics[width=\textwidth]{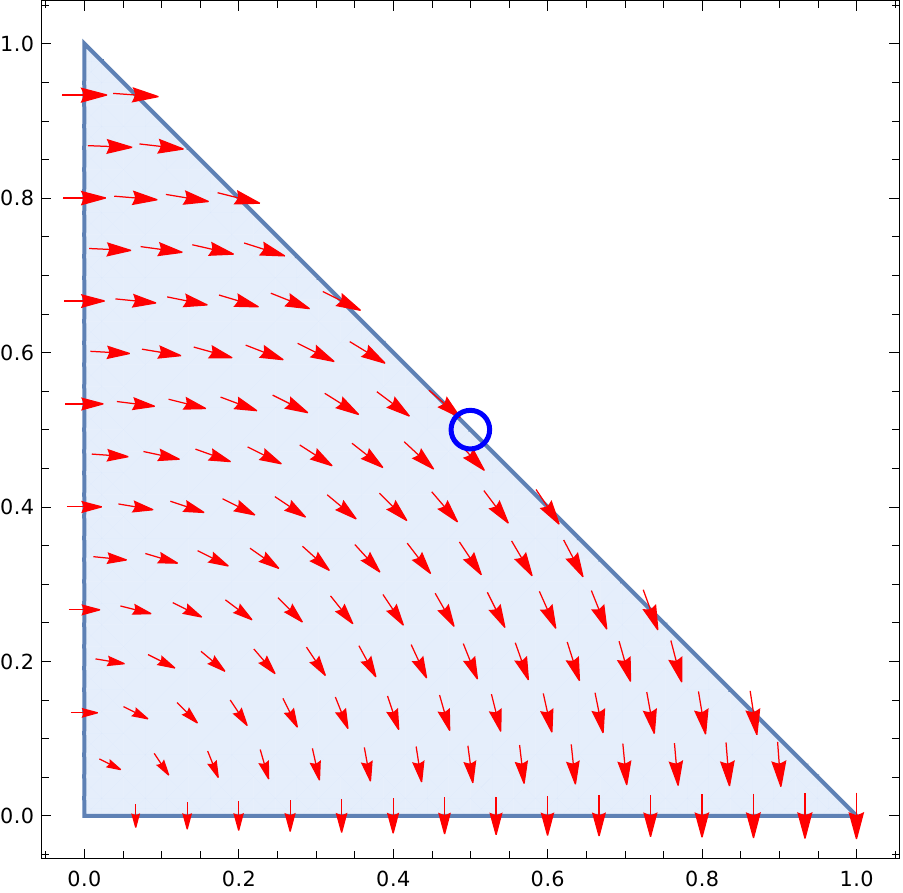}
        \caption{$\bv^1_1$}
     \end{subfigure}
     \begin{subfigure}[b]{0.32\textwidth}
         \centering
         \includegraphics[width=\textwidth]{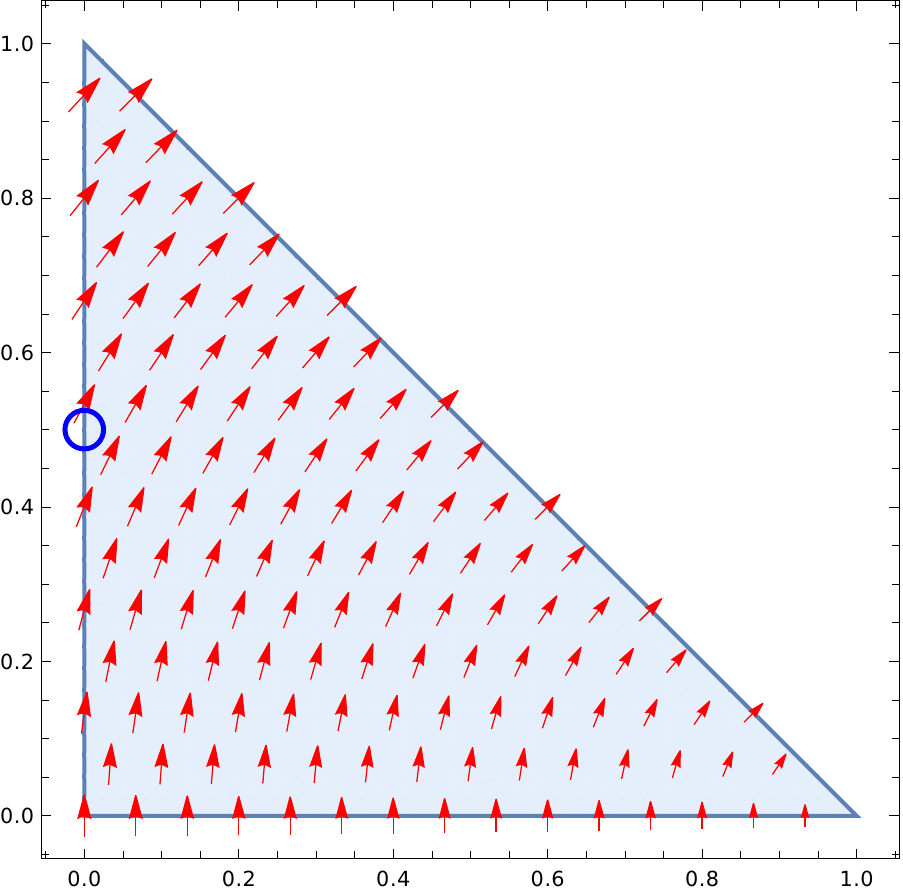}
  \caption{$\bv^1_2$}
     \end{subfigure}
          \begin{subfigure}[b]{0.32\textwidth}
         \centering
         \includegraphics[width=\textwidth]{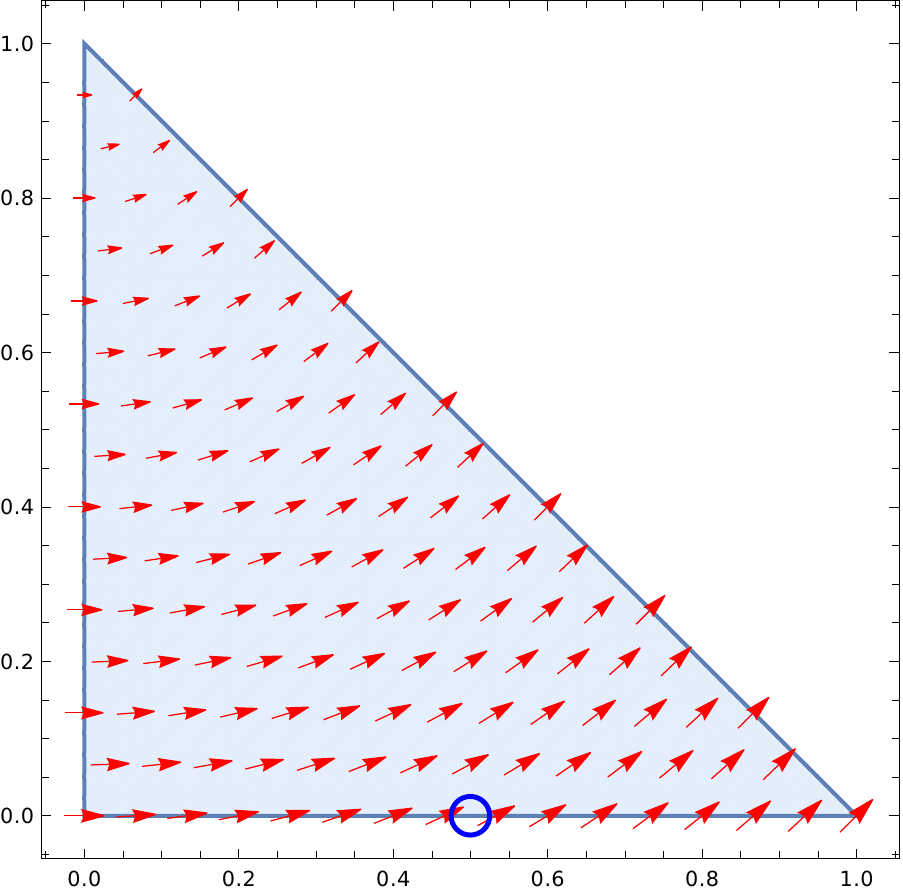}
  \caption{$\bv^1_3$}
     \end{subfigure}
        \caption{Tangential-conforming vectorial shape functions of NT1 element. Blue circles indicate the position where the dofs are defined. }
        \label{Fig:shape_function_NT1}
\end{figure}

\sssect{\hspace{-5mm}. First-order triangular element NT1 \\} 
\label{sec:app:nt1}
The  N\'ed\'elec space of this element reads 
 \begin{equation}
  \left[ \mathcal{ND}^\triangle \right]^{2}_1 = \bigg\{ \left[ \begin{array}{c}
1 \\
0  
  \end{array} \right] \,, \left[ \begin{array}{c}
0 \\
1
  \end{array} \right]  \,, \left[ \begin{array}{c}
- \eta \\
\xi
  \end{array} \right]  \bigg\}\,,
 \end{equation}
 and the general form of the vectorial shape function is  
 \begin{equation}
 \label{eq:app:gf1}
 \bv^1  = \left[\begin{array}{c}
a_1 - a_3 \, \eta \\ 
a_2 + a_3 \, \xi
\end{array}\right],
\end{equation}
where $a_i, i=1,2,3$ are coefficients yet to be defined based on the dofs. Starting from definition in Equation (\ref{eq:t:edge_dof}), we set $r_j =1$ for all edges and the tangential vectors of the edges, see Figure \ref{Figure:nedelec_elements} (right), are 
\begin{equation}
\bt_1  = \frac{1}{\sqrt{2}}\left[ \begin{array}{c}
1 \\
-1 \\
\end{array} \right], \quad \bt_2  = \left[ \begin{array}{c}
0 \\
1 \\
\end{array} \right], \quad \bt_3  =
\left[ \begin{array}{c}
1 \\
0 \\
\end{array} \right].
\end{equation}
We calculate the dofs  following Equation (\ref{eq:t:edge_dof}) using $\xi + \eta = 1$ on the first edge, $\xi=0$ on the second edge and $\eta=0$ on the third edge and obtain
\begin{equation}
 m^{e_1}_{1} = a_1- a_2 - a_3, \qquad m^{e_2}_{1} = a_2, \qquad  m^{e_3}_{1} = a_1\,.
\end{equation}
In order to obtain the three vectorial shape functions $ \bv^1_1,  \bv^1_2$ and $ \bv^1_3$ form the general function in  (\ref{eq:app:gf1}), we have to compute the three associated combinations for $a_1, a_2$ and $a_3$. We derive the explicit expressions for the three vectorial shape functions by enforcing for the function $\bv_j^k$ associated to the edge $e_j$ 
\begin{equation}
 m^{e_i}_{1} (\bv_j^k) = \delta_{ij}\,.
\end{equation}
 The evaluation of the dofs for each edge leads with
\begin{equation}
\begin{aligned}
\textrm{edge 1:} \quad &m^{e_1}_{1} = 1, \qquad m^{e_2}_{1} = 0, \qquad  m^{e_3}_{1} = 0 \quad \Rightarrow \quad a_1=0, \quad a_2=0, \quad a_3=-1 \\ 
\textrm{edge 2:}  \quad  &m^{e_1}_{1} = 0, \qquad m^{e_2}_{1} = 1, \qquad  m^{e_3}_{1} = 0 \quad \Rightarrow \quad a_1=0, \quad a_2=1, \quad a_3=-1 \\ 
\textrm{edge 3:} \quad   &m^{e_1}_{1} = 0, \qquad m^{e_2}_{1} = 0, \qquad  m^{e_3}_{1} = 1 \quad \Rightarrow \quad a_1=1, \quad a_2=0, \quad a_3=\,\,\,\,1 \\ 
 \end{aligned}
\end{equation}

to the shape vectors 
\begin{equation}
\bv^1_1 =  \left( \begin{array}{c}
\eta \\ 
- \xi
\end{array} \right) \,, \quad
 \bv^1_2 =  \left( \begin{array}{c}
\eta \\
 1 - \xi
\end{array} \right)\,, \quad
\bv^1_3 =  \left( \begin{array}{c}
1 - \eta \\ 
\xi
\end{array} \right) \,.
\end{equation}
A visualization of them is depicted in Figure \ref{Fig:shape_function_NT1}.

\begin{figure}[ht]
     \begin{subfigure}[b]{0.32\textwidth}
         \centering
         \includegraphics[width=\textwidth]{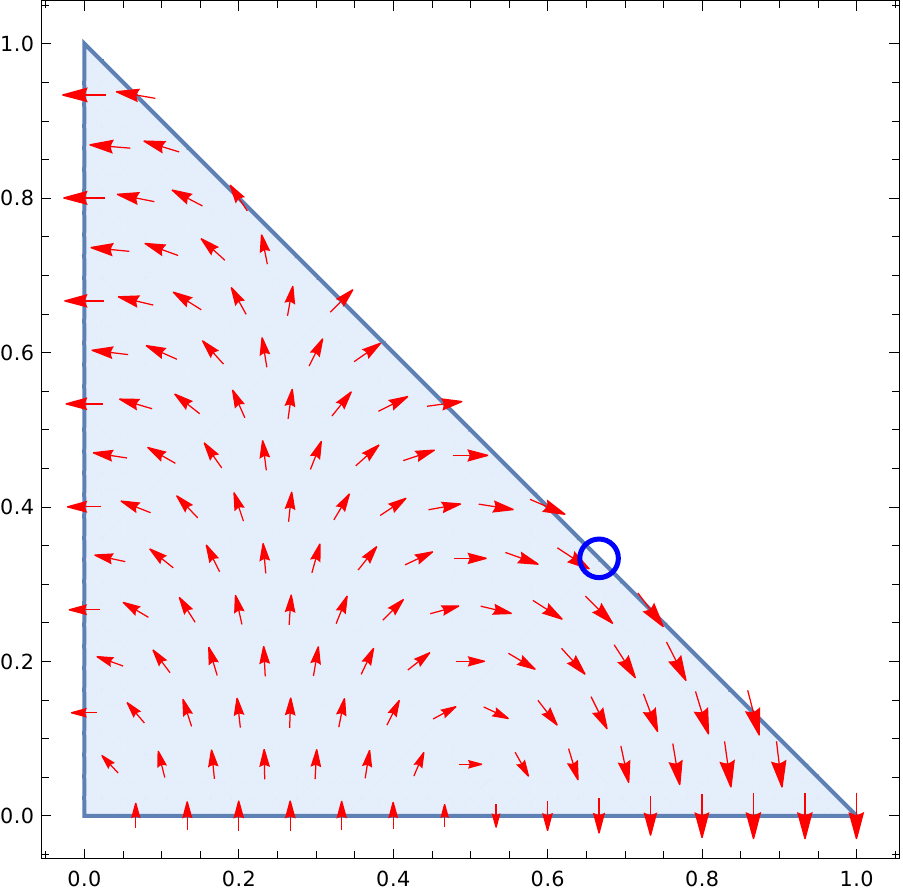}
        \caption{$\bv	^2_1$}
     \end{subfigure}
     \begin{subfigure}[b]{0.32\textwidth}
         \centering
         \includegraphics[width=\textwidth]{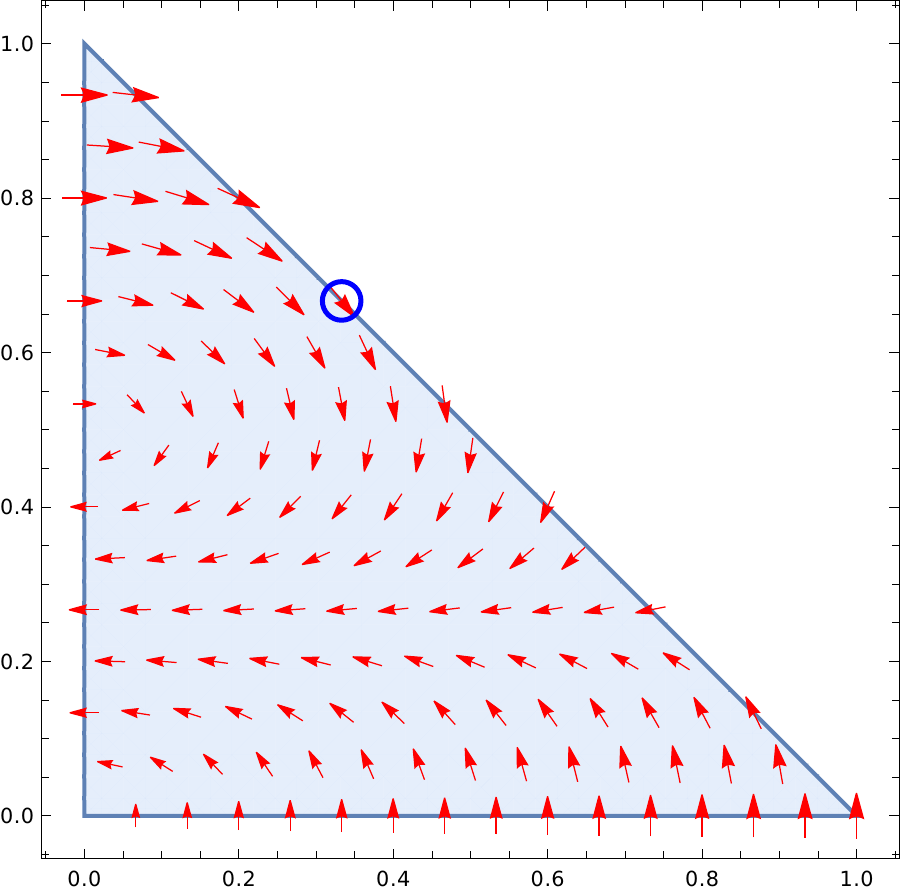}
  \caption{$\bv^2_2$}
     \end{subfigure}
          \begin{subfigure}[b]{0.32\textwidth}
         \centering
         \includegraphics[width=\textwidth]{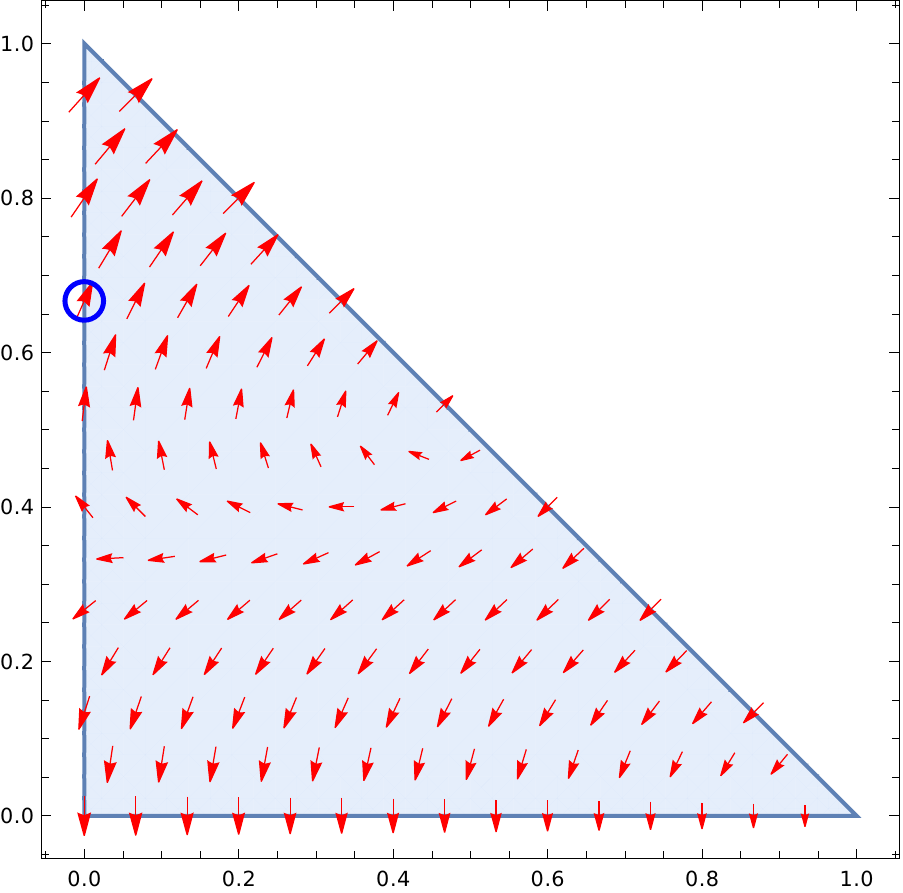}
  \caption{$\bv^2_3$}
     \end{subfigure}
          \begin{subfigure}[b]{0.32\textwidth}
         \centering
         \includegraphics[width=\textwidth]{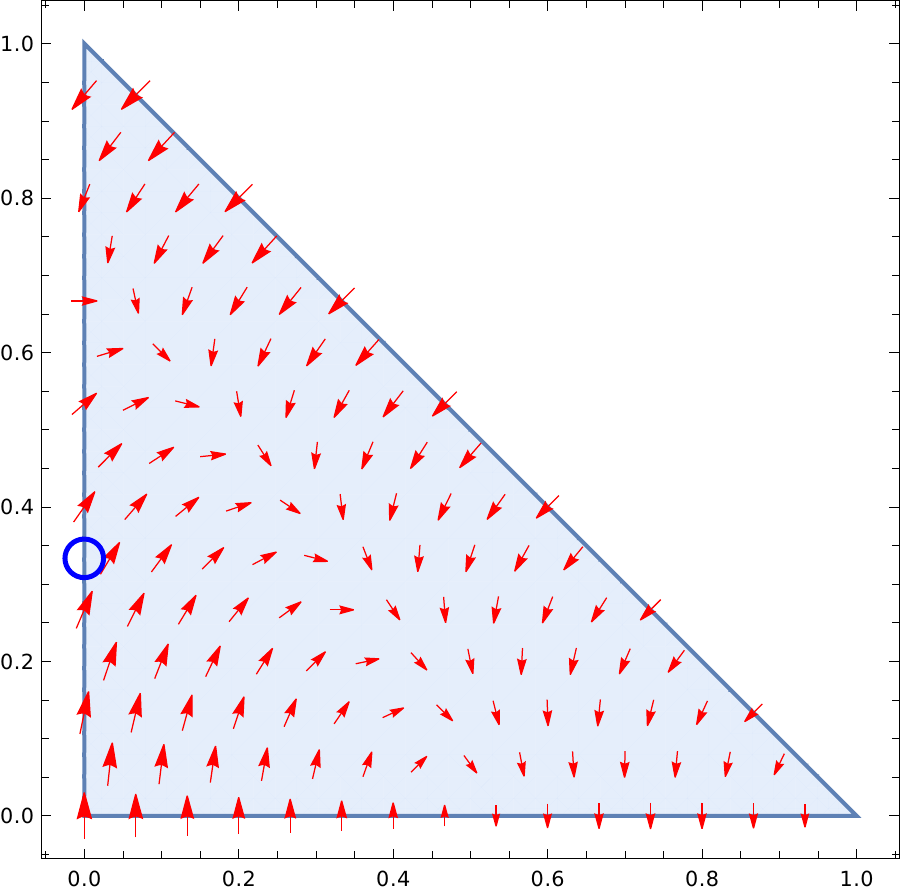}
        \caption{$\bv^2_4$}
     \end{subfigure}
     \begin{subfigure}[b]{0.32\textwidth}
         \centering
         \includegraphics[width=\textwidth]{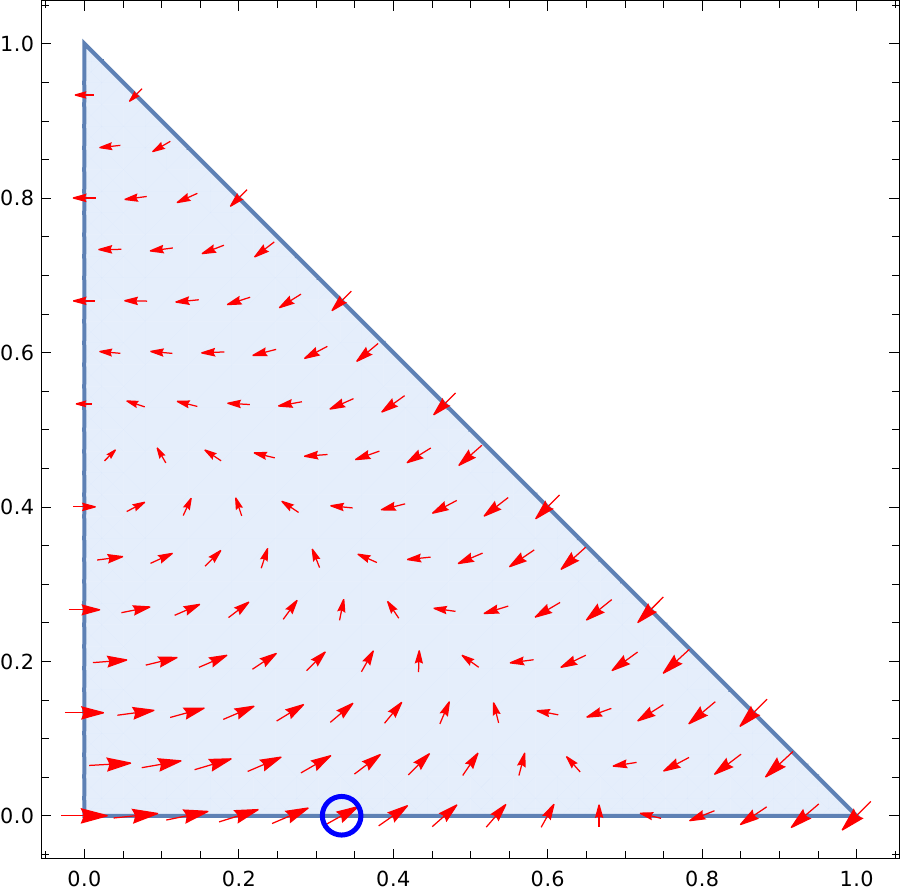}
  \caption{$\bv^2_5$}
     \end{subfigure}
          \begin{subfigure}[b]{0.32\textwidth}
         \centering
         \includegraphics[width=\textwidth]{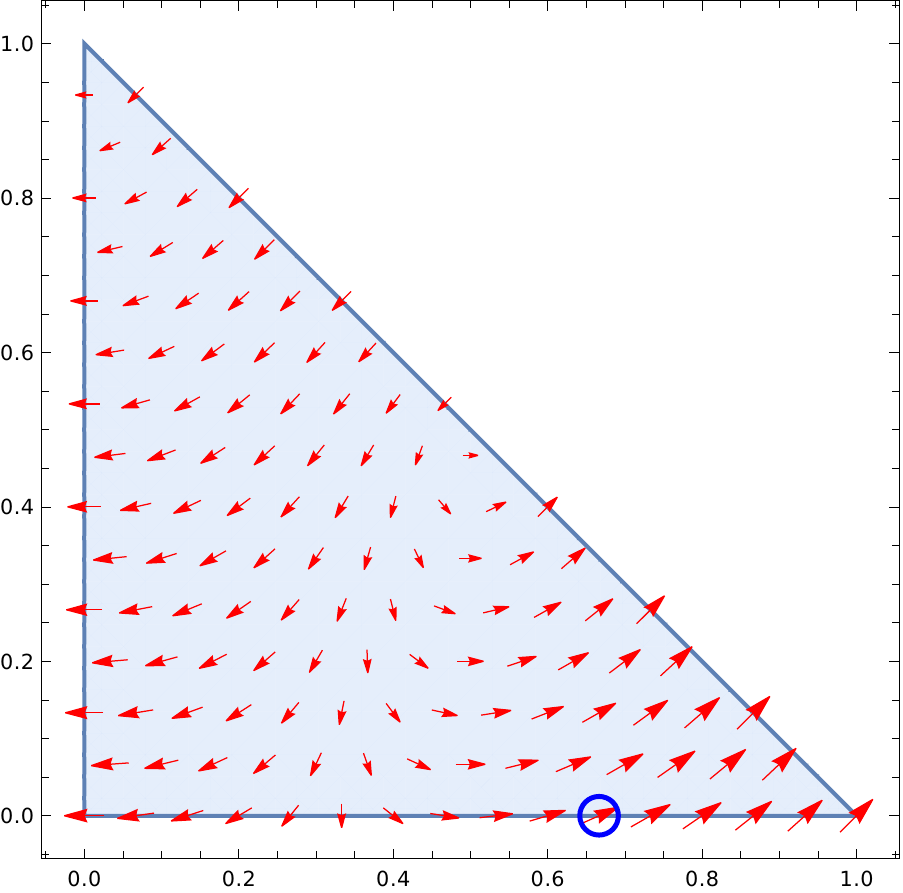}
  \caption{$\bv^2_6$}
     \end{subfigure}
          \begin{subfigure}[b]{0.32\textwidth}
         \centering
         \includegraphics[width=\textwidth]{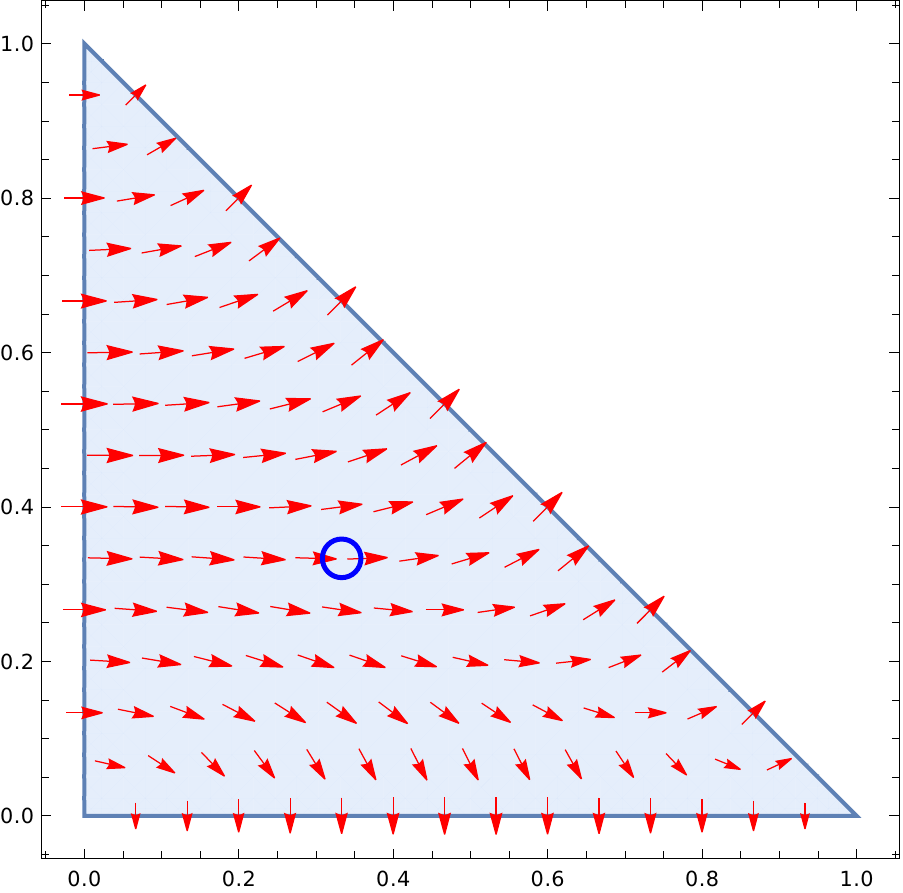}
        \caption{$\bv^2_7$}
     \end{subfigure}
  \begin{subfigure}[b]{0.32\textwidth}
         \centering
         \includegraphics[width=\textwidth]{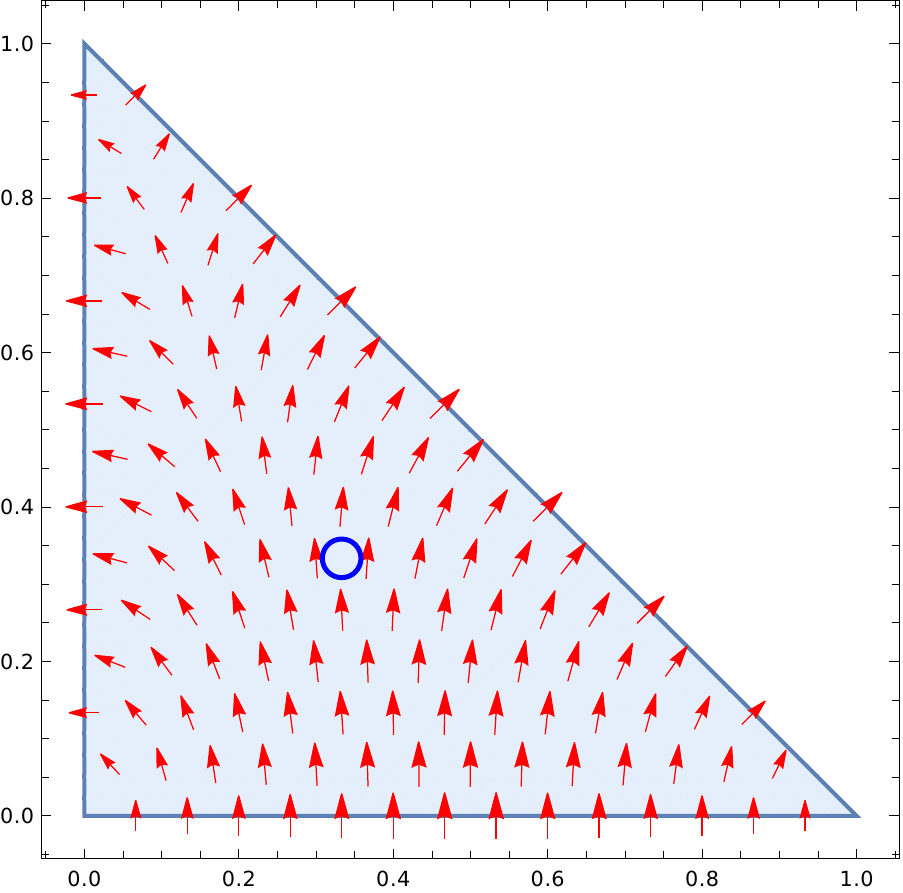}
        \caption{$\bv^2_8$}
     \end{subfigure}
       \caption{Tangential-conforming vectorial shape functions of NT2 element. Blue circles indicate the position where the dofs are defined. }
        \label{Fig:shape_function_NT2}
\end{figure}

\sssect{\hspace{-5mm}. Second-order triangular element NT2 \\} 
The  N\'ed\'elec space of this element reads 
 \begin{equation}
  \left[ \mathcal{ND}^\triangle \right]^{2}_2 = \bigg\{ \left[ \begin{array}{c}
1 \\
0  
  \end{array} \right] \,, \left[ \begin{array}{c}
 \xi \\
0
  \end{array} \right]  \,,  \left[ \begin{array}{c}
 \eta \\
0  \end{array} \right] \,, \left[ \begin{array}{c}
0 \\
1  
  \end{array} \right] \,, \left[ \begin{array}{c}
0 \\
\xi
  \end{array} \right]  \,,  \left[ \begin{array}{c}
0 \\
\eta  \end{array} \right]   \,,   \left[ \begin{array}{c}
- \eta^2 \\
\xi \eta
  \end{array} \right]   \,,   \left[ \begin{array}{c}
- \xi \eta \\
\xi^2
  \end{array} \right] \bigg\}\,,
 \end{equation}
and the general form of the shape functions reads
\begin{equation}
\bv^2 =  \left( \begin{array}{c}
a_1 + a_2 \, \xi + a_3 \, \eta - a_7 \, \eta^2 - a_8 \, \xi\eta \\
a_4 + a_5 \, \xi + a_6 \, \eta + a_7 \, \xi\eta + a_8 \, \xi^2
\end{array} \right) ,
\end{equation}
where $a_i, i=1,...,8$ are coefficients yet to be defined based on the dofs. 
The explicit functions $r_j$ and $\bq_i$ in Equations (\ref{eq:t:edge_dof}) and (\ref{eq:t:inner_dof_t}) are assumed as  
\begin{equation}
\begin{aligned}
&\textrm{edge 1 :} \qquad  &r_1 =\xi\,,&   \qquad   &r_2 = \eta\,,  \\
&\textrm{edge 2 :} \qquad  &r_1 =  \eta\,,& \qquad &r_2 = 1- \eta\,, \\
&\textrm{edge 3 :} \qquad  &r_1 = 1-\xi\,,& \qquad &r_2 = \xi\,, \\  
&\textrm{inner  \, :} \qquad   &\bq_1 =    \left[ \begin{array}{c}
1 \\
0 
  \end{array} \right],& \qquad &\bq_2 = \left[ \begin{array}{c}
0 \\
1  
  \end{array} \right], 
  \end{aligned}
\end{equation}
and the tangential vectors of the edges are same as in NT1 element. \\ 
The inner and outer dofs are calculated according to Equations
(\ref{eq:t:edge_dof}) and (\ref{eq:t:inner_dof_t}) using  $\xi + \eta = 1$ on the first edge, $\xi=0$ on the second edge and $\eta=0$ on the third edge
 \begin{equation}
 \begin{aligned}
 &m^{e_1}_{1} = \frac{1}{6}(3 a_1 + 2 a_2 + a_3 - 3 a_4 - 2 a_5 - a_6 - a_7 - 2 a_8) \,, \\ 
 &m^{e_1}_{2} =\frac{1}{6}  (3 a_1 + a_2 + 2 a_3 - 3 a_4 - a_5 - 2 a_6 - 2 a_7 - a_8) \,,\\ 
  &m^{e_2}_{1} = \frac{1}{6} (3a_4 + 2 a_6) \,, \qquad
  m^{e_2}_{2} = \frac{1}{6} ( 3a_4 + a_6) \,, \\
& m^{e_3}_{1} = \frac{1}{6} (3a_1 + a_2) \,,  \qquad \,\,\, m^{e_3}_{2} = \frac{1}{6} (3a_1 + 2a_2) \,, \\ 
  &m^\textrm{inner}_2 = \frac{1}{24} (12 a_1 + 4 a_2 + 4 a_3 - 2 a_7 - a_8)  \,, \\
 &m^\textrm{inner}_1 =  \frac{1}{24} (12 a_4 + 4 a_5 + 4 a_6 + a_7 + 2 a_8) \,, 
 \end{aligned}
\end{equation}
and the resulting shape functions shown in Figure \ref{Fig:shape_function_NT2} are obtained by an analogous procedure as before
\begin{equation}
\begin{aligned}
&\textrm{edge 1 :} &\bv_1^2&= 2 \left( \begin{array}{c} 
-\eta + 4 \eta \xi \\
 2 \xi - 4 \xi^2 
 \end{array} \right),& \quad 
  &\bv_2^2=2 \left( \begin{array}{c} 
-2 \eta + 4 \eta^2 \\
\xi - 4 \eta \xi
 \end{array} \right),& \\
&\textrm{edge 2 :}  &\bv_3^2 &= 2\left( \begin{array}{c} 
-2 \eta + 4 \eta^2 \\
 -1 + 3 \eta + \xi - 4 \eta \xi
 \end{array} \right),& \quad 
 &\bv_4^2 = 2\left( \begin{array}{c} 
3 \eta - 4 \eta^2 - 4 \eta \xi \\
 2 - 3 \eta - 6 \xi + 4 \eta \xi + 4 \xi^2
 \end{array} \right),& \\ 
&\textrm{edge 3 :}  &\bv_5^2 &= 2\left( \begin{array}{c} 
2 - 6 \eta + 4 \eta^2 - 3 \xi + 4 \eta \xi \\
 3 \xi - 4 \eta \xi - 4 \xi^2
 \end{array} \right),& \quad
 &\bv_6^2 = 2\left( \begin{array}{c} 
-1 + \eta + 3 \xi - 4 \eta \xi \\
 -2 \xi + 4 \xi^2
 \end{array} \right),& \\ 
&\textrm{inner \, :}  &\bv_7^2 &= 2\left( \begin{array}{c} 
8 \eta - 8 \eta^2 - 4 \eta \xi \\
-4 \xi + 8 \eta \xi +  4 \xi^2
 \end{array} \right),& \quad
  &\bv_8^2 = 2\left( \begin{array}{c} 
-4 \eta + 4 \eta^2 + 8 \eta \xi \\
8 \xi - 4 \eta \xi - 8 \xi^2
 \end{array} \right).&
 \end{aligned}
\end{equation}

\ssect{Construction of quadrilateral N\'ed\'elec shape functions}
 Quadrilateral elements have the domain $\B_e^\square = \{-1 \le \xi \le 1, -1 \le \eta \le 1\}$. 

\begin{figure}[ht]
    \centering
     \begin{subfigure}[b]{0.33\textwidth}
         \centering
         \includegraphics[width=\textwidth]{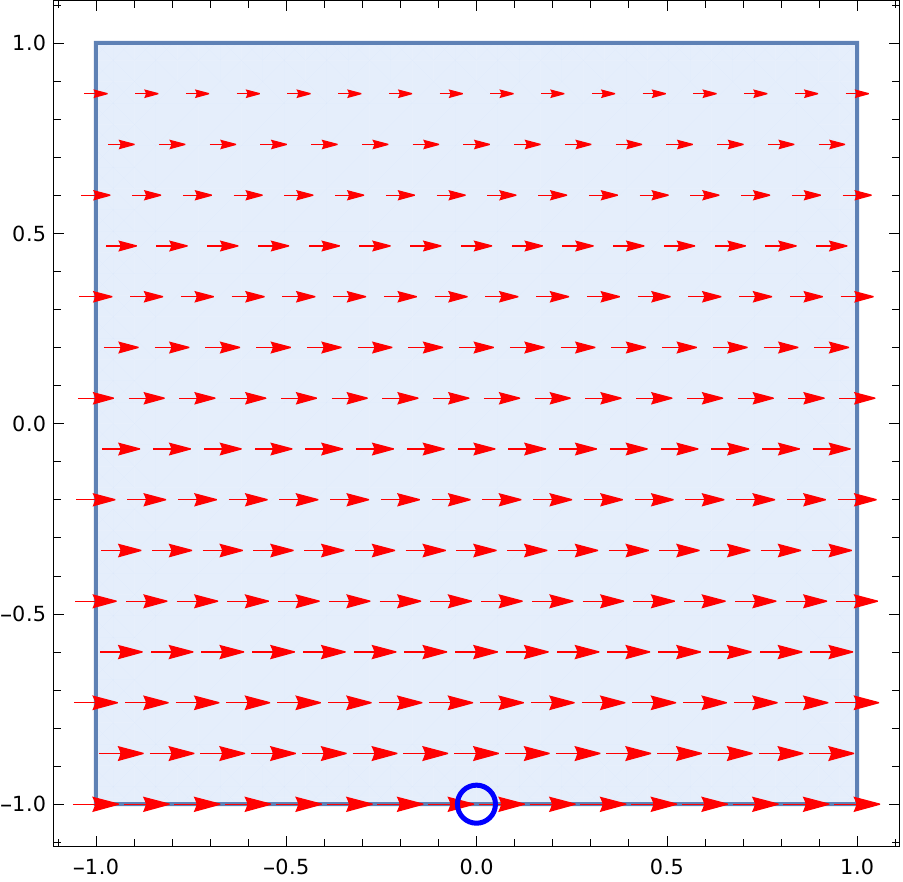}
        \caption{$\bv_1^1$}
     \end{subfigure}
     \begin{subfigure}[b]{0.33\textwidth}
         \centering
         \includegraphics[width=\textwidth]{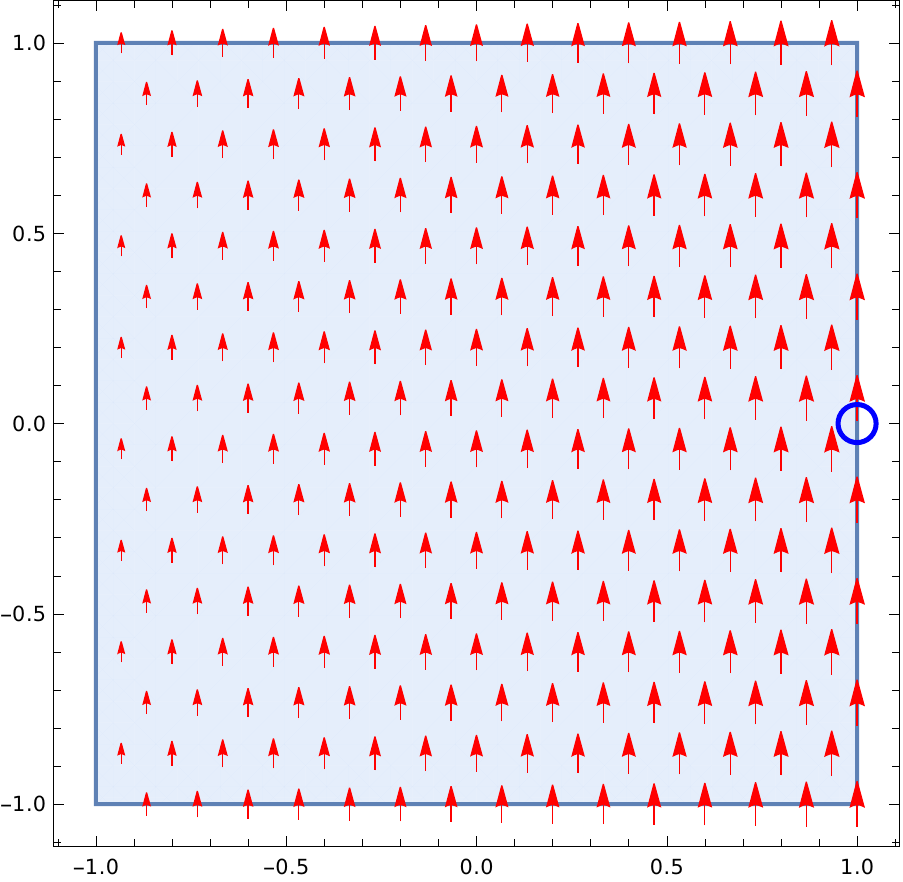}
  \caption{$\bv_2^1$}
     \end{subfigure}
          \begin{subfigure}[b]{0.33\textwidth}
         \centering
         \includegraphics[width=\textwidth]{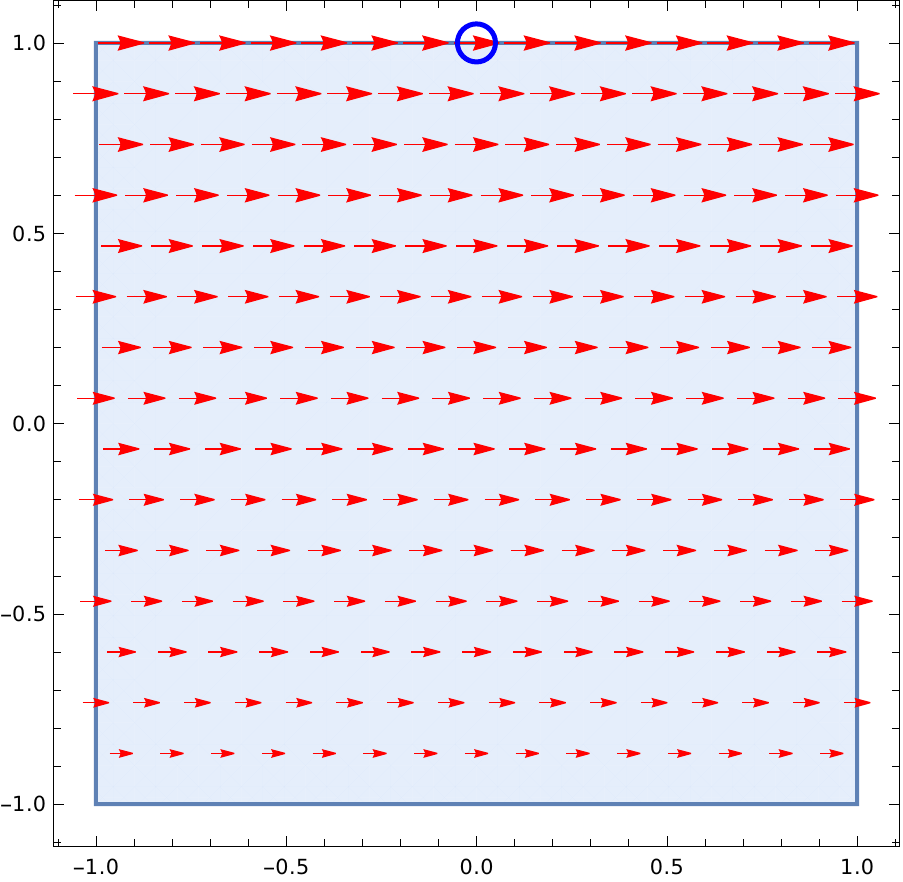}
  \caption{$\bv_3^1$}
     \end{subfigure}
          \begin{subfigure}[b]{0.33\textwidth}
         \centering
         \includegraphics[width=\textwidth]{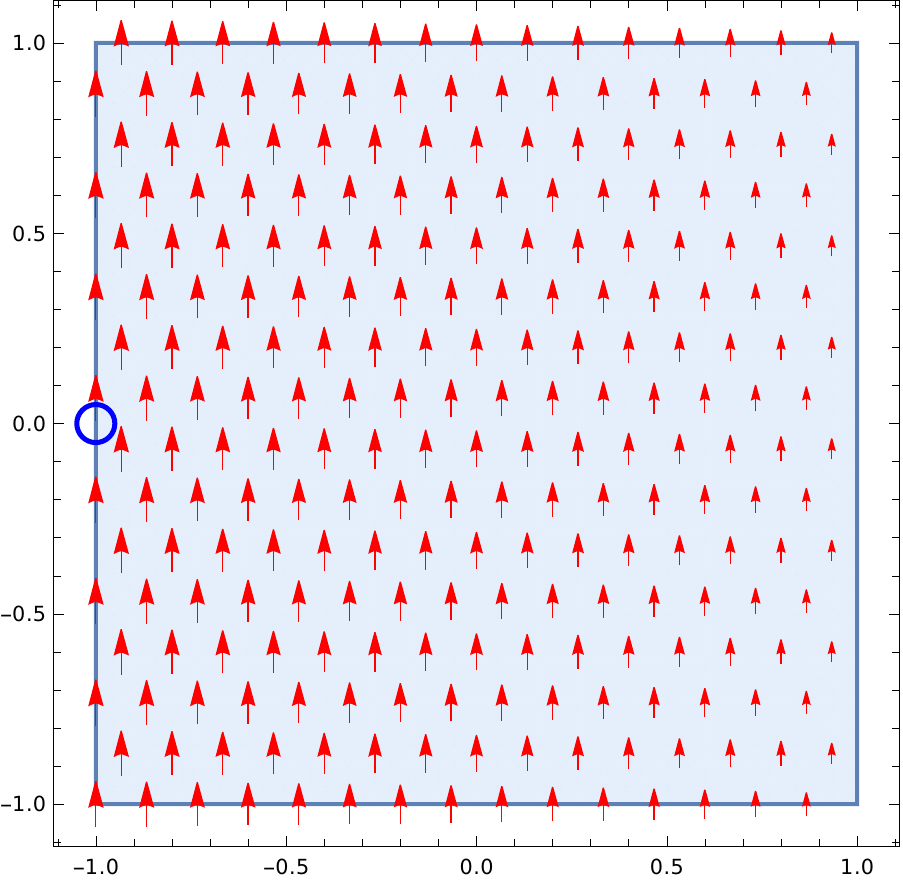}
        \caption{$\bv_4^1$}
     \end{subfigure}
       \caption{Tangential-conforming vectorial shape functions of NQ1 element. Blue circles indicate the position where the dofs are defined. }
        \label{Fig:shape_function_NQ1}
\end{figure}

\sssect{\hspace{-5mm}. First-order quadrilateral element NQ1 \\} 
The  N\'ed\'elec space of this element reads 
 \begin{equation}
  \left[ \mathcal{ND}^\square \right]^{2}_1 = \bigg\{ \left[ \begin{array}{c}
1 \\
0  
  \end{array} \right] \,, \left[ \begin{array}{c}
\eta \\
0
  \end{array} \right]  \,, \left[ \begin{array}{c}
0 \\
1
  \end{array} \right]  \,, \left[ \begin{array}{c}
0 \\
\xi
  \end{array} \right]\bigg\} \,,
 \end{equation}
and the general form of the shape vectors reads 
\begin{equation}
\bv^1 =  \left( \begin{array}{c}
a_1 + a_2 \, \eta  \\ a_3 + a_4 \, \xi
\end{array} \right),
\end{equation}
where $a_i, i=1,..,4$ are coefficients yet to be defined based on the dofs. Starting from definition in Equation (\ref{eq:t:edge_dof}), we set $r_j = 1$ for all edges. The tangential vectors for the first and third edges are $\bt_1 = \bt_3 = (1,0)^T$ and for the second and fourth edges are $\bt_2 = \bt_4 = (0,1)^T$, see Figure \ref{Figure:nedelec_elements} (left).  We calculate the edge dofs taking into consideration $\eta = -1$ on the first edge, $\xi = 1$ on the second edge, $\eta = 1$ on the third edge and $\xi = -1$ on the fourth edge
\begin{equation}
m^{e_1}_{1} = 2(a_1 - a_2)\,, \quad
m^{e_2}_{1} = 2(a_3 + a_4)\,, \quad
m^{e_3}_{1} = 2(a_1 + a_2)\,, \quad
m^{e_4}_{1} = 2(a_3 - a_4) \,.
\end{equation}
We solve the system of equations obtained by an  analogous procedure to Section \ref{sec:app:nt1} leading to the following shape functions demonstrated in Figure \ref{Fig:shape_function_NQ1} where $\bv^1_i$ is associated with the edge $e_i$ for $i=1,..,4$
\begin{equation}
\begin{aligned}
\bv_{1}^1 =&  \left( \begin{array}{c}
(-\eta + 1)/4 \\
0 
\end{array} \right), \quad 
\bv_{2}^1 =  \left( \begin{array}{c}
0 \\
 (\xi + 1)/4 
\end{array} \right) , \\
\bv_{3}^1 =&  \left( \begin{array}{c}
(\eta + 1)/4 \\
0
\end{array} \right) , \quad 
\bv_{4}^1 =  \left( \begin{array}{c}
0 \\
 (-\xi + 1)/4
\end{array}\right).
\end{aligned}
\end{equation}

\sssect{\hspace{-5mm}. Second-order quadrilateral element NQ2 \\} 
The  N\'ed\'elec space of this element reads 
 \begin{equation}
 \begin{aligned}
  \left[ \mathcal{ND}^\square \right]^{2}_2 = \bigg\{ & \left[ \begin{array}{c}
1 \\
0  
  \end{array} \right] \,, \left[ \begin{array}{c}
\xi \\
0
  \end{array} \right]  \,, \left[ \begin{array}{c}
\eta \\
0
  \end{array} \right]  \,, \left[ \begin{array}{c}
\xi \eta \\
0
  \end{array} \right]  \,,  \left[ \begin{array}{c}
\eta^2 \\
0
  \end{array} \right]  \,,    \left[ \begin{array}{c}
\xi \eta^2 \\
0
  \end{array} \right]  \,,  \\ & \left[ \begin{array}{c} 
0 \\
1 
  \end{array} \right] \,, \left[ \begin{array}{c}
0 \\
\xi
  \end{array} \right]  \,, \left[ \begin{array}{c}
0 \\
\eta
  \end{array} \right]  \,, \left[ \begin{array}{c}
0 \\ 
\xi \eta \\
  \end{array} \right]  \,,  \left[ \begin{array}{c}
0 \\
\xi^2
  \end{array} \right]  \,,  \left[ \begin{array}{c}
0 \\
\eta \xi^2
  \end{array} \right] \bigg\}  \,,
   \end{aligned}
 \end{equation}
and the vectorial shape functions have the following general form 
\begin{equation}
\bv^2 =  \left( \begin{array}{c}
a_1 + a_2 \,\xi + a_3\, \eta + a_4\, \xi \eta + a_5\, \eta^2 +  a_6 \, \xi \eta^2 \\ a_7 + a_8 \, \xi + a_9 \, \eta +  a_{10} \, \xi \eta + a_{11} \, \xi^2 + a_{12}  \,\eta \xi^2
\end{array} \right) \,,
\end{equation}
where $a_i, i=1,..,12$ are coefficients yet to be defined based on the dofs. Starting from Equations  (\ref{eq:t:edge_dof}) and (\ref{eq:t:inner_dof_q}), the explicit functions $r_j$ and $\bq_i$ are set as 
\begin{equation}
 \begin{aligned}
\textrm{edge 1:} \quad  r_1 =& \frac{1}{2}  (1- \xi)\,, \quad &r_2 =&\frac{1}{2}  (1+ \xi)\,, \\
\textrm{edge 2:} \quad r_1 =& \frac{1}{2}  (1- \eta)\,,  \quad &r_2 =&\frac{1}{2}  (1+ \eta)\,, \\
\textrm{edge 3:} \quad r_1 =& \frac{1}{2}  (1+ \xi)\,, \quad &r_2 =&\frac{1}{2}  (1- \xi)\,, \\
\textrm{edge 4:} \quad r_1 =& \frac{1}{2}  (1+ \eta)\,,\quad &r_2 =&\frac{1}{2}  (1- \eta)\,, \\
\textrm{inner :} \quad \bq_1 =& \left[ \begin{array}{c}
\frac{1}{2}  (1+ \xi) \\
0 
  \end{array} \right] , \quad &\bq_2 =&\left[ \begin{array}{c}
\frac{1}{2}  (1- \xi) \\
0 
  \end{array} \right] , \\ \bq_3 =& 
  \left[ \begin{array}{c}
0 \\
 \frac{1}{2}  (1+ \eta) 
  \end{array} \right] , \quad &\bq_4 =&\left[ \begin{array}{c}
0 \\
 \frac{1}{2}  (1- \eta) 
  \end{array} \right]
 \,,
  \end{aligned}
\end{equation}
The edge and inner dofs are calculated according to Equations  (\ref{eq:t:edge_dof}) and (\ref{eq:t:inner_dof_q})  considering that the tangential vector and the coordinates coloration are same as in NQ1 element
\begin{equation}
 \begin{aligned}
  m^{e_1}_{1} =& \frac{1}{3} \, (3a_1 - a_2 - 3a_3 + a_4 + 3a_5 - a_6) \,,  &m^{e_1}_{2} &= \frac{1}{3}  \, (3a_1 + a_2 - 3a_3 - a_4 + 3a_5 + a_6)\,, 
  \\  m^{e_2}_{1} =& \frac{1}{3}  \, (-a_{10} +3a_{11} - a_{12} + 3a_7 + 3a_8 - a_9) \,,   &m^{e_2}_{2} &= \frac{1}{3}  \, (a_{10} +3a_{11} + a_{12} + 3a_7 + 3a_8 + a_9) \,, 
  \\  m^{e_3}_{1} =& \frac{1}{3}  \,(3a_1 + a_2 + 3a_3 + a_4 + 3a_5 + a_6) \,,  &m^{e_3}_{2} &= \frac{1}{3} \,(3a_1 - a_2 + 3a_3 - a_4 + 3a_5 - a_6) \,, 
  \\  m^{e_4}_{1} =& \frac{1}{3} \,( -a_{10} + 3a_{11} + a_{12} + 3a_7 - 3a_8 + a_9) \,,  &m^{e_4}_{2} &= \frac{1}{3} \,( a_{10} + 3a_{11} - a_{12} + 3a_7 - 3a_8 - a_9) \,, 
  \\  m^\textrm{inner}_1 =& \frac{2}{9} (9 a_1 + 3 a_2 + 3 a_5 + a_6) \,, &m^\textrm{inner}_2 &= \frac{2}{9} (9 a_1 - 3 a_2 + 3 a_5 - a_6) \,, 
  \\ m^\textrm{inner}_3 =& \frac{2}{9} (3 a_{11} + a_{12} + 9 a_7 + 3 a_9) \,, &m^\textrm{inner}_4 &= \frac{2}{9} (3 a_{11} - a_{12} + 9 a_7 - 3 a_9) \,.
 \end{aligned}
\end{equation}

\begin{figure}[ht]
     \begin{subfigure}[b]{0.32\textwidth}
         \centering
         \includegraphics[width=\textwidth]{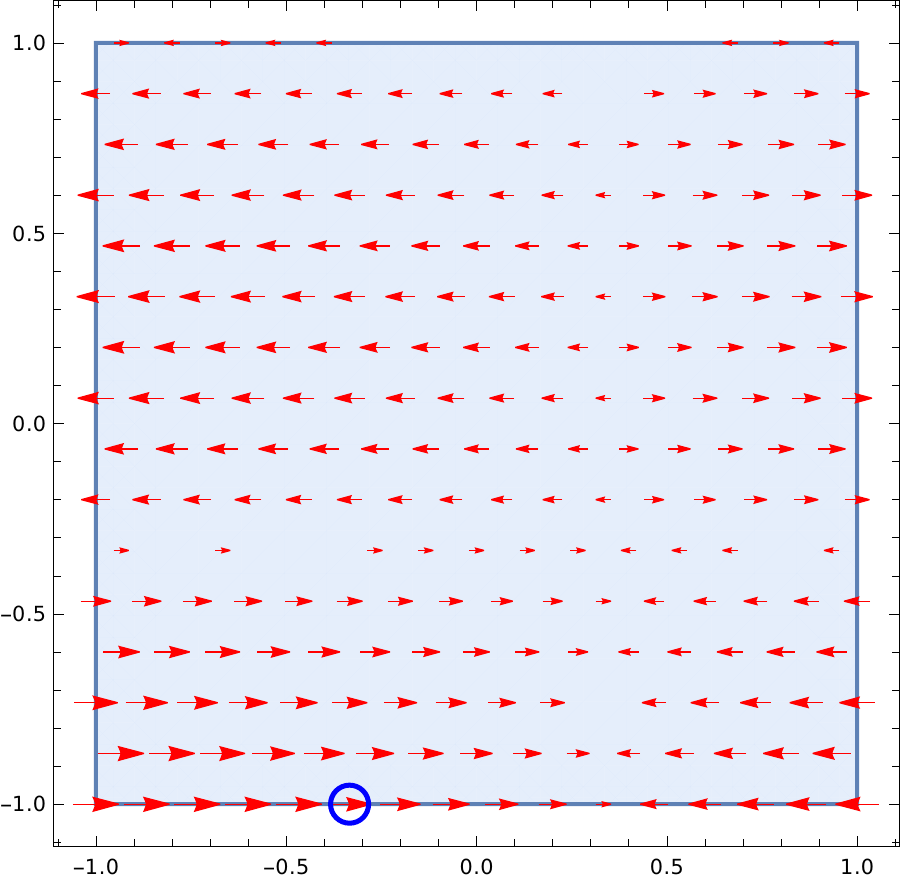}
        \caption{$\bv^2_1$}
     \end{subfigure}
     \begin{subfigure}[b]{0.32\textwidth}
         \centering
         \includegraphics[width=\textwidth]{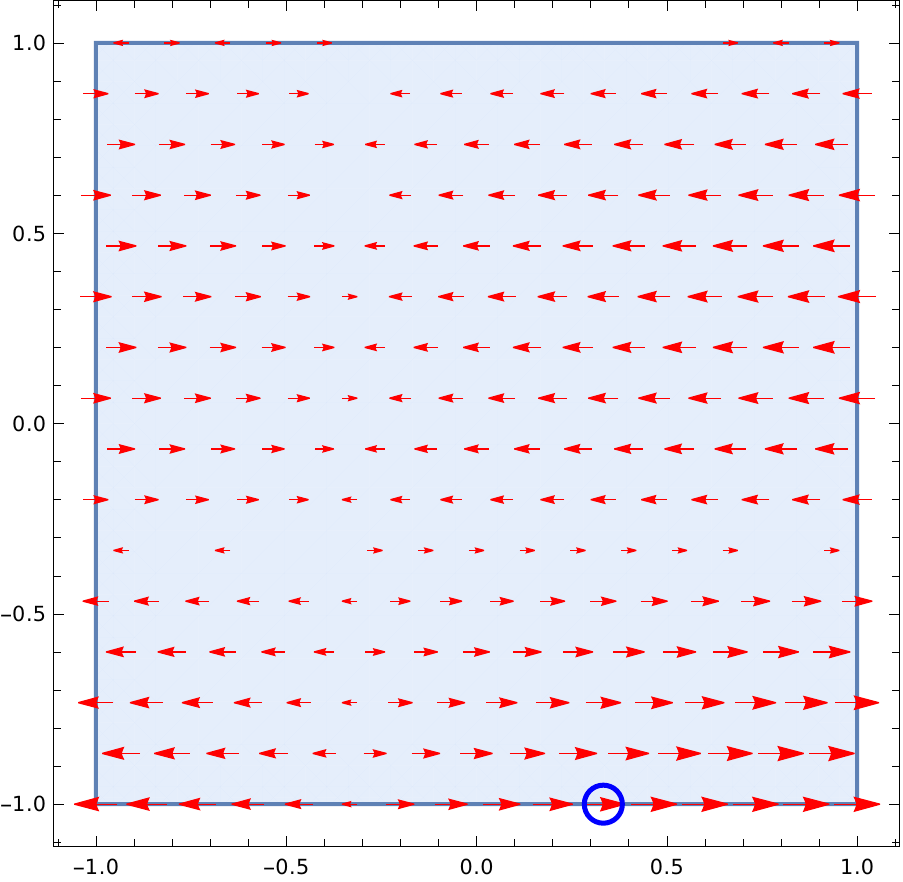}
  \caption{$\bv^2_2$}
     \end{subfigure}
          \begin{subfigure}[b]{0.32\textwidth}
         \centering
         \includegraphics[width=\textwidth]{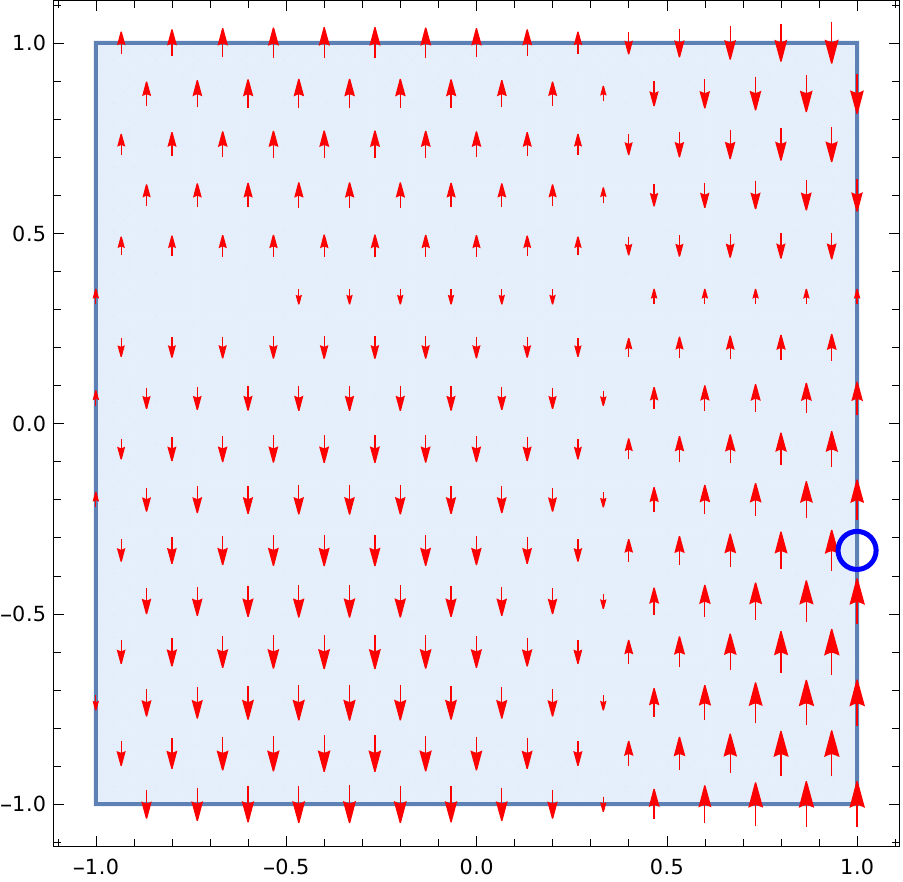}
  \caption{$\bv^2_3$}
     \end{subfigure}
          \begin{subfigure}[b]{0.32\textwidth}
         \centering
         \includegraphics[width=\textwidth]{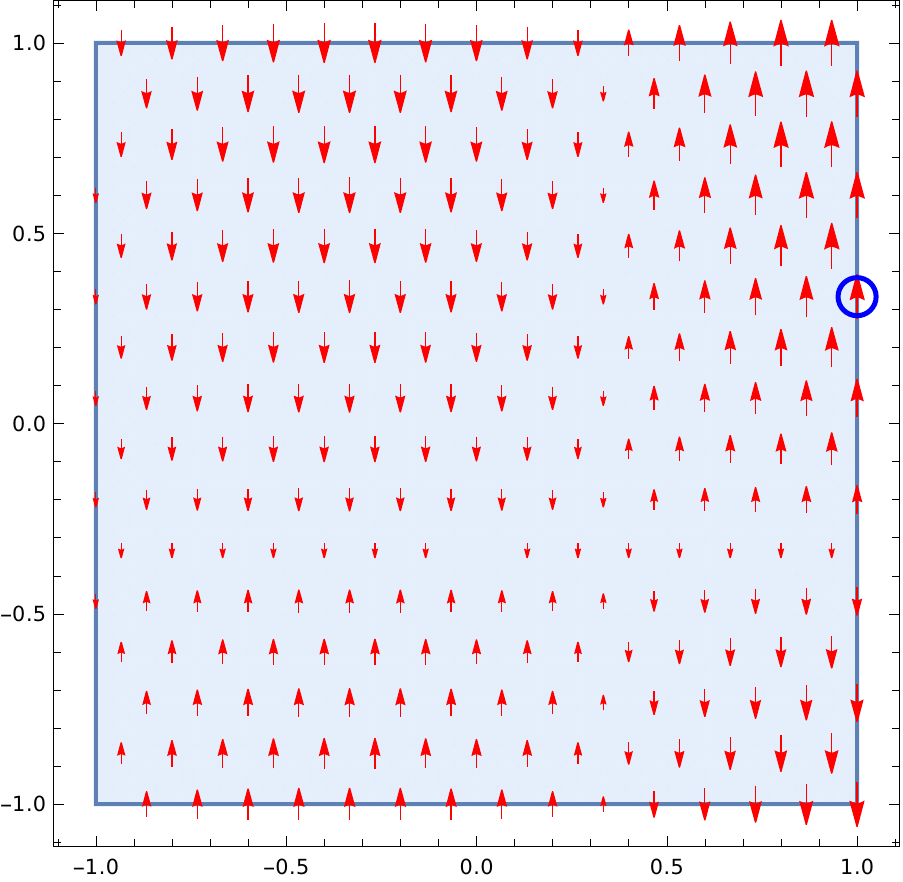}
        \caption{$\bv^2_4$}
     \end{subfigure}
       \begin{subfigure}[b]{0.32\textwidth}
         \centering
         \includegraphics[width=\textwidth]{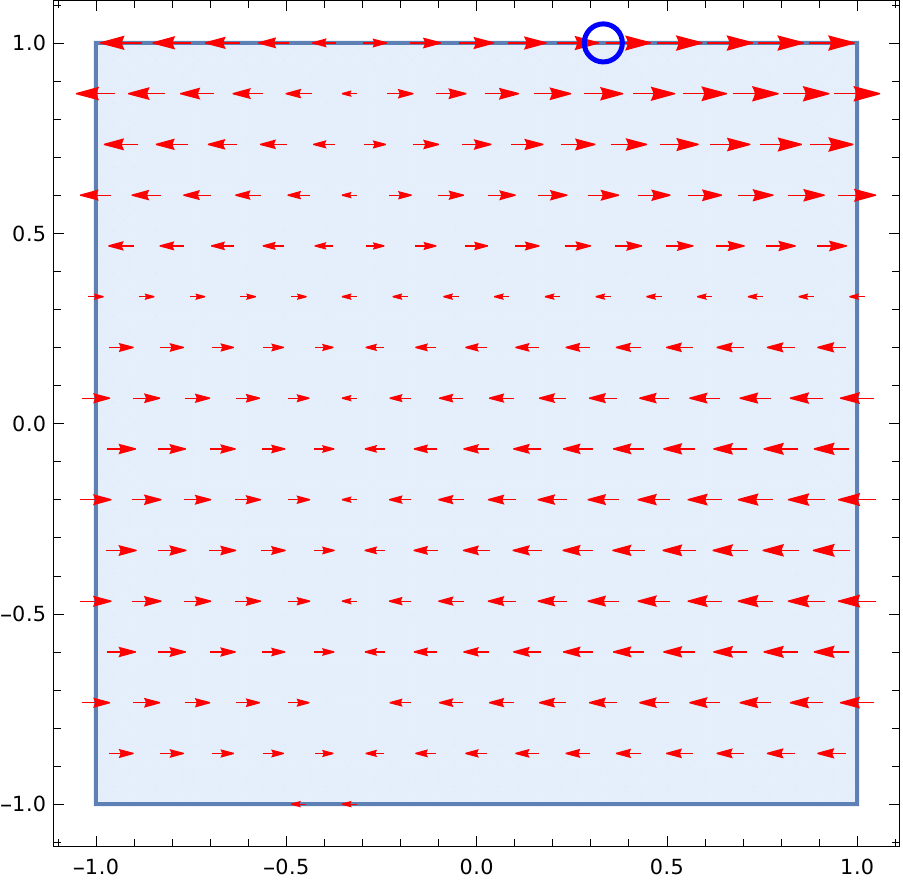}
        \caption{$\bv^2_5$}
     \end{subfigure}
     \begin{subfigure}[b]{0.32\textwidth}
         \centering
         \includegraphics[width=\textwidth]{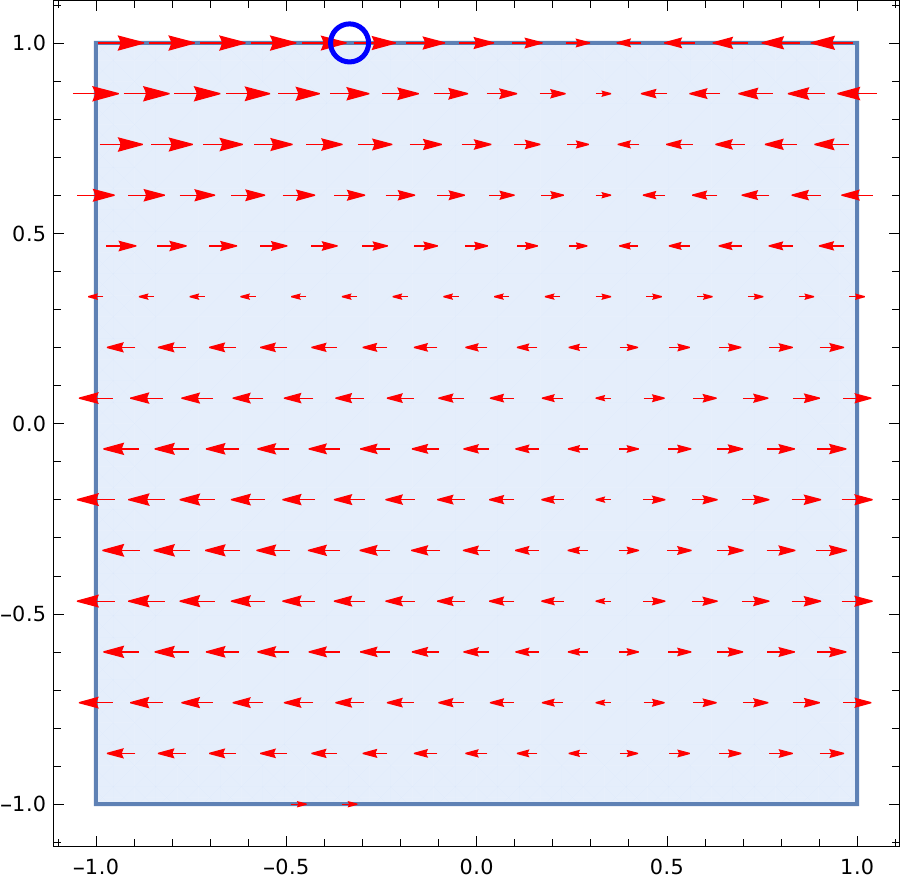}
  \caption{$\bv^2_6$}
     \end{subfigure}
          \begin{subfigure}[b]{0.32\textwidth}
         \centering
         \includegraphics[width=\textwidth]{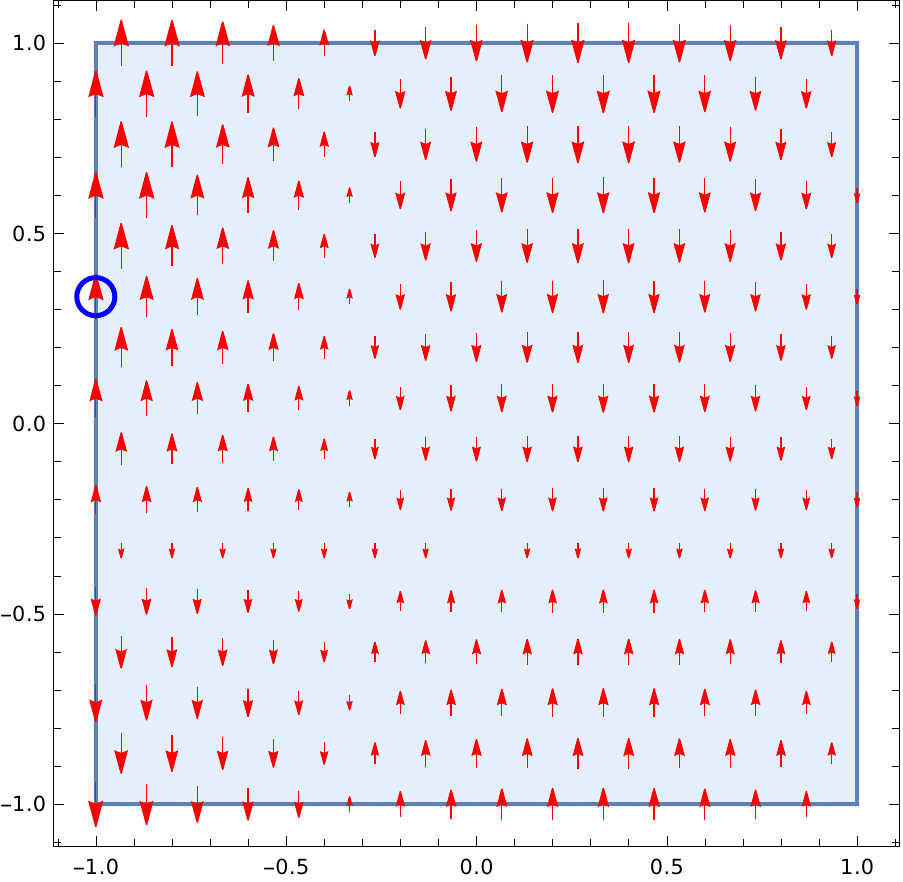}
  \caption{$\bv^2_7$}
     \end{subfigure}
          \begin{subfigure}[b]{0.32\textwidth}
         \centering
         \includegraphics[width=\textwidth]{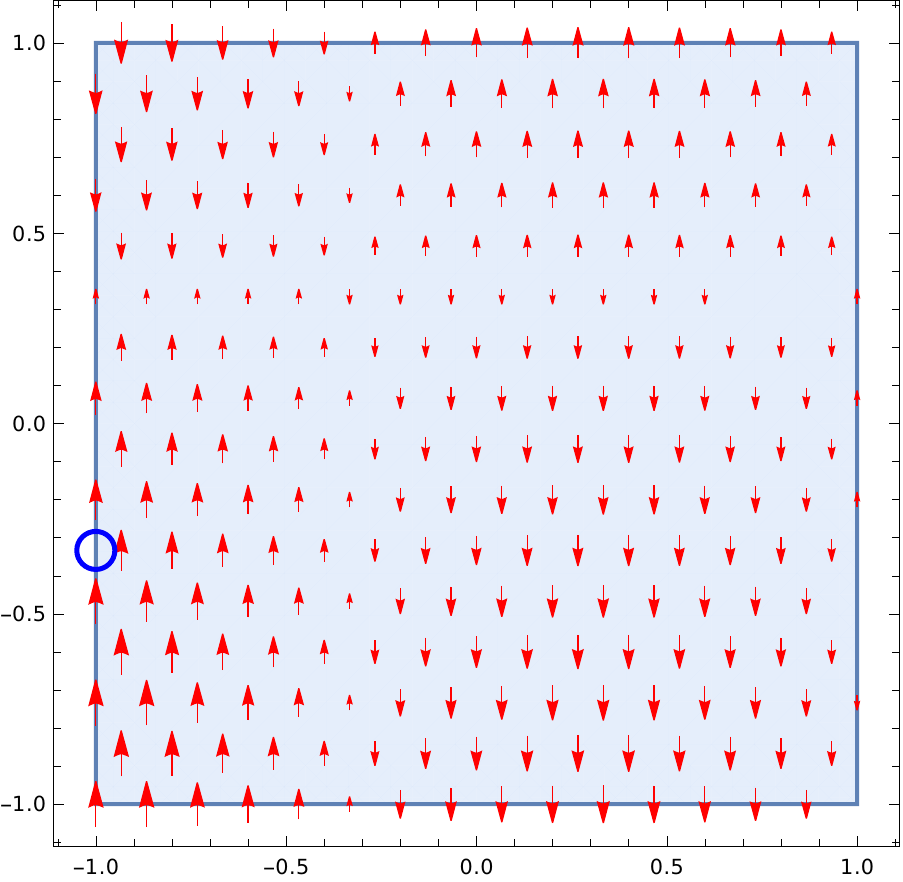}
        \caption{$\bv^2_8$}
     \end{subfigure}
        \begin{subfigure}[b]{0.32\textwidth}
         \centering
         \includegraphics[width=\textwidth]{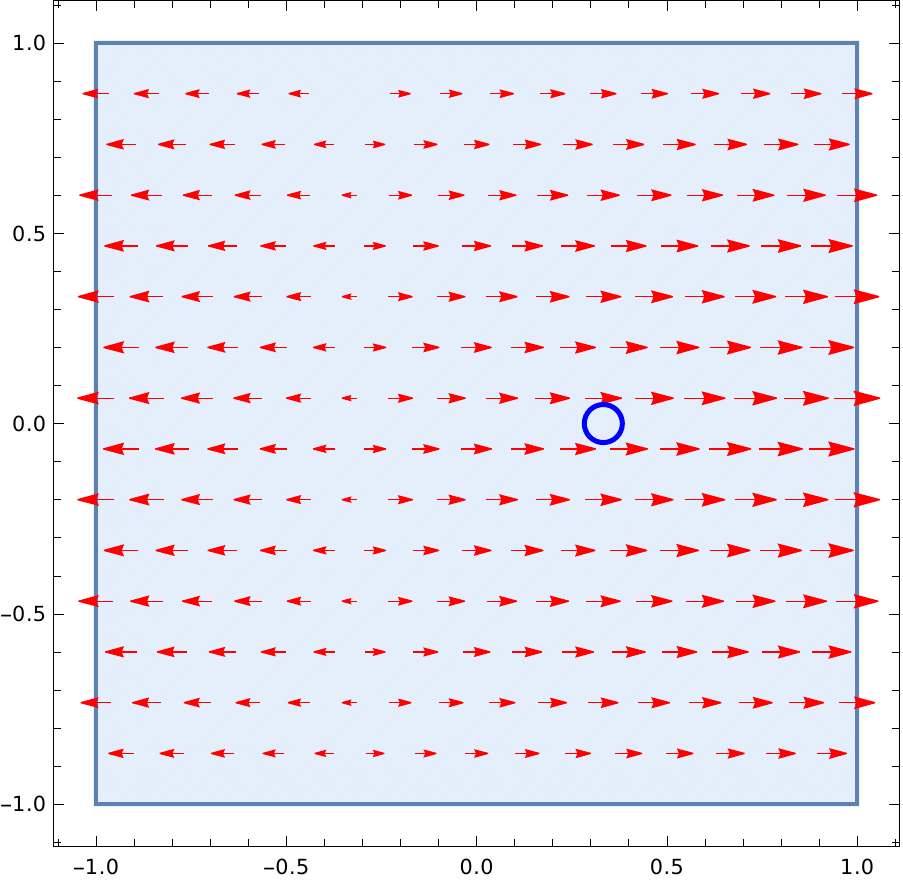}
        \caption{$\bv^2_9$}
     \end{subfigure}
     \begin{subfigure}[b]{0.32\textwidth}
         \centering
         \includegraphics[width=\textwidth]{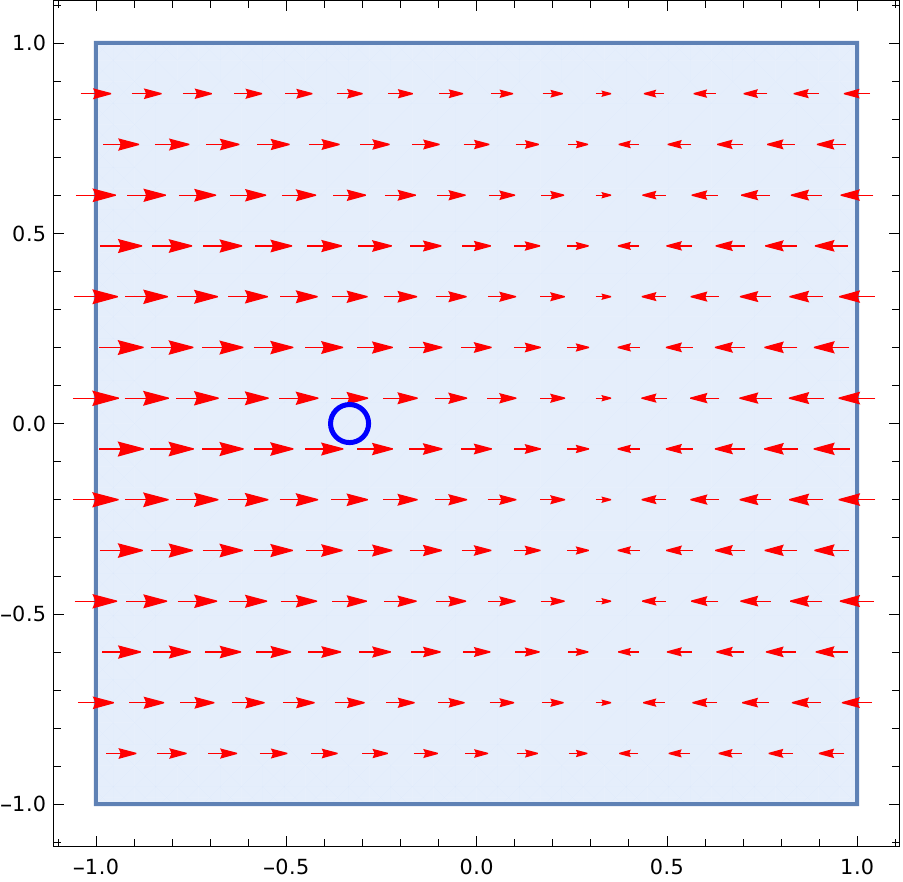}
  \caption{$\bv^2_{10}$}
     \end{subfigure}
          \begin{subfigure}[b]{0.32\textwidth}
         \centering
         \includegraphics[width=\textwidth]{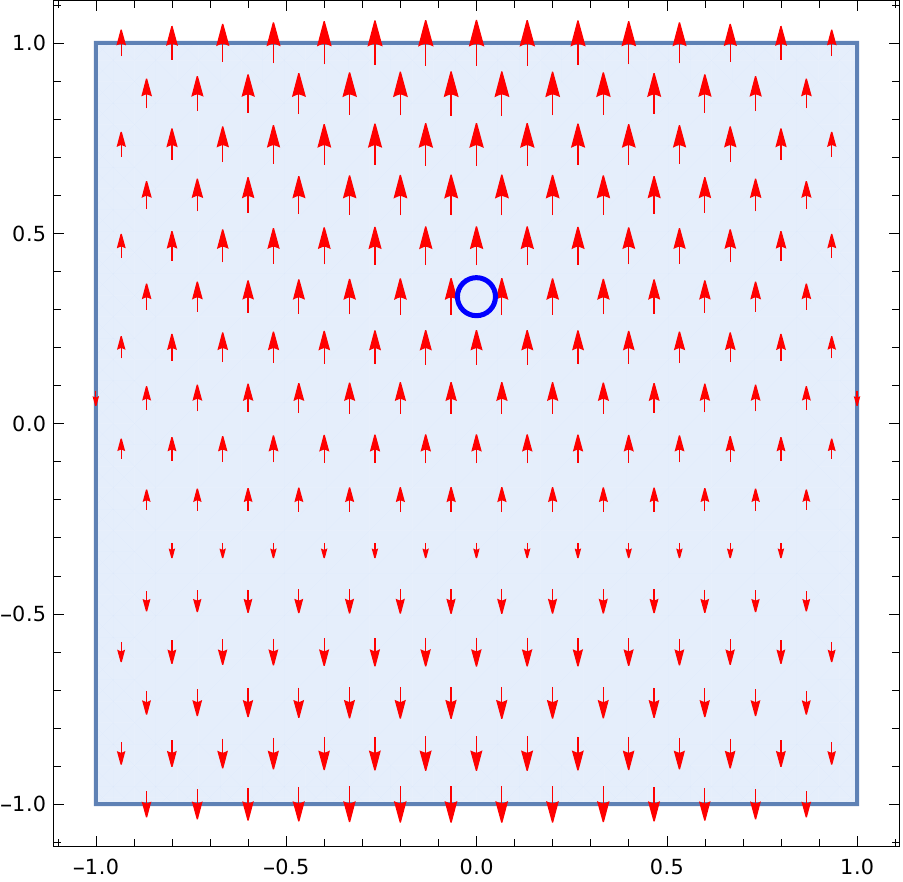}
  \caption{$\bv^2_{11}$}
     \end{subfigure}
          \begin{subfigure}[b]{0.32\textwidth}
         \centering
         \includegraphics[width=\textwidth]{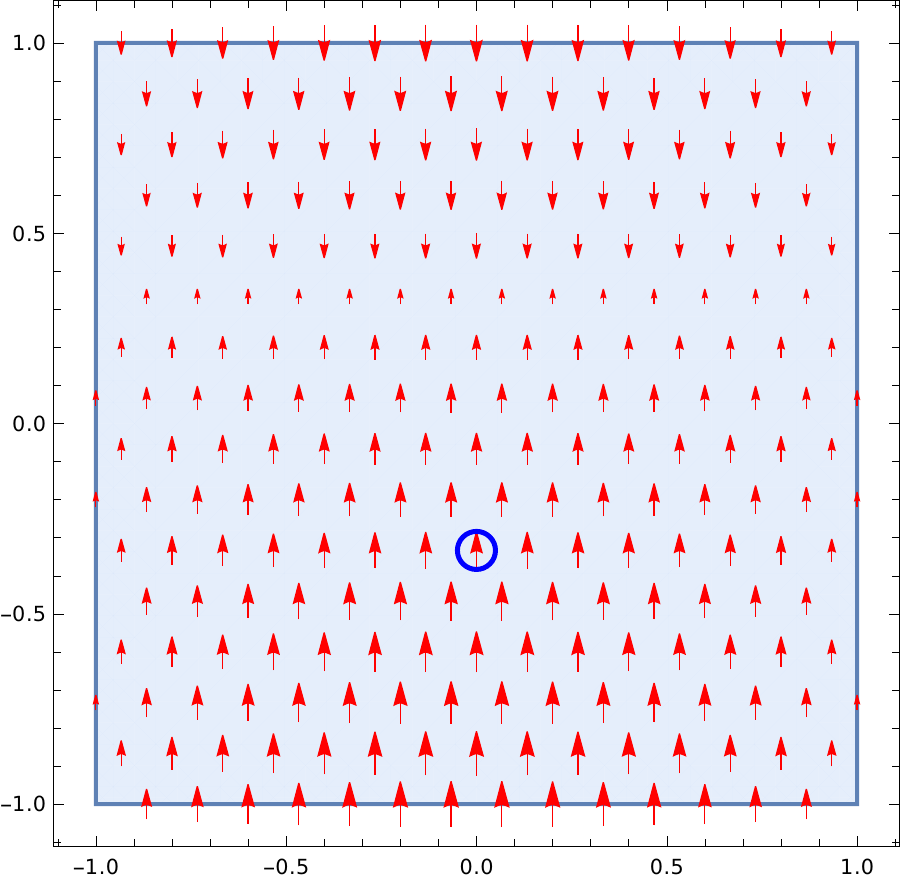}
        \caption{$\bv^2_{12}$}
     \end{subfigure}
       \caption{Tangential-conforming vectorial shape functions of NQ2 element.  Blue circles indicate the position where the dofs are defined. }
        \label{Fig:shape_function_NQ2}
\end{figure}

The  basis functions demonstrated in Figure  \ref{Fig:shape_function_NQ2}  are obtained by an analogous procedure as before

\begin{equation}
\begin{aligned}
\textrm{edge 1:} \quad \bv^2_1 &= \left( \begin{array}{c} 
-1/8 - \eta/4 + 3 \eta^2 /8 + 3 \xi/8 +  3 \eta \xi/4 - 9 \eta^2 \xi/8 \\
0
 \end{array} \right), \\
  \bv^2_2 &= \left( \begin{array}{c} 
-1/8 - \eta/4 + 3 \eta^2 /8 - 3 \xi/8 -  3 \eta \xi/4 + 9 \eta^2 \xi/8 \\
0
 \end{array} \right), \\
\textrm{edge 2:} \quad   \bv^2_3& = \left( \begin{array}{c} 
0 \\
- 1/8 + 3 \eta/8 + \xi/4 - 3 \eta \xi/4 +  3 \xi^2/8 - 9 \eta \xi^2/8
 \end{array} \right), \\
 \bv^2_4 &= \left( \begin{array}{c} 
0 \\
- 1/8 - 3 \eta/8 + \xi/4 + 3 \eta \xi/4 +  3 \xi^2/8 + 9 \eta \xi^2/8
 \end{array} \right), \\
\textrm{edge 3:} \quad   \bv^2_5 &= \left( \begin{array}{c} 
-1/8 + \eta/4 + 3 \eta^2/8 - 3 \xi/8 +  3 \eta \xi/4 + 9 \eta^2 \xi/8 \\
0
 \end{array} \right), \\
 \bv^2_6 &= \left( \begin{array}{c} 
-1/8+ \eta/4 + 3 \eta^2/8 +  3 \xi/8 - 3 \eta \xi/4 - 9 \eta^2 \xi/8 \\
0
 \end{array} \right), \\
\textrm{edge 4:} \quad   \bv^2_7 &= \left( \begin{array}{c} 
0 \\
-1/8 - 3 \eta/8 - \xi/4 - 3 \eta \xi/4 +  3 \xi^2/8 + 9 \eta \xi^2/8
 \end{array} \right), \\
  \bv^2_8 &= \left( \begin{array}{c} 
0 \\
-1/8 + 3 \eta/8 - \xi/4 + 3 \eta \xi/4 +  3 \xi^2/8 - 9 \eta \xi^2/8
 \end{array} \right)\,.
  \\
\textrm{inner:} \quad   \bv^2_9 &= \left( \begin{array}{c} 
3/8 - 3 \eta^2/8 + 9 \xi/8 - 9 \eta^2 \xi/8 \\
0
 \end{array} \right)\,, 
 \\ 
  \bv^2_{10} &= \left( \begin{array}{c} 
3/8 - 3 \eta^2/8 - 9 \xi/8 + 9 \eta^2 \xi/8 \\
0
 \end{array} \right)\,, 
 \\
   \bv^2_{11} &= \left( \begin{array}{c} 
0 \\
3/8 +9 \eta/8 - 3 \xi^2/8 - 9 \eta \xi^2/8
 \end{array} \right)\,, 
 \\
   \bv^2_{12} &= \left( \begin{array}{c} 
0 \\
3/8 -9 \eta/8 - 3 \xi^2/8 + 9 \eta \xi^2/8
 \end{array} \right)\,. 
 \end{aligned}
\end{equation}

\end{appendix}

\end{document}